\documentclass[11pt,psamsfonts,leqno,twoside]{amsart}
\usepackage{amssymb}

\usepackage{epic}
\usepackage{eepicemu}
\usepackage{diagrams}
\usepackage{graphicx}

\usepackage{xfm}

\newcommand\inlinegraphic[2][{scale=1.0}]{\begin{array}{c} 
\includegraphics[#1]{./EPS/#2}\end{array}}
\newcommand\inlinegraphicee[1]
{\begin{array}{c} 
\input{./EEPIC/#1.eepic}\end{array}}

\numberwithin{equation}{section}
\numberwithin{figure}{section}

\newtheorem{theorem}{Theorem}[section]
\newtheorem{proposition}[theorem]{Proposition}
\newtheorem{corollary}[theorem]{Corollary}
\newtheorem{lemma}[theorem]{Lemma}
\newtheorem{definition}[theorem]{Definition}

\newtheorem{remark}[theorem]{Remark}

\font\nnbf=plbx9 at9pt
\def\section#1{
\vskip-\lastskip\vskip18pt plus2pt minus2pt
\refstepcounter{section}
\centerline{{\uppercase{\nnbf\thesection. #1}}}
\vskip-\lastskip\vskip12pt plus2pt minus2pt}

\newcounter{ssection}
\def\ssection#1#2{
\vskip-\lastskip\vskip18pt plus2pt minus2pt
\refstepcounter{ssection}
\setcounter{theorem}{0}
\setcounter{subsection}{0}
\centerline{{\uppercase{\nnbf\thessection. #1}}}
\centerline{{\uppercase{\nnbf #2}}}
\vskip-\lastskip\vskip12pt plus2pt minus2pt}

\def\subsection#1{\vskip-\lastskip\vskip12pt 
plus2pt minus2pt
\refstepcounter{subsection}
{\indent\bf\thesubsection. #1}\hskip.55em\ignorespaces
}

\def\R{{\mathbb R}}
\def\Z{{\mathbb Z}}
\def\nat{{\mathbb N}}

\def\C{{\mathbb C}}
\def\la{\lambda}
\def\eps{\varepsilon}
 
\def\S{{\mathfrak{S}}}
\def\zdeltahat{\Z[\deltabold^{\pm 1}]^\wedge}

\def\u #1 #2{\mathcal U(#1, #2)}
\def\uhat #1 #2{\widehat{\mathcal U}(#1, #2)}
\def\varLambdaHat{{\widehat{\varLambda}}}
\def\RHat{{\widehat R}}
\def\kt #1{{{\rm KT}_{#1}}}
\def\akt #1{\widehat{{\rm KT}}_{#1}}
\def\abmw #1{\widehat{W}_{#1}}
\def\bmw #1{W_{#1}}
\def\w #1 #2{\bmw {#1}^{(#2)}}
\def\V #1 #2{V_{#1}^{(#2)}}
\def\k #1 #2{ KT_{#1}^{(#2)}}
\def\br{D}
\def\abr{\widehat{D}}
\def\p #1{\mathbold {#1}}
\def\pbar #1{\overline{\p #1}}

\def\Zdeltabold{\Z[\deltabold^{\pm 1}]}
\def\perm{{\rm perm}}
\def\bdry{\partial}
\def\pmone{^{\pm 1}}
\def\inv{^{-1}}

\def\mathbold{\pmb}
\def\lambdabold{{\mathbold \la}} 
\def\labold{\lambdabold}
\def\deltabold{{\mathbold \delta}}
\def\zbold{{\mathbold z}}
\def\qbold{{\mathbold q}}

\begin{document}


\inicaut1{F. M.}
\imaut1{Frederick M.}
\nazaut1{Goodman}
\adraut1{Department of Mathematics\cr
University of Iowa\cr
Iowa City, Iowa, U.S.A.}
\emailaut1{goodman@math.uiowa.edu}
\inicaut2{H.}
\imaut2{Holly}
\nazaut2{Hauschild}
\adraut2{Department of Mathematics\cr
 University of California San Diego\cr
La Jolla, California, U.S.A.}
\emailaut2{hhauschild@math.ucsd.edu}
\abbrevtitle{Affine BWM algebras}
\dedykacja{Dedicated to William G. Bade on his 80th birthday}
\title{Affine Birman--Wenzl--Murakami  algebras and\cr tangles in
the solid torus}
\tyt{\im1\ \naz1\ {{\rm (Iowa City, IA)\ and}}\cr \im2\ \naz2}{La Jolla, CA}
\ab{The affine Birman--Wenzl--Murakami algebras can be defined algebraically, 
via generators and relations, or geometrically as algebras of tangles in the 
solid torus, modulo Kauffman skein relations.  We prove that the two versions 
are isomorphic, and we show that these algebras are free over any ground ring, 
with a basis similar to a well known basis of the affine Hecke algebra.}


\def\hhh#1#2{{\noindent\footnotesize#1}\hfill{\footnotesize#2}\vskip-2pt}
\hhh{{\bf Contents}}{}\vskip4pt
\hhh{1. Introduction}{\pageref{intro}}
\hhh{2. The affine Kauffman tangle algebra}{\pageref{sec2}}
\hhh{3. The affine Birman--Wenzl--Murakami algebra}{\pageref{sec3}}
\hhh{4. The affine Brauer algebra}{\pageref{sec4}}
\hhh{5. Isomorphism of ordinary BMW and Kauffman tangle algebras}{\pageref{sec5}}
\hhh{6. Isomorphism of affine BMW and Kauffman tangle algebras}{\pageref{sec6}}
\hhh{References}{\pageref{biblio}}
\vskip18pt plus2pt minus2pt

\normalsize

\section{Introduction}\label{intro}
\setcounter{equation}{0}
\setcounter{figure}{0}
The purpose of this paper is to establish an isomorphism between the affine 
  Birman--Wenzl--Murakami algebras (defined by
generators and relations) and the algebras of $(n,n)$-tangles in the
solid torus, modulo Kauffman skein relations.

The Birman--Wenzl--Murakami  (or BMW)  algebras were conceived 
in~\cite{Birman-Wenzl} 
and~\cite{Murakami-BMW} as an algebraic framework for the Kauffman link 
invariant 
of~\cite{Kauffman}.  The definition of these algebras by generators and 
relations was
motivated by  certain canonical elements, and relations satisfied by these 
elements, in
Kauffman tangle algebras in $D^2 \times I$, that is,  algebras of framed 
$(n,n)$-tangles
in $D^2 \times I$, modulo Kauffman skein relations.  
 It is therefore not surprising, and is sometimes taken as evident,  that the 
BMW algebras
are isomorphic
 to the Kauffman tangle algebras.   A proof of the isomorphism was given by 
Morton and
Wassermann~\cite{Morton-Wassermann} in a paper from 1989 that unfortunately was 
never
published. (A related  result is given by Kauffman~\cite[Theorem~4.4]{Kauffman}.
  See Remark~\ref{remark-on the definition of BMW} below.)  Some
aspects of the BMW algebras can be understood more clearly in the Kauffman 
tangle picture;
for example, the existence of the Markov  trace and conditional expectations is 
evident
from the tangle picture,
 and an argument of Morton and Traczyk~\cite{Morton-Traczyk}  gives the freeness 
of the
Kauffman tangle algebras over an arbitrary ring.
  
The affine BMW algebras are related to the ordinary BMW algebras
as the affine Hecke algebras of type $A$ are related to  the
ordinary Hecke algebras of type $A$.  The affine algebras have
an ``extra'' generator $x_1$ which satisfies the braid relation
$x_1 g_1 x_1 g_1 =  g_1 x_1 g_1 x_1$ with the first ``ordinary''
braid generator~$g_1$.  The extra generator of the affine
algebras  can be imagined geometrically as a strand looping
around the hole in $A \times I$, where $A$ denotes the annulus
$S^1
\times I$.
The full set of relations for the affine BMW algebras (due to
H\"aring-Oldenburg~\cite{H-O2}) 
are modeled after relations which hold in the Kauffman tangle
algebras in $A \times I$.  So the affine BMW algebras ``should"
be  isomorphic to the Kauffman tangle algebras in $A \times I$,
assuming that a sufficient list of relations has been
discovered.  This isomorphism is our main result.  On the way to
proving the isomorphism, we also show that the affine BMW
algebras, over any ring,  are free with a basis generalizing a
well known basis of the affine Hecke algebras.
 
As for the ordinary BMW algebras, certain properties of the
affine BMW algebras are most easily obtained from the Kauffman
tangle picture---in particular the existence of the Markov trace
and conditional expectations.  Moreover,  we require the tangle
picture even to describe our basis of the affine BMW algebras.
  
In proving our main theorems, we rely on the results of Morton
and Wassermann for the ordinary BMW and Kauffman tangle
algebras, as well as on techniques from their paper.  Since
their paper is not generally available, we have reproduced their
results in Section~5 of the present paper.
 
The outline of the paper is as follows.  In Sections~2 and 3 we
introduce the ordinary and affine Birman--Wenzl--Murakami
algebras and the Kauffman tangle algebra and derive basic
properties of these algebras.  In Section~4, we introduce an
affine version of the Brauer centralizer algebra, which is a
homomorphic image of the  Kauffman tangle algebra in $A \times
I$.  Section~5 is devoted to the Morton--Wassermann proof of the
isomorphism of the ordinary BMW and Kauffman tangle algebras.
Finally, in Section~6 we obtain the main results of the paper.
  
In a future paper, we  intend to study the  cyclotomic BMW
algebras, which are quotients of the affine BMW algebras in
which the affine generator $x_1$ satisfies a polynomial
relation.  Cyclotomic BMW algebras have previously appeared
in~\cite{H-O1}, \cite{H-O2},~\cite{RO}.
 
Let us mention some additional antecedents for the algebras
studied here.  The affine and cyclotomic Brauer algebras, which
we use in Section~4, were introduced by H\"aring-Oldenburg
\cite{H-O2}, and the cyclotomic case was studied
in~\cite{Rui-Yu}. Formally similar algebras determined by
weighted Brauer diagrams are used in the approach of Olshanski
and Okounkov to the representation theory of the infinite
symmetric group; see, for
example,~\cite{Olshanski},~\cite{Okounkov}.  Nazarov introduced
a ``degenerate'' affine BMW algebra, under the name ``degenerate
affine Wenzl algebra'', in~\cite{Nazarov}.  The cyclotomic
version of this algebra was recently studied
in~\cite{Ariki-Mathas-Rui}.   Ram and Orellana~\cite{RO} have studied certain representations of the affine BMW algebras.

We would like to thank Hans Wenzl for finding a significant error in a previous version of this paper.

\section{The affine Kauffman tangle algebra}\label{sec2}
\setcounter{equation}{0}
\setcounter{figure}{0}
In this section we introduce the geometric versions of the ordinary and
affine Birman--Wenzl--Murakami  algebras.

\subsection{Framed tangles.}
\hskip-2pt The objects considered in this subsection, namely  framed tangles in 
$S \times I$,
where $S$ is an oriented surface, do not enter directly into the definition of 
the
Kauffman tangle algebra and its affine counterpart, but they provide
important motivation for the definition.

 Let $S$ be  a smooth oriented surface (with boundary).  Let
$I$ denote the unit interval $[0,1]$.   Choose once and for all a
countable family of mutually disjoint oriented intervals
$\{J_i : i \in \nat\}$ in $S$.

\newcounter{blicz}
\renewenvironment{enumerate}{\begin{list}{\rm(\arabic{blicz})}
{\usecounter{blicz}%
\def\makelabel##1{\hss\llap{##1}}%
\setlength{\topsep}{4pt}%
\setlength{\parsep}{0pt}%
\setlength{\itemsep}{0pt}%
\setbox0=\hbox{(1)\enspace}%
\setlength{\leftmargin}{\parindent}%
\addtolength{\leftmargin}{\wd0}
}}
{\end{list}}

\begin{definition}\rm
Fix integers $k, n \ge 0$.
A  {\em framed $(k,n)$-tangle}\/ in 
$S \times I$  is a  family of  non-intersecting embedded piecewise smooth
rectangles (ribbons)  in $S \times I$ such that: 

\begin{enumerate}
\item with $F$ denoting the union of the  ribbons,
\begin{eqnarray*}
\qquad F \cap \bdry(S\times I) &=& F \cap (S \times \{0, 1\}) \\ 
 &=& J_1\times \{0, 1\} \cup \dots \cup
J_k\times \{0, 1\}  ,
 \end{eqnarray*}
\item any ribbon in the family intersects the boundary of $S \times I$
transversally.
\end{enumerate}
\end{definition}

\begin{remark}\rm
Note that $k = 0$, or $n = 0$, or both, are allowed in the
definition of tangles.  If one of $k$ or $n$ is positive, then a $(k,
n)$-tangle must be a  non-empty family  of ribbons, but the empty tangle,
with no ribbons, is an allowable $(0,0)$-tangle.
\end{remark}

We will identify  framed tangles which are isotopic via an isotopy fixing the
intersections of the curves with the boundary of $S \times I$; so
strictly speaking, a tangle is an isotopy class of families of ribbons.

{\spaceskip 0.33em plus0.2em minus0.17em One can compose $(k,n)$-tangles and 
$(n, m)$-tangles by ``stacking''.
Namely, the first  tangle (the $(k,n)$-tangle) is placed above the second
tangle (the $(n,m)$-tangle),  the $n$ oriented intervals at which the ribbons
intersect  the lower boundary of the first tangle are identified with the
$n$ oriented intervals at which the ribbons
intersect  the upper boundary of the second tangle, and the resulting
family of ribbons is compressed into $S \times I$.}  

Since composition of tangles is evidently associative, framed tangles in $S
\times I$ may be regarded as morphisms in a certain category. Namely,
the  objects of the category are $\{0, 1, \dots\}$, and the morphisms
from $k$ to $n$ are the $(k,n)$-tangles in $S \times I$. 
In
particular, the set of $(n,n)$-tangles forms a monoid,  whose identity is
the tangle $J_1 \times I \cup \dots \cup J_n \times I $.

We are  interested here in two possibilities for $S$, namely the
plane $\R^2$ and  the annulus $A^2 = \{(x,y) \in \R^2 : 1/4 \le
x^2 + y^2 \}$.  We  refer to tangles in $\R^2
\times I$ as ``ordinary"  framed tangles and tangles in 
$A^2 \times I$ as ``affine" framed tangles.  For both cases we
take our intervals $J_i$ to lie on the positive $x$-axis and to
be ordered by their order on the $x$-axis.  For convenience, we
can take $J_i = [i, i + 1/2] \times \{0\}$.

\subsection{Tangle diagrams.}
Fix points $a_i$ in $\R$, for $i \ge 0$, with $0 = a_0 < a_1 < a_2 \cdots$.

\begin{definition}\rm
An {\em ordinary $(k,n)$-tangle diagram}  is a piece of a
knot diagram in $R = \R \times I$ such that:
\begin{enumerate}
\item with $F$ denoting the union of the curves comprising the knot
diagram, 
\begin{eqnarray*}
\qquad\ F \cap \bdry(R) &=& F \cap (I \times \{0, 1\}) \\
&=& \{(a_0,1), \dots, (a_{k-1},1)\} \cup \{ (a_0,0), \dots, (a_{n-1},0)\},
\end{eqnarray*}
\item  any curve in the family intersects the boundary of $R$
transversally.
\end{enumerate}
\end{definition}

Recall that a knot diagram means a collection of piecewise smooth curves
which may have intersections and self-intersections, but only simple
transverse intersections.  At each intersection or crossing, one of the
two strands (pieces of curves) which intersect is indicated as crossing
over the other.  

\begin{definition}\rm
  Two ordinary tangle diagrams are said to be {\em
ambient isotopic} if they are related by a sequence of {\em Reidemeister
moves} of types I, II, and III, followed by an isotopy of $R$ fixing the
boundary.  Two ordinary tangle diagrams are said to be {\em regularly
isotopic} if they are related by a sequence of Reidemeister moves of
types~II and III only, followed by an isotopy of $R$ fixing the
boundary.  See the following figure for the Reidemeister moves of types I, II, 
and~III.
\begin{eqnarray*}
\text{I}&\quad\ &\inlinegraphic{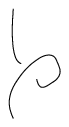} \quad \longleftrightarrow
 \quad \inlinegraphic{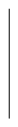} \quad
  \longleftrightarrow  \quad \inlinegraphic{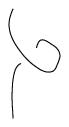}\\
\text{II}&\quad\ &\inlinegraphic[scale =.5] {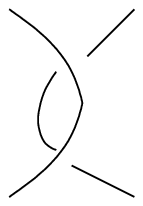} \quad 
\longleftrightarrow 
\quad \inlinegraphic[scale=1.75]{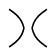} \\
  \text{III}&\quad\ &\inlinegraphic{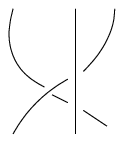} \quad \longleftrightarrow \quad 
\inlinegraphic{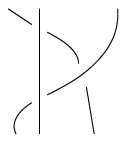} 
\end{eqnarray*} 
\end{definition}

Tangle diagrams can be composed, similarly to tangles, and composition
respects ambient isotopy or regular isotopy.  Thus one obtains a product
on ambient isotopy classes of tangle diagrams, respectively on regular
isotopy classes of tangle diagrams. 

As is well known, isotopy classes of ordinary framed tangles correspond
bijectively to regular isotopy classes of ordinary tangle diagrams,
by a framed version of Reidemeister's theorem.  The correspondence is as 
follows:  
A representative of an isotopy class of framed  tangles (in $(x, y, z)$-space)
can always be chosen that lies close to and almost parallel with the $xz$ plane; 
twists
in ribbons are converted to ``kinks" (or almost planar loops) in such a 
representative:

$$
\inlinegraphic[scale=.3]{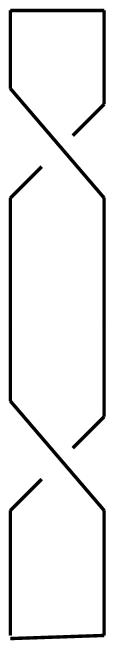} \quad \longleftrightarrow \quad 
\inlinegraphic{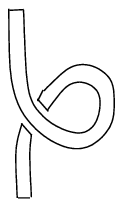}
$$

Replace each ribbon in such a representative by its ``core", a curve running 
lengthwise along
the center of the ribbon, and project the collection of these curves in the $xz$ 
plane. The
result is a tangle diagram.  Any two such projections are  regularly  isotopic 
tangle diagrams. 
Thus one has a well defined map from isotopy classes of ordinary framed tangles 
to regular
isotopy classes of tangle diagrams.  This map is bijective.

It follows that 
the monoid of (isotopy classes of) ordinary
$(n,n)$-framed tangles is isomorphic to the monoid of regular isotopy classes
of $(n,n)$-tangle diagrams.

\begin{definition}\rm\label{definition affine tangle}
An {\em affine  $(k,n)$-tangle diagram}\/ is an
ordinary
$(k+1,\break n+1)$-tangle diagram which includes a distinguished curve $P$
connecting $(a_0, 1)$ and $(a_0, 0)$ such that the height coordinate (the
second coordinate) varies monotonically along the curve $P$.
\end{definition}

Note that an affine tangle diagram is always equivalent by an isotopy of
$R$ fixing the boundary to a diagram including the curve $\{a_0\} \times
I$.  We will draw affine tangle diagrams with the distinguished curve
drawn as a thickened vertical segment.   We refer to the distinguished
curve as the ``flagpole''.  The figure below shows
some fundamental affine tangle diagrams from which more complex affine
tangle diagrams can be built by composition.\vadjust{\vskip6pt}
$$
\begin{array}{c}
\setlength{\unitlength}{0.00062500in}
\begingroup\makeatletter\ifx\SetFigFont\undefined%
\gdef\SetFigFont#1#2#3#4#5{%
  \reset@font\fontsize{#1}{#2pt}%
  \fontfamily{#3}\fontseries{#4}\fontshape{#5}%
  \selectfont}%
\fi\endgroup%
{\renewcommand{\dashlinestretch}{30}
\begin{picture}(659,1281)(0,-10)
\thicklines
\path(197,783)(197,33)
\path(197,1233)(197,933)
\thinlines
\path(647,1233)(647,1232)(647,1229)
	(647,1224)(646,1216)(646,1206)
	(645,1193)(643,1178)(641,1162)
	(638,1144)(633,1125)(628,1106)
	(621,1087)(612,1067)(601,1046)
	(587,1024)(570,1002)(550,979)
	(525,956)(497,933)(466,912)
	(434,893)(404,878)(375,867)
	(350,859)(327,855)(306,854)
	(287,854)(270,856)(253,858)
	(237,860)(220,862)(202,862)
	(183,861)(162,857)(140,849)
	(116,838)(92,823)(68,804)
	(47,783)(29,757)(18,732)
	(13,706)(12,682)(15,659)
	(22,636)(31,614)(42,593)
	(54,572)(67,553)(80,535)
	(92,519)(103,506)(111,496)
	(117,489)(120,485)(122,483)
\path(272,408)(274,408)(278,409)
	(286,409)(296,410)(310,411)
	(325,411)(342,410)(360,409)
	(379,406)(398,401)(417,395)
	(437,385)(457,372)(477,354)
	(497,333)(514,308)(529,282)
	(540,256)(549,230)(556,204)
	(561,179)(565,154)(567,129)
	(569,106)(570,85)(571,66)
	(572,52)(572,42)(572,36)(572,33)
\end{picture}
} \\[-6pt] \vrule height24pt width0pt depth0pt 
{\scriptstyle X_1} \end{array} 
\qquad
\begin{array}{c} \setlength{\unitlength}{0.00062500in}
\begingroup\makeatletter\ifx\SetFigFont\undefined%
\gdef\SetFigFont#1#2#3#4#5{%
  \reset@font\fontsize{#1}{#2pt}%
  \fontfamily{#3}\fontseries{#4}\fontshape{#5}%
  \selectfont}%
\fi\endgroup%
{\renewcommand{\dashlinestretch}{30}
\begin{picture}(1272,1546)(0,-10)
\put(903,40){\makebox(0,0)[lb]{{\SetFigFont{8}{9.6}{\rmdefault}{\mddefault}{\updefault}$i+1$}}}
\path(333,1498)(333,298)
\path(483,1498)(483,298)
\path(933,1498)(633,298)
\path(633,1498)(708,973)
\path(858,748)(933,298)
\path(1083,1498)(1083,298)
\path(1233,1498)(1233,298)
\put(558,73){\makebox(0,0)[lb]{{\SetFigFont{9}{10.8}{\rmdefault}{\mddefault}{\updefault}$i$}}}
\thicklines
\path(33,298)(33,598)(33,748)(33,1498)
\end{picture}
} \\[-6pt]  {\scriptstyle G_i} 
\end{array}
\qquad
\begin{array}{c} \setlength{\unitlength}{0.00062500in}
\begingroup\makeatletter\ifx\SetFigFont\undefined%
\gdef\SetFigFont#1#2#3#4#5{%
  \reset@font\fontsize{#1}{#2pt}%
  \fontfamily{#3}\fontseries{#4}\fontshape{#5}%
  \selectfont}%
\fi\endgroup%
{\renewcommand{\dashlinestretch}{30}
\begin{picture}(1245,1588)(0,-10)
\put(483,65){\makebox(0,0)[lb]{{\SetFigFont{8}{9.6}{\rmdefault}{\mddefault}{\updefault}$i$}}}
\path(333,340)(333,1540)
\path(483,340)(483,1540)
\path(1083,340)(1083,1540)
\path(1233,340)(1233,1540)
\path(558,340)(558,341)(559,344)
	(562,353)(566,367)(572,387)
	(579,411)(587,437)(595,462)
	(603,486)(611,508)(619,527)
	(627,544)(634,558)(642,570)
	(650,581)(658,590)(667,598)
	(676,605)(687,611)(697,617)
	(709,621)(721,624)(733,626)
	(746,626)(758,626)(770,624)
	(782,621)(794,617)(804,611)
	(815,605)(824,598)(833,590)
	(841,581)(849,570)(857,558)
	(864,544)(872,527)(880,508)
	(888,486)(896,462)(904,437)
	(912,411)(919,387)(925,367)
	(929,353)(932,344)(933,341)(933,340)
\path(558,1540)(558,1539)(559,1536)
	(562,1527)(566,1513)(572,1493)
	(579,1469)(587,1443)(595,1418)
	(603,1394)(611,1372)(619,1353)
	(627,1336)(634,1322)(642,1310)
	(650,1299)(658,1290)(667,1282)
	(676,1275)(687,1269)(697,1263)
	(709,1259)(721,1256)(733,1254)
	(746,1254)(758,1254)(770,1256)
	(782,1259)(794,1263)(804,1269)
	(815,1275)(824,1282)(833,1290)
	(841,1299)(849,1310)(857,1322)
	(864,1336)(872,1353)(880,1372)
	(888,1394)(896,1418)(904,1443)
	(912,1469)(919,1493)(925,1513)
	(929,1527)(932,1536)(933,1539)(933,1540)
\put(783,40){\makebox(0,0)[lb]{{\SetFigFont{8}{9.6}{\rmdefault}{\mddefault}{\updefault}$i+1$}}}
\thicklines
\path(33,1540)(33,1240)(33,1090)(33,340)
\end{picture}
} \\[-6pt]  {\scriptstyle  E_i} 
\end{array}
$$

Two affine $(k,n)$-tangle diagrams are ambient (respectively  regularly)
isotopic if they are ambient (resp. regularly) isotopic as ordinary $(k+1,
n+1)$-tangle diagrams.  Note that in an ambient isotopy, no Reidemeister
move of type I is ever applied to the flagpole, as the flagpole is
required to be represented by a monotonic path.

The set of affine tangle diagrams, regarded as a subset of ordinary
tangle diagrams,  is evidently closed under the composition of tangle
diagrams.  Thus one obtains a product
on ambient isotopy classes of affine tangle diagrams, respectively on
regular isotopy classes of affine tangle diagrams. 

We can extend the correspondence between isotopy classes of ordinary
framed tangles and regular isotopy classes of ordinary tangle diagrams to a
correspondence between isotopy classes of affine framed tangles and regular
isotopy classes of affine tangle diagrams, as follows.
First for each affine $(k,n)$-framed tangle $T$ in $A^2 \times I$, we can
consider  $T \cup  (\{0\} \times I)$ in $\R^2
\times I$.  The latter projects to an ordinary $(k+1,n+1)$-tangle with
distinguished curve $P = \{0\} \times I$.  

It follows that 
the monoid of (isotopy classes of) affine
$(n,n)$-tangles is isomorphic to the monoid of regular isotopy classes
of affine $(n,n)$-tangle diagrams. 

\subsection{The ordinary  Kauffman tangle algebras.}
Let $\u k n$ denote the family of ordinary $(k,n)$-tangle
diagrams modulo regular isotopy.  Likewise, let $\uhat k n$ denote the
family of affine  $(k,n)$-tangle
diagrams modulo regular isotopy.  Then $\u n n$ and $\uhat n n$ are monoids for 
each~$n$.

\begin{definition}\rm\label{definition Kauffman tangle algebra}
 Let $R$ be any (commutative, unital) ring with distinguished elements
$\lambda$, $z$ and $\delta$, with $\lambda$ and $\delta$ invertible,
satisfying the relation
$$
\la\inv - \la= z(\delta - 1).
$$
 The (ordinary) {\em Kauffman tangle
algebra} $\kt {n, R}$  over $R$
 is the monoid algebra $R \ \u n n$ modulo the following
relations:
\begin{enumerate}
\item (Kauffman skein relation)
$$
\qquad\quad \inlinegraphic{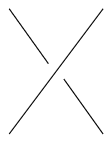} - \inlinegraphic[scale=.5]{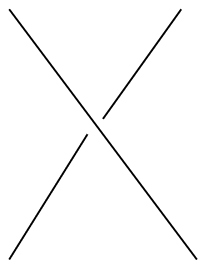} 
\quad = 
\quad
z\left( \inlinegraphic[scale=2]{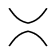} \quad - \quad
\inlinegraphic[scale=2]{id_smoothing}\right).
$$
Here, the figures indicate tangle diagrams which differ only in the
region shown and are identical outside this region.
\item (Untwisting relation)
$$\qquad\quad 
\inlinegraphic{right_twist} \quad = \quad \la \quad
\inlinegraphic{vertical_line} \quad\ \text{and} \quad\ 
\inlinegraphic{left_twist} \quad = \quad \la\inv \quad
\inlinegraphic{vertical_line}. 
$$
\item (Free loop relation)
 $$
T \cup \bigcirc = \delta T,
$$
where $T \cup\bigcirc $ is a tangle diagram consisting of the union
 of $T$ and an additional
closed loop having no crossing with $T$ and no self-crossings.
\end{enumerate}
\end{definition}

The {\em  generic Kauffman tangle algebra} (as considered   by Morton
and Tra\-czyk~\cite{Morton-Traczyk} and Morton and 
Wassermann~\cite{Morton-Wassermann})
 is the Kauffman tangle  algebra $\kt {n, \varLambda}$
 over the ring
$$
\varLambda = \Z[\labold^{\pm 1}, \zbold, \deltabold^{\pm 1}]/\langle \labold\inv 
- \labold=
\zbold(\deltabold - 1) \rangle.
$$
Here $\labold$, $\zbold$, $\deltabold$ are indeterminates.
For any ring $R$ as above, there is a homomorphism from $\varLambda$ to $R$ 
determined by
$\labold \mapsto \la$, $\zbold \mapsto z$, and $\deltabold \mapsto \delta$.
Thus the specialization $ \kt {n, \varLambda} \otimes_{\varLambda} R  $ makes 
sense.
It follows from the freeness of Kauffman tangle algebras~\cite{Morton-Wassermann} 
that
 for any $R$,
$$
\kt {n, R}  \cong \kt {n, \varLambda} \otimes_{\varLambda} R;
$$
see Corollary \ref{corollary-  freeness of ktn and isomorphism of
specializations} below.  This is actually a special case of a general universal
coefficient theorem for skein modules due to J. Przytycki; 
see~\cite[Lemma~5]{Przy1}.

\begin{remark}\rm
If $R$ is a ring with distinguished elements
$\lambda$, $z$, and $\delta$ as above, and $S \supseteq R$ is a ring
containing $R$, then the inclusion $R \ \u n n \rightarrow S\ \u n n$
induces an $R$-algebra homomorphism $\kt {n, R} \rightarrow \kt{n,S}$. 
This map is always injective; see Corollary~\ref{corollary-imbedding 
of ordinary BMW algebras}.
\end{remark}

\subsection{The affine Kauffman tangle algebras.}\label{subsection- affine KT 
algebra}
Let $R$ be a ring with distinguished elements $\la$, $z$ and $\delta$ as
above.  Our {\em preliminary} definition of the affine Kauffman tangle
algebra over $R$ is the monoid algebra $R\  \uhat n n$ modulo the Kauffman
skein relation, the untwisting relation, and the free loop relation, as
for the ordinary Kauffman tangle algebras.    Here it is
understood that none of the curves in the diagrams for these relations
represent a part of the flagpole.  Denote this algebra temporarily by
$\widehat K_{n, R}$.

For $r \ge 1$, let $\varTheta_r$  (resp. $\varTheta_{-r}$)  denote the (regular isotopy class of) the closed
curve with no self--crossings that winds $r$ times around the flagpole in the positive sense  (resp. in the negative sense). 
$$
\begin{array}{c c}
\includegraphics[scale=1.5]{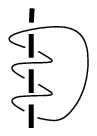} \\[-6pt] \vrule height20pt width0pt depth0pt 
{\scriptstyle \varTheta_3} \end{array} 
\qquad
\begin{array}{c}
\includegraphics[scale=1.5]{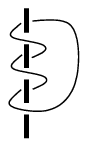}  \\[-6pt]  \vrule height18pt width0pt depth0pt  {\scriptstyle \varTheta_{-3}} 
\end{array}
$$

It is not
difficult to see that the curves $\varTheta_r$ generate $\widehat K_{0, R}$.  A
theorem of Turaev~\cite{Turaev-Kauffman-skein}  says that 
$\{ \varTheta_r : r \ge 1\}$ is algebraically independent in $\widehat
K_{0,\varLambda}$.   Hence $\widehat K_{0,\varLambda}$ is the polynomial algebra
over $\varLambda$ in the infinitely many variables $\varTheta_r$ for $r \ge 1$.

For any~$R$, one has an algebra map from $\widehat K_{0, R} $ to
$\widehat K_{0, \varLambda} \otimes_{\varLambda} R$; since the former is 
generated
by $\{ \varTheta_r : r \ge 1\}$ and the latter is the polynomial algebra
over $R$ in the variables $\{ \varTheta_r : r \ge 1\}$, it follows that the
map is an isomorphism.  

 Thus,  $\widehat K_{0, R} $  is the polynomial
algebra over $R$ in the variables $\varTheta_r$  for $r \ge 1$.

 We have $\varTheta_1 = \varTheta_{-1}$ in $\widehat K_{0, R} $,  but
$\varTheta_r$ is not ambient isotopic to  $\varTheta_{-r}$ for $r \ge 2$.
According to Turaev's theorem,  $\varTheta_{-r}$ is a polynomial  in the variables $\varTheta_k$,  $k \ge 1$,
\begin{equation} \label{theta minus r polynomial in thetas}
\varTheta_{-r} =  f_r(\varTheta_1, \varTheta_2, \dots).
\end{equation}
Later we will need  a little more information about the polynomials $f_r$.
For $a \ge 1, b \ge 0$,  let $\varTheta_{a, b}$ be the curve with $a$ positive 
windings around the flagpole, and one positive crossing, and $b$ negative windings, as in the first of the following figures.   Let  $\varTheta_{a, b}^{-}$ be the curve with the crossing reversed.  
$$
\begin{array}{c c}
\includegraphics[scale=1.3]{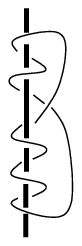} \\[-6pt] \vrule height20pt width0pt depth0pt 
{\scriptstyle \varTheta_{2, 3}} \end{array} 
\qquad
\begin{array}{c}
\includegraphics[scale=1.3]{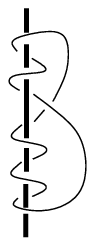}  \\[-6pt]  \vrule height18pt width0pt depth0pt  {\scriptstyle \varTheta_{2, 3}^{-}} 
\end{array}
$$
In particular, $\varTheta_{r, 0} = \la\inv \varTheta_r$; we interpret  $\varTheta_0$ as 
$\delta$.

\begin{lemma}  \label{lemma - recursion for f_r} 
Let $r, a, b \ge 1$.
\begin{enumerate}
\item $\varTheta_{-r} = \la \varTheta_{1, r-1}.$
\item $\varTheta_{a, b} =  \la^2 \varTheta_{a+1, b-1} + z(\varTheta_{a}\varTheta_{-b} - \varTheta_{a-b}).
$
\item
$
f_r(\varTheta_1, \varTheta_2, \dots) = 
\lambda^{2r-2}  \varTheta_r + z f'_r(\varTheta_1, \dots, \varTheta_{r-1}),
$
where $f'_r$ is a polynomial in $\varTheta_1, \dots, \varTheta_{r-1}$.
\end{enumerate}
\end{lemma}

\begin{proof}  Point (1) follows from introducing
 a twist at the top of $\varTheta_{-r}$:
 $$
\inlinegraphic[scale=1.5]{Theta3inv} = \la \inlinegraphic[scale=1.5]{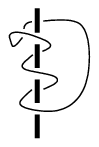} = \la \inlinegraphic[scale=1.5]{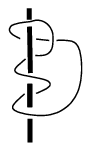}.
 $$
 The tangle  obtained by smoothing the crossing  in $\varTheta_{a, b}$
horizontally is  $\varTheta_a \varTheta_{-b}$,  and the tangle obtained by smoothing the crossing vertically is
$\varTheta_{a-b}$, while $\varTheta_{a, b}^{-} = \la^2 \varTheta_{a+1, b-1}$.
Thus the Kauffman skein relation gives
\begin{equation} \label{theta2}
\varTheta_{a, b} =  \la^2 \varTheta_{a+1, b-1} + z(\varTheta_{a}\varTheta_{-b} - \varTheta_{a-b}).
\end{equation}
An induction based on points (1) and (2) yields (3).
\end{proof}

We now return to the definition of the affine Kauffman tangle algebra.
Because $\widehat K_{0, R}$ is a polynomial algebra over $R$ and
$\widehat K_{n, R}$ is a $\widehat K_{0, R}$ module, it makes sense to just
absorb $\widehat K_{0, R}$ into the ground ring for the affine Kauffman tangle algebra, and this is what we~do.

\begin{definition}\rm\label{definition affine Kauffman tangle algebra}
Let $S$ be a ring containing distinguished elements 
$\la$, $z$,  $\delta$, and $q_r$  ($r \ge 1$), with 
$\lambda$ and $\delta$ invertible, such that 
the relation
$
\la\inv - \la= \break  z(\delta - 1)
$
holds.
The {\em affine Kauffman tangle
algebra} 
$$
\akt {n, S} = \akt {n, S}(\lambda, z, \delta, q_1, q_2, \dots)
$$
is the monoid algebra $S \uhat n n$ modulo the following relations:
\begin{enumerate}
\item The Kauffman skein relation, the untwisting relations, and the free
loop relation, as for the ordinary Kauffman tangle algebra. 
\item   $T\cup \varTheta_r = q_r T$, where $T\cup \varTheta_r$ is the union
of a tangle $T$ and a copy of the curve $\varTheta_r$, such that  there are no crossings between $T$ 
and~$\varTheta_r$.
\end{enumerate}
\end{definition}

\begin{remark}\rm
 Note that $\akt {0, S} \cong S$.
\end{remark}

If $R$ is a ring with distinguished elements $\lambda$, $z$, and
$\delta$ as above, and 
$\qbold_r$, $r
\ge 1$, are indeterminates, we   denote by  $\RHat$  the polynomial
ring $R[\qbold_1, \qbold_2, \dots ]$. 

\begin{remark}\rm\label{remark-injections of tangle algebras in affine tangle 
algebras}
 Let  $R$ be a ring with distinguished elements
$\lambda$,\ $z$, and $\delta$ as above, and $S \supseteq R$  a ring with 
additional
elements $q_1, q_2, \dots$.  Adding a flagpole to ordinary $(n,n)$-tangle
diagrams induces 
a
homomorphism of $R$-algebras $i_n : \kt {n,R} \rightarrow \akt {n, S}$. 
This homomorphism is always injective.  For the moment, we verify
injectivity in two special cases.
\begin{enumerate}
\item  $i_n : \kt {n,S} \rightarrow \akt {n, S}$ is injective.  In fact
the image of $i_n$ is the span of affine $(n,n)$-tangle diagrams having
no intersection with the flagpole.  We can define a map $o_n$  of such
diagrams to ordinary $(n,n)$-tangle diagrams by removing the flagpole.
The map $o_n$ induces an $S$-algebra homomorphism inverse to $i_n$; that is,
 $o_n\circ i_n$ is the identity on $\kt {n, S}$.
\item $i_n : \kt {n,R} \rightarrow \akt {n, \RHat}$ is injective.  In
fact, 
 removing the flagpole from affine $(n,n)$-tangle diagrams  and
mapping $\qbold_r \mapsto \delta$ determines a homomorphism of $R$-algebras 
$f_n :\akt {n,\RHat}
\rightarrow \kt {n, R}$, and $f_n \circ i_n$ is the identity on $\kt {n, R}$.
\end{enumerate}
Once we have verified that $\kt {n, R}$ imbeds in $\kt {n, S}$ (see Corollary 
\ref{corollary-imbedding of ordinary
BMW algebras}), it will follow  from (1) that $\kt {n, R}$  imbeds in $\akt {n, 
S}$ in general.
\end{remark}

The {\em  generic affine Kauffman tangle algebra}\/ is the Kauffman tangle  
algebra 
$$
\akt {n,\varLambdaHat} = \akt {n,
\varLambdaHat}(\lambdabold, \zbold, \deltabold, \qbold_1, \qbold_2,\dots)
$$
 over the ring $\varLambdaHat = \varLambda[\qbold_1, \qbold_2, \dots]$.

Since for any ring $S$ as above, we have a homomorphism $\varLambdaHat$ to
$S$, determined by the assignments $\lambdabold \mapsto \lambda$, 
$\zbold \mapsto z$, $\deltabold \mapsto \delta$, and $\qbold_i \mapsto
q_i$,  the specialization
$\akt {n, \varLambdaHat} \otimes_{\varLambdaHat} S$ makes sense.
We will eventually show that 
$$
\akt {n, S}  \cong \akt {n, \varLambdaHat} \otimes_{\varLambdaHat} S;
$$
see Corollary \ref{corollary- freeness affine specializations}.

\begin{remark}\rm %
\label{remark - Brauer specialization of affine algebra}
We will  use the following specialization below, in connection with an affine 
version of the Brauer
algebra: Set $R =  \Z[\delta\pmone]$.  We have a homomorphism of $e: \varLambda 
\rightarrow R$ determined by
 $\la \mapsto 1$, $z \mapsto 0$, $\delta \mapsto \delta$.
 Then $\akt {n,\RHat} $ is the quotient of $\RHat\ \uhat n n$ by the relations:
\begin{enumerate}
\item (Kauffman skein relation)
$
\inlinegraphic{pos_crossing} = \inlinegraphic[scale=.5]{neg_crossing}
$
for a crossing of strands neither of which represents a part of the flagpole.
\bigskip
\item (Untwisting relation)\vadjust{\vskip4pt}
$$\qquad\ 
\inlinegraphic{right_twist} \quad = \quad  \quad
\inlinegraphic{vertical_line} \quad\ \text{and} \quad\
\inlinegraphic{left_twist} \quad = \quad  \quad
\inlinegraphic{vertical_line} 
$$
\item (Removal of closed loops)
 $T \cup \bigcirc = \delta T,$  where $T \cup
\bigcirc $ is a tangle diagram consisting of the union of $T$ and an additional
closed loop having no crossings with the flagpole, and 
 $T\cup \varTheta_r = \qbold_r T$,  for $r\ge 1$.
 \end{enumerate}
 Since in this specialization,  two tangle diagrams related by crossing changes are equivalent,  for any $r \ge 1$ and or any closed curve $c$ which loops $r$ times around the flagpole,  $c = \varTheta_r = \qbold_r$.
\end{remark}

\subsection{The elements $ X_r $.}
We define elements $X_r \in \akt {n, S}$ for $1 \le r \le n$ by
$$
X_r = (G_{r-1} \cdots G_2 G_1) X_1 (G_1 G_2 \cdots G_{r-1}).
$$
For example, 
$$
X_4 =  \inlinegraphic{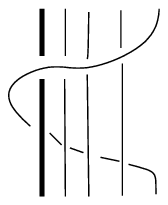} \quad\ \text{and} \quad\ X_4\inv = 
\inlinegraphic{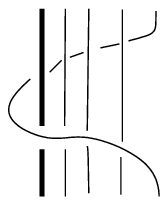}
$$

\subsection{Symmetries.}\label{subsection-symmetries of kauffman tangle 
algebra}
The map on ordinary (respectively, affine) tangle diagrams which flips diagrams 
top to bottom
induces an anti-automorphism  $\alpha$ of $\kt {n, R}$ (respectively, of 
$\akt{n, S}$).
 The anti-automorphism $\alpha$
fixes  $E_i$, $G_i$, and $X_1$.

There is an isomorphism  $\beta : \kt {n, R}(\lambda, z, \delta) \rightarrow  
\kt {n, R}(\la\inv, -z, \delta)  $
determined by the map of ordinary tangle diagrams that reverses all crossings.

For $r \ge 1$, define $q_{-r} = f_r(q_1, \dots, q_r)$,  where
$f_r$ is the polynomial of equation \ref{theta minus r polynomial in thetas}.
There is an isomorphism  
$$\beta : \akt {n, S}(\lambda, z, \delta, q_1, q_2, \dots) 
\rightarrow  
\akt {n, S}(\la\inv, -z, \delta, q_{-1}, q_{-2}, \dots)  $$
determined by the map of affine tangle diagrams that reverses all crossings 
(also crossings of ordinary strands with the flagpole).

The ordinary Kauffman tangle algebra $\kt {n, R}$ has an automorphism 
$\varrho_n$ which
flips ordinary tangle diagrams left to right.  The automorphism $\varrho_n$ 
satisfies
$\varrho_n(E_i) = E_{n+1 -i}$ and $\varrho_n(G_i) = G_{n+1 -i}$. 

\subsection{Inclusions, conditional expectations, and  
trace.}\label{subsection- inclusions and trace}
The  observations in this section apply equally to the ordinary and affine
Kauffman tangle algebras.   We fix a ring $S$ containing distinguished
elements
$\la$, $z$,  $\delta$, and $q_r$ ($r  \ge 1)$ as above, and we let $\akt n$ 
denote $\akt {n,S}$.

For $n \ge 2$, the map $\iota$ from affine $(n-1,n-1)$-tangle diagrams to
affine
$(n, n)$-tangle diagrams that adds an additional strand on the right
without adding any crossings:
$$
\iota: \quad \inlinegraphic{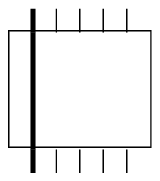} \quad \mapsto \quad 
\inlinegraphic{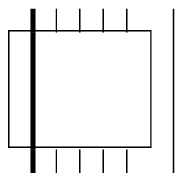}
$$
respects regular isotopy, composition of tangle diagrams, and the
relations of the  affine Kauffman tangle algebras, so induces a
homomorphism $\iota : \akt {n-1} \rightarrow \akt {n}$.
Note that $\alpha \circ \iota = \iota \circ \alpha$, where $\alpha$ is the 
anti-automorphism
of  the affine Kauffman tangle algebras described in 
Section~\ref{subsection-symmetries of kauffman tangle algebra}.

\def\cl{\mathop{\rm cl}\nolimits}

The map of affine $(n,n)$-tangle diagrams to affine
$(n-1, n-1)$-tangle diagrams that ``closes" the rightmost strand, 
without adding any crossings:
$$
{\rm cl}_n: \quad \inlinegraphic{tangle_box} \quad \mapsto \quad 
\inlinegraphic{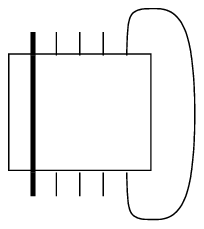}
$$
respects regular isotopy and the
relations of the  affine Kauffman tangle algebras, so induces an
$S$-linear map $\akt n \rightarrow \akt {n-1}$.  We define
$\eps_n : \akt n \rightarrow \akt {n-1}$ by 
$$
\eps_n(T) = \delta\inv{\rm cl}_n(T).
$$
 Note that $\eps_n\circ\alpha = \alpha\circ\eps_n$, where
 $\alpha$ is the anti-automorphism
of  the affine Kauffman tangle algebras described above.

We have
$$
\eps_n \circ \iota(x) = x
$$
for $x \in \akt {n-1}$.  In particular, the map $\iota$ of $\akt {n-1}$
into $\akt {n}$ is injective.  Consequently, we will drop the notation
$\iota$ and regard $\akt {n-1}$ as a subalgebra of $\akt n$.

More
generally,
$$
\eps_n( x y) =  x \eps_n(y),
\quad\
\eps_n( y x ) =   \eps_n(y) x
$$
for $x \in \akt {n-1}$ and $y \in \akt n$;  that is, 
$\eps_n$ is a $\akt {n-1}$-$\akt {n-1}$-bimodule map, or
{\em conditional expectation}.

\begin{remark}\rm
One has the formula
$$
E_n x E_n = \cl_n(x) E_n
$$
for $x$ an (affine) $(n,n)$-tangle diagram.  Equivalently, the idempotent 
$E_n/\delta$
implements the conditional expectation:
$$
(E_n/\delta) x (E_n/\delta) = \eps_n(x) (E_n/\delta)
$$
for $x \in \akt n$.
\end{remark}

We define $\eps : \akt n \rightarrow \akt 0 \cong S$ by
$$
\eps = \eps_1 \circ \cdots \circ \eps_n.
$$
Equivalently, if we define the closure of an affine $(n,n)$-tangle
diagram to be the affine $(0,0)$-diagram obtained by closing all strands
without introducing any new crossings:
$$
{\rm cl} :  \inlinegraphic{tangle_box} \quad \mapsto \quad 
\inlinegraphic{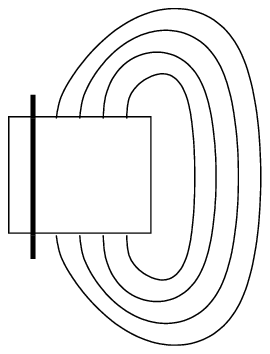}
$$
then $\eps(x) = \delta^{-n} {\rm cl}(x)$.  Note that
$$
\eps \circ \iota = \eps, \quad\ \eps\circ \eps_n = \eps.
$$

We will show below that $\akt n$ is generated as a unital algebra by
$X_1^{\pm 1}$, $G_i^{\pm 1}$ and $E_i$ ($1 \le i \le n-1$).  In particular
$\akt 1$ is generated by $X_1^{\pm 1}$, so is commutative.  

For $n \ge 2$, let $\eps_{n,1}$ denote the map
$$
\eps_{n,1} = \eps_2 \circ \dots \circ \eps_n : \akt n \rightarrow
\akt 1.
$$
Then $\eps_{n,1}$ is a conditional expectation, so for $x \in \akt n$,
$$
\eps_{n,1}( X_1^{\pm 1} x) = X_1^{\pm 1}  \eps_{n,1}(x) =   \eps_{n,1}(x)
X_1^{\pm 1} =
\eps_{n,1}(x  X_1^{\pm 1}).
$$
Since $\eps = \eps_1 \circ \eps_{n,1}$, we have
$$
\eps(X_1^{\pm 1} x) = \eps( x X_1^{\pm 1}).
$$

For any affine $(n,n)$-tangle  diagram $T$ and any $r\le n-1$, the closure of 
the
affine tangle diagram
$ E_r T$ is isotopic to the closure of the affine tangle diagram $T E_r$, and
similarly for $E_r$ replaced by $G_r^{\pm 1}$.  The following figure
illustrates this for $G_r\inv$:
$$
\inlinegraphic{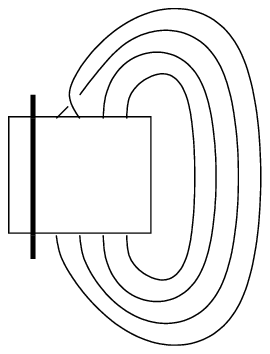} \quad = \quad 
\inlinegraphic{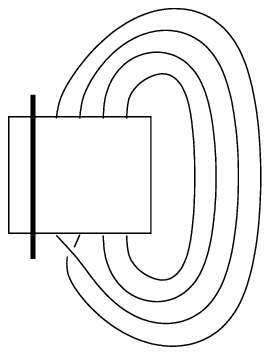}
$$
It follows that for $x \in \akt n$ and $r \le n-1$,
$$
\eps( E_r x) =  \eps( x E_r) ,\quad\ \eps( G_r^{\pm 1} x) =  \eps( x G_r^{\pm 
1})
$$

\begin{proposition}
 $\eps : \akt n \rightarrow \akt 0 \cong S$
is a trace.  That is,  $\eps$ is $S$-linear and $\eps(x y ) = \eps(y x)$ for all
$x, y \in \akt n$.
\end{proposition}

\begin{proof}
We have shown that the trace property
$$
\eps(x y ) = \eps(y x)
$$
holds when $y$ is arbitrary and $x$ is one of $X_1^{\pm 1}$, $G_i^{\pm
1}$ and $E_i$ ($1 \le i \le n-1$).   It follows that the trace property
also holds when $x$ is any product of the elements $X_1^{\pm 1}$, $G_i^{\pm
1}$ and $E_i$.  Since these elements generate $\akt n$ as a unital
algebra, the property holds for all $x$.
\end{proof}

\begin{remark}\rm\label{remark- Markov property of trace}
 One easily checks by picture proofs that 
$\eps_{r+1}(T G_r^{\pm 1}) = (\la^{\pm 1}/\delta)  T, $
and hence 
$$
\eps(T G_r^{\pm 1}) =  \eps(T) \la^{\pm 1}/\delta
$$
for $T \in \akt {r}$. Likewise,
$ \eps_{r+1}(T  E_r) =(1/\delta) T ,$
and hence
$$
\eps(T E_r) = \eps(T)/\delta
$$
for $T \in \akt {r}$.  Moreover, if $X_r'$
 denotes
$$
X_r' = G_{r-1} \cdots G_1 X_1 G_1\inv \cdots G_{r-1}\inv,  
$$
then for $s\ge 1$ and   $T \in \akt {r}$, one has
$\eps_{r+1}( T( X_r')^s) = q_{s} T$.
Consequently, 
$$
\eps( T( X_r')^s) = q_{s} \eps(T).
$$
It follows from Lemma \ref{lemma - recursion for f_r} that
$$
\eps( T( X_r')^{-s}) = f_{s}(q_1, \dots, q_s) \eps(T).
$$
A trace with these properties is usually called a {\em Markov trace}.  The 
terminology originated 
in~\cite{Jones}.
\vadjust{\eject}
\end{remark}

\def\sectionmark#1{.}
\section{The affine Birman--Wenzl--Murakami algebra}\label{sec3}

\vskip0pt\vskip-12pt\vskip0pt
\subsection{The Birman--Wenzl--Murakami algebra.}
The Birman--Wenzl--Murakami algebras were introduced independently by
Birman and Wenzl~\cite{Birman-Wenzl} and by Murakami~\cite{Murakami-BMW}
as an algebraic setting for the Kauffman link invariant~\cite{Kauffman}.
These algebras are known to appear as centralizer algebras for the
quantum universal enveloping algebras of ${\rm sp}(2n, \C)$ or ${\rm so}(n, \C)$
acting on tensor powers of the vector representation.  They are
deformations of Brauer's centralizer algebras (see~\cite{Wenzl-Brauer}),
and are extensions of the Hecke algebras of type $A$, as will be explained
below.

{\spaceskip 0.33em plus0.2em minus0.17em The presentation which we give here 
follows Morton and 
Wassermann~\cite{Morton-Wassermann}.  The parameters differ slightly from those
used by Birman and Wenzl.}

As before, let $R$ be a commutative unital ring with
invertible elements $\la$ and $\delta$ and an element $z$ satisfying
$$
\la \inv - \la = z(\delta - 1).
$$

\begin{definition}\rm
The {\em Birman--Wenzl--Murakami}\/ algebra 
$\bmw {n, R}$  is the $R$-algebra with generators $g_i^{\pm 1}$  and
$e_i$ ($1 \le i \le n-1$) and relations:
\begin{enumerate}
\item (Inverses) \hods $g_i g_i\inv = g_i\inv g_i = 1$.
\item (Idempotent relation)\hods $e_i^2 = \delta e_i$.
\item (Braid relations) \hods $g_i g_{i+1} g_i = g_{i+1} g_i g_{i+1}$ 
and $g_i g_j = g_j g_i$ if $|i-j|  \ge 2$.
\item (Commutation relations)  \hods $g_i e_j = e_j g_i$  and
$e_i e_j = e_j e_i$  if $|i-j|\ge 2$. 
\item (Tangle relations)\hods $e_i e_{i\pm 1} e_i = e_i$, $g_i
g_{i\pm 1} e_i = e_{i\pm 1} e_i$, and $ e_i  g_{i\pm 1} g_i=   e_ie_{i\pm 1}$.
\item (Kauffman skein relation)\hods  $g_i - g_i\inv = z(e_i -1)$.
\item (Untwisting relations)\hods $g_i e_i = e_i g_i = \la\inv e_i$,
and $e_i g_{i \pm 1} e_i = \la e_i$.
\end{enumerate}
\end{definition}

\begin{remark}\rm\label{remark-on the definition of BMW}
$\hbox{}$

\begin{enumerate}
\item If $z$ is taken to be invertible, then $e_i = 1 + z\inv(g_i\inv -g_i)$. 
 In this case, several of the relations are redundant.  Moreover,
the algebra is a quotient of the braid group algebra.
\item  In the parametrization of~\cite{Wenzl-BCD}, $z$ is replaced by 
$  q - q\inv$, where $q$ is another indeterminate, and
$q - q\inv$ is assumed to be invertible.  The Kauffman skein relation
then becomes a cubic relation for  the braid generators~$g_i$.
\item The main result of~\cite{Morton-Wassermann} is that the
Birman--Wenzl--Murakami  algebra $\bmw {n,\varLambda}$ defined over $\varLambda$ 
is
isomorphic to the Kauffman tangle algebra $\kt {n, \varLambda}$.  

\quad\ In \cite{Kauffman}, Theorem 4.4,  Kauffman gives a presentation for the 
algebra of tangles
 generated by the $E_i$'s, $G_i$'s, and free loops.  It is not clear to us how 
this result
 is related to the Morton--Wassermann theorem.
\end{enumerate}
\end{remark}

\begin{remark}\rm\label{remark- injection-  shift- and reversal}
$\hbox{}$
\begin{enumerate}
\item
The assignment $e_i \mapsto e_i$, $g_i \mapsto g_i$ defines a homomorphism 
$\iota$ from $\bmw {n, R} $ to $\bmw {n+1, R}$, since the relations are 
preserved.   It is not evident
that $\iota$ is injective, but this will follow eventually from the isomorphism 
of $\bmw {n, R}$ with
$\kt {n, R}$;   see Theorem~\ref{theorem- isomorphism for ordinary BMW}.
\item Define $S(x) : \bmw {n, R} \rightarrow  \bmw {n+1, R}$ by 
$$
S(x) = {\rm Ad}(g_{n}g_{n-1}\cdots g_{1})(x).
$$
Then $S$ is an injective homomorphism.  It follows from  the relations
that
$S(e_{i}) = e_{i+1}$ and $S(g_{i}) = g_{i+1}$ for $i \le n-1$. 
The map $S$ is called the 
{\em shift endomorphism}.
\item  The assignment $e_i \mapsto e_{n-i}$, 
$g_i \mapsto g_{n-i}$  determines an automorphism $\varrho_n$ of $\bmw {n, R} $, 
as
the relations are preserved by this assignment.  The shift map $S: \bmw {n, R} 
\rightarrow \bmw
{n+1, R}$ satisfies $S = \varrho_{n+1}\circ\varrho_{n}$.
\item  The assignment $g_i \mapsto g_i$ and $e_i \mapsto e_i$ determines
an anti-auto\-mor\-phism $\alpha$ of $\bmw {n , R} $.
\item  The assignment $g_i \mapsto g_i\inv$, $e_i \mapsto e_i$ determines
an isomorphism $\beta$ from $\bmw {n , R} (\la, z, \delta)$ to 
$\bmw {n , R} (\la\inv,  -z, \delta)$. 
\end{enumerate}
\end{remark}

We have the following relation of the BMW algebras to the Hecke algebras of type 
$A$:
 
\begin{proposition}\label{proposition- ideal of en and Hecke quotient}
Let $J_{n, R}$ denote the ideal $ \bmw {n, R}  e_{n-1} \bmw {n, R}  $
in $\bmw {n, R}$.
\begin{enumerate}
\item  $J_{n, R}$ is the ideal in $\bmw {n,
R}
$ generated by
$\{e_1, \dots, e_{n-1}\}$.
\item  $\bmw {n, R} /J_{n, R}$ is isomorphic to the Hecke
algebra over $\varLambda$
 with generators
$g_i\pmone$  $(1 \le i \le n-1)$ satisfying the braid relations and
the quadratic relation
$$
g_i^2 = 1 - z g_i.
$$
\end{enumerate}
\end{proposition} 

\begin{proof}  For (1),  it follows from the relation $e_j e_{j+1} e_j = e_j$ 
that all the $e_i$'s are
elements of  $\bmw n e_{n-1} \bmw n$.    Statement (2) is evident, since the 
$e_i$'s are zero in
the quotient. 
\end{proof}

\begin{remark}\rm
For $\omega \in \S_n$, let  $\omega = s_{i_1} \cdots s_{i_r}$ be
a reduced expression for $\omega$ (in terms of the generators $s_i$ of the
symmetric group). Let
$g_\omega$ be the corresponding word in the generators $g_i$ of $\bmw {n, R}$.
It follows from the proposition and well known facts concerning the Hecke 
algebra
that 
$g_\omega$ is well defined modulo $J_{n, R} $.
In fact, $g_\omega$ is well defined in $\bmw {n, R}$, not merely well
defined modulo the ideal   $J_{n,R} $;  see 
Proposition~\ref{proposition-positive permutation braid}.
\end{remark}

Recall that  $\varLambda$ denotes  
 the ring
$$
\varLambda = \Z[\labold^{\pm 1}, \zbold, \deltabold^{\pm 1}]/\langle \labold\inv 
- \labold=
\zbold(\deltabold - 1) \rangle,
$$
where $\labold$, $\zbold$, $\deltabold$ are indeterminates.
The {\em generic Birman--Wenzl--Murakami algebra} is the BMW algebra $\bmw 
{n, \varLambda}$ over the
ring $\varLambda$.  Given a ring $R$ as above, we have two possible 
specializations to an algebra over
$R$, namely $\bmw {n, R}$ and $\bmw {n, \varLambda} \otimes_\varLambda R$.  We 
will show later that these
are in fact isomorphic; see Corollary \ref{corollary- isomorphic specializations 
ordinary BMW}.
 
\subsection{The affine Birman--Wenzl--Murakami  
algebra.}\label{subsection-affine BMW}
Let $S$ be a ring with distinguished elements $\lambda$, $z$, and $\delta$ as 
above, and with additional
elements $q_1, q_2, \dots$.

\begin{definition}\rm\label{definition affine BMW}
 The {\em affine
Birman--Wenzl--Murakami} algebra
$\abmw  {n, S}$ is the
$S$ algebra with generators $x_1^{\pm 1}$, $g_i^{\pm 1}$  and
$e_i$ ($1 \le i \le n-1$) and relations:
\begin{enumerate}
\item (Inverses)\hods $g_i g_i\inv = g_i\inv g_i = 1$ and 
$x_1 x_1\inv = x_1\inv x_1= 1$.
\item (Idempotent relation)\hods $e_i^2 = \delta e_i$.
\item (Affine braid relations) 
\begin{enumerate}
\item[\rm(a)] $g_i g_{i+1} g_i = g_{i+1} g_ig_{i+1}$ and 
$g_i g_j = g_j g_i$ if $|i-j|  \ge 2$.
\item[\rm(b)] $x_1 g_1 x_1 g_1 = g_1 x_1 g_1 x_1$ and $x_1 g_j =
g_j x_1 $ if $j \ge 2$.
\end{enumerate}
\item[\rm(4)] (Commutation relations) 
\begin{enumerate}
\item[\rm(a)] $g_i e_j = e_j g_i$  and
$e_i e_j = e_j e_i$  if $|i-
j|
\ge 2$. 
\item[\rm(b)] $x_1 e_j = e_j x_1$ if $j \ge 2$.
\end{enumerate}
\item[\rm(5)] (Affine tangle relations)\vadjust{\vskip-2pt\vskip0pt}
\begin{enumerate}
\item[\rm(a)] $e_i e_{i\pm 1} e_i = e_i$,
\item[\rm(b)] $g_i g_{i\pm 1} e_i = e_{i\pm 1} e_i$ and
$ e_i  g_{i\pm 1} g_i=   e_ie_{i\pm 1}$.
\item[\rm(c)\hskip1.2pt] For $r \ge 1$, $e_1 x_1^{ r} e_1 = q_r e_1$. 
\vadjust{\vskip-
2pt\vskip0pt}
\end{enumerate}
\item[\rm(6)] (Kauffman skein relation)\hods  $g_i - g_i\inv = z(e_i -1)$.
\item[\rm(7)] (Untwisting relations)\hods $g_i e_i = e_i g_i = \la \inv e_i$
 and $e_i g_{i \pm 1} e_i = \la  e_i$.
\item[\rm(8)] (Unwrapping relations)\hods $e_1 x_1 g_1 x_1 = \la\inv e_1 = x_1 
g_1 x_1 e_1$.
\end{enumerate}
\end{definition}

\vskip0pt\vskip-3pt\vskip0pt
\begin{remark}\rm
\label{remark- inclusions and symmetries affine BMW}
 Let  $R$ be a ring with distinguished
elements $\lambda$, $z$, and $\delta$ as above, and $S \supseteq R$  a ring with 
additional
elements $q_1, q_2, \dots$.\vadjust{\vskip-2pt\vskip0pt}  
\begin{enumerate}
\item The assignment $e_i \mapsto e_i$, $g_i \mapsto g_i$, 
$x_1 \mapsto x_1$ defines a homomorphism $\iota$ from $\abmw  {n, S}$ to $\abmw  
{n+1,
S}$, since the relations are preserved.   It is not evident that $\iota$ is 
injective,
but this will follow from the isomorphism $\abmw  {n, S} \cong \akt {n, S}$; 
see Corollary \ref{corollary- injectivity of iota from affine bmw n-1 to affine 
bmw n}.
\item    The assignment $e_i \mapsto e_i$, $g_i \mapsto g_i$ 
defines an $R$-algebra 
homomorphism 
$i_n: \bmw {n, R} \rightarrow \abmw  {n, S}$.  This map is always injective, as 
we will
verify later.  (For the moment we verify that  $i_n: \bmw {n, R} \rightarrow 
\abmw  {n,
\RHat}$ is injective.  In fact, the assignment $e_i \mapsto e_i$, $g_i \mapsto 
g_i$, $x_1
\mapsto 1$, $\qbold_i \mapsto \delta$  defines an $R$-algebra  homomorphism from 
$f_n:
\abmw  {n,
\RHat}  \rightarrow \bmw {n, R}$, and $f_n \circ i_n$ is the identity on $\bmw 
{n,R}$.)

\quad\ Note that  the following diagram commutes:\vadjust{\vskip-4pt}
$$\qquad\ 
\begin{diagram}
\bmw {n, R}      &\rTo^{\scriptstyle i_n}    & \abmw  {n, S} \\
\dTo<{\scriptstyle \iota}                  &        & \dTo<{\scriptstyle \iota}  
\\
\bmw {n +1, R}                  &\rTo^{\scriptstyle i_{n+1} }   & \abmw  {n+1, 
S}
\end{diagram}
$$
\item  The assignment $g_i \mapsto g_i$, $e_i \mapsto e_i$, $x_1 \mapsto x_1$ 
determines
an anti-automorphism $\alpha$ of $\abmw {n , S} $.  One has $\alpha \circ i_n = 
i_n \circ \alpha$, where
$i_n : \bmw {n, R} \rightarrow \abmw {n, S}$ and $\alpha$ also denotes the 
anti-automorphism of $\bmw
{n, R}$ of Remark~\ref{remark- injection-  shift- and reversal}.

\item  For $r \ge 1$, let $q_{-r} = f_r(q_1, \dots, q_r)$,  where
$f_r$ is the polynomial of equation \ref{theta minus r polynomial in thetas}.
The assignment $g_i \mapsto g_i\inv$, $e_i \mapsto e_i$, $x_1 \mapsto 
x_1\inv$ determines
an isomorphism 
$$\beta : \abmw {n , S } (\la, z, \delta, , q_1, q_2, \dots ) \to 
\abmw {n , S} (\la\inv,  -z, \delta,  q_{-1}, q_{-2}, \dots).$$   
The following diagram 
commutes:\vadjust{\vskip-4pt}
$$\qquad\ 
\begin{diagram}
\bmw {n, R}(\la, z, \delta)     &\rTo^{\scriptstyle i_n}    & \abmw  {n, S} 
(\la, z, \delta, q_1, q_2, \dots)\\
\dTo<{\scriptstyle \beta}                  &        & \dTo<{\scriptstyle \beta}  
\\
\bmw {n , R}  (\la\inv,  -z, \delta)                &\rTo^{\scriptstyle i_{n} }   
& \abmw  {n, S}(\la\inv,  -z,
\delta, q_{-1}, q_{-2}, \dots)
\end{diagram}
$$

\end{enumerate}
\end{remark}

\begin{lemma}\label{lemma-quadratic relation}
 The Kauffman skein relation
implies
$$
g_i^2 = 1 + \la\inv z e_i - z g_i,\quad\
g_i^{-2} = 1 - \la z e_i + z g_i.
$$
\end{lemma}

\begin{proof}  Multiply by $g_i$ or $g_i\inv$ and simplify using the
untwisting relation.
\end{proof}

\begin{lemma}\label{lemma-commutation relations 1}
$\hbox{}$
\begin{enumerate}
\item  $g_{i}^{\pm 1} (g_{i+1} g_{i}) = (g_{i+1} g_{i}) g_{i+1}\pmone$.
\item   $e_{i} (g_{i+1} g_{i}) = (g_{i+1} g_{i}) e_{i+1}$.
\item $g_i g_{i+1} e_i e_{i+2} = g_{i+2}  g_{i+1} e_i e_{i+2}$.
\item  If $i < m$, then
\begin{eqnarray*}
\qquad\ g_{i}\pmone (g_{m} g_{m-1}\cdots g_{1}) &=&  
(g_{m} g_{m-1}\cdots g_{1}) g_{i+1}\pmone ,\\
g_{i+1}\pmone (g_{1} g_{2}\cdots g_{m}) &=&  
(g_{1} g_{2}\cdots g_{m}) g_{i}\pmone.
\end{eqnarray*}
\item  If $i < m$, then 
\begin{eqnarray*}
\qquad\ e_{i} (g_{m} g_{m-1}\cdots g_{1}) &=&  
(g_{m} g_{m-1}\cdots g_{1}) e_{i+1},\\
e_{i+1} (g_{1} g_{2}\cdots g_{m}) &=&  
(g_{1} g_{2}\cdots g_{m}) e_{i} .
\end{eqnarray*}
\end{enumerate}
\end{lemma}

\begin{proof}  The first statement is just the braid relation.  The second 
results from
two applications of the tangle relation (5b)  from Definition 
\ref{definition affine BMW}:
$$
 e_{i} (g_{i+1} g_{i}) = e_{i} e_{i+1} = (g_{i+1} g_{i}) e_{i+1}.
$$
Statement (3) also results from two applications of the tangle relation (5b)  
from 
Definition \ref{definition affine BMW}:
$$
 g_{i}  g_{i+1} e_{i} e_{i+2} =  e_{i+1} e_{i} e_{i+2} = e_{i+1}e_{i+2} e_{i}  =
 g_{i+2} g_{i+1} e_{i+2} e_i.
$$
The first parts of statements (4) and (5) follow from  statements  (1) and 
(2) 
and the commutation relations (3a) and
(4a) of Definition \ref{definition affine BMW}.  The second parts of  
statements (4) and (5)  follow by applying the
anti-auto\-morphism~$\alpha$.
\end{proof}

\subsection{The elements  $x_i$.}
For $2 \le r \le n$ define
$$
 x_r = (g_{r-1} \cdots g_1) x_1( g_1 \cdots g_{r-1}).
$$

\begin{proposition}\label{proposition-xr relations1}
\mbox{}
\begin{enumerate}
\item For all $r$ and $j \not\in  \{r, r-1\}$,  $g_j  x_r =  x_r g_j$. 
\item  For all $r$ and $j \not\in  \{r, r-1\}$,   
$e_j  x_r =  x_r e_j$.
\item  For all $r$ and $j$, $x_j x_r = x_r x_j$.
\end{enumerate}
\end{proposition}

\begin{proof}  For $j < r-1$, 
relations  (1) and (2)   follow by applying Lemma 
\ref{lemma-commutation relations 1}(4,5) and 
Definition \ref{definition affine BMW}(3b,4b).

For $j > r$,  relations (1) and (2)  follow from  Definition \ref{definition 
affine BMW}(3,4).

For part (3), we first observe that $x_1$ commutes with $x_r$ for $r \ge 2$, 
by the braid relations
Definition \ref{definition affine BMW}(3b).  Namely, 
\begin{eqnarray*}
x_1 x_r &=&x_1(g_{r-1} \cdots g_2 g_1) x_1 (g_1 g_2 \cdots g_{r-1})
=  (g_{r-1} \cdots g_2) x_1 g_1 x_1 g_1 (g_2 \cdots g_{r-1})\\
&=&  (g_{r-1} \cdots g_2) g_1 x_1 g_1 x_1 (g_2 \cdots g_{r-1})
= (g_{r-1} \cdots g_2 g_1) x_1 (g_1g_2 \cdots g_{r-1}) x_1 \\
&=& x_r x_1.
\end{eqnarray*}
Finally, if $j < r$, then
\[
x_j x_r = (g_{j-1}\cdots g_{1}) x_{1}(g_{1}\cdots g_{j-1}) x_{r} 
= x_{r}(g_{j-1}\cdots g_{1}) x_{1}(g_{1}\cdots g_{j-1}) = x_{r} x_{j},
\]
since $x_{r}$ commutes with $g_{i}$  for $i \le j-1$  and with $x_{1}$.
\end{proof}

\begin{proposition}\label{proposition-xr relations2} 
\mbox{}
\begin{enumerate}
\item  For $r < n$, $g_r x_r = x_{r+1} g_r \inv$ and $g_r \inv x_{r+1} = x_r 
g_r$.
\item For $r < n$, 
\begin{eqnarray*}
\qquad\ g_r x_{r+1} &=& x_r g_r -z x_{r+1} + z\lambda e_r x_{r+1},\\
g_r \inv x_r &=& x_{r+1} g_r\inv + z x_r - z e_r x_r.
\end{eqnarray*}
\item  For $r < n$, $e_r  x_r = \lambda^{-2}  e_r x_{r+1}\inv$ and $x_r e_r = 
\la^{-2} x_{r+1}\inv e_r$.
\item For $r < n$, $e_{r} x_r \inv = \lambda^2  e_{r} x_{r+1} $ and 
$ x_r \inv e_{r} = \lambda^2  x_{r+1}  e_{r} $.
\end{enumerate}
\end{proposition}

\begin{proof}
Statement (1) follows from the definition of $x_r$.  Statement (2) uses the 
Kauffman skein relation, the quadratic relation of Lemma \ref{lemma-quadratic 
relation}, and part~(1).
Statement (3) is equivalent to 
$$
e_r x_r x_{r+1} = \lambda^{-2}  e_r.
$$
For $r = 1$, this follows from 
Definition \ref{definition affine BMW}(7,8).
 The proof is completed by induction on $r$.  For $r \ge 2$, one has (using 
braid and tangle relations)
$$\eqalign{
e_r x_r x_{r+1}={}& e_r (g_{r-1} \cdots g_1) x_1 (g_1 \cdots g_{r-2}(g_{r-1} g_r  
g_{r-1}) g_{r-2} \cdots g_1) x_1( g_1 \cdots g_r) \cr
={}& e_r (g_{r-1} \cdots g_1) x_1 (g_1 \cdots g_{r-2}(g_{r} g_{r-1}  g_{r}) 
g_{r-2} \cdots g_1) x_1( g_1 \cdots g_r) \cr
={}& (e_r g_{r-1}g_r)( g_{r-2} \cdots g_1) x_1 (g_1 \cdots g_{r-2}g_{r-1}   
g_{r-2} \cdots g_1)\cr
&\cdot x_1( g_1 \cdots g_{r-2})(g_{r} g_{r-1}g_r)\cr
={}& (e_r e_{r-1})( g_{r-2} \cdots g_1) x_1 (g_1 \cdots g_{r-2}g_{r-1}   g_{r-2} 
\cdots g_1)\cr
&\cdot x_1( g_1 \cdots g_{r-2})(g_{r-1} g_{r}g_{r-1})  \cr
={}& e_r e_{r-1} x_{r-1} x_r   g_r g_{r-1}   \cr
={}& \lambda^{-2} e_r e_{r-1} g_r g_{r-1}   \quad\ 
\mbox{by the induction assumption}\cr
={}& \lambda^{-2} e_r e_{r-1} e_r = \lambda^{-2} e_r.\cr}
$$
Statement (4) follows from (3) by applying the isomorphism 
$\beta$ of Remark
\ref{remark- inclusions and symmetries affine BMW}.
\end{proof}

The next proposition contains general unwrapping and affine braid relations:

\begin{proposition}\label{proposition- unwrapping 1}
\label{proposition- affine braid relation 1}
For all $n \ge 1$ and all $r \ge 1$, 
\begin{enumerate}
\item  $g_{n}x_{n}g_{n} x_{n} = x_{n}g_{n} x_{n} g_{n}$.
\item  $e_{n} x_{n} g_{n} x_{n} = \la\inv e_{n}$.
\item  $e_{n} x_{n}^{r} g_{n} x_{n} =  \la^{-2} e_{n} x_{n}^{r-1} g_{n}\inv$.
\item  $e_{n} x_{n}^{-r} g_{n}\inv x_{n}\inv = \la^2 e_{n} x_{n}^{-r+1} g_{n}$.
\end{enumerate}
\end{proposition}

\begin{proof}  
 Statement (1) is equivalent to 
$x_{n+1} x_n = x_n x_{n+1}$, so follows from 
Proposition \ref{proposition-xr 
relations1}.
 Statement (2) follows from  (3) and the untwisting relation of Definition 
\ref{definition affine BMW}.   For statement (3), we have
$$
\begin{aligned}
e_{n} x_{n}^r g_{n} x_{n}  &= \la^{-2} e_n x_n^{r-1}x_{n+1}\inv g_n x_n\\
&=  \la^{-2} e_n  x_n^{r-1} g_n\inv x_n\inv x_n = 
 \la^{-2} e_{n} x_{n}^{r-1} g_{n}\inv .
\end{aligned}
$$
by Proposition \ref{proposition-xr relations2}.  Statement (4) follows by 
applying the isomorphism $\beta$ to~(3).
\end{proof}

\begin{corollary}\label{corollary -- trace of neg. powers of x1}
For $r \ge 1$,  $e_1 x_1^{-r} e_1 =  f_r(q_1, \dots, q_r) e_1$, where 
$f_r$ is the polynomial of  equation \ref{theta minus r polynomial in thetas} in Section
\ref{subsection- affine KT 
algebra}.
\end{corollary}

\begin{proof}   For $r \ge 1$,  set  $\psi_r = e_1 x_1^r e_1 = q_r e_1$, and $\psi_{-r} = e_1 x_1^{-r} e_1$.   For $a \ge 1$ and $b \ge 0$,  define
$\psi_{a, b} = e_1 x_1^{-b} g_1 x_1^{a} e_1$, and 
$\psi_{a, b}^{-} = e_1 x_1^{-b} g_1\inv x_1^{a} e_1$.   In particular, 
$\psi_{a, 0} = e_1 g_1 x_1^a e_1 = \la\inv q_a e_1$.
We have
\begin{equation} 
\label{equation1 -- corollary -- trace of neg. powers of x1}
\psi_{-r} = e_1 x_1^{-r} e_1 = \la\inv e_1 x_1^{-r} g_1\inv e_1 = \la \psi_{1, r-1},
\end{equation}
and
\begin{equation}
 \label{equation2 -- corollary -- trace of neg. powers of x1}
\psi_{a, b}^{-} = \la^2 \psi_{a+1, b-1},
\end{equation}
using Proposition \ref{proposition- affine braid relation 1}(4).
Moreover, for $a, b \ge 1$,
\begin{equation}
\label{equation3 -- corollary -- trace of neg. powers of x1}
\begin{aligned}
\psi_{a, b} &=   e_1 x_1^{-b} g_1 x_1^{a} e_1 \\
&=  e_1 x_1^{-b} g_1\inv x_1^{a} e_1 + z(e_1 x_1^{-b} e_1 x_1^a e_1 - e_1 x_1^{a-b} e_1) \\
&= \psi_{a, b}^{-} + z( \psi_{-b} q_a - \psi_{a-b}) \\
&= \psi_{a+1, b-1} + z( \psi_{-b} q_a - \psi_{a-b}),
\end{aligned}
\end{equation}
by the Kauffman skein relation, and equation \ref{equation2 -- corollary -- trace of neg. powers of x1}. 
Comparing the recursion relations of equations  \ref{equation1 -- corollary -- trace of neg. powers of x1} and  \ref{equation3 -- corollary -- trace of neg. powers of x1} with those of
Lemma \ref{lemma - recursion for f_r} gives that $\psi_{-r} = f_r(q_1, \dots, q_r) e_1$.
\end{proof}

\subsection{$\abmw  {n}$ as a 
$\abmw{n-1}$-bimodule.}\label{subsection-bimodule}
$\bmw {n, R}$ has a simple structure as a\break $\bmw {n-1, R}$-bimodule, 
described by the following
proposition, which is Lemma 3.1 from~\cite{Birman-Wenzl}.

\begin{proposition}\label{proposition-Birman-Wenzl lemma}
$\bmw {n, R}$ is the span of elements of the form $a \chi b$, where $a, b \in 
\bmw {n-1, R}$ and
$\chi \in \{ g_{n-1}^{\pm 1}, e_{n-1}, 1\}$.
\end{proposition}

We generalize this result to the affine Birman--Wenzl--Murakami  algebras in 
this section.
Fix $S$ and let $\abmw n$ denote $\abmw {n, S}$.

\begin{lemma}\label{lemma- bimodule 1}
\mbox{}
\begin{enumerate}
\item For $n \ge 1$ and $s  \in \Z$, there exists an element 
$b \in \abmw  {n-1}$ such that
$$
\qquad\ e_{n} x_{n}^{s} e_{n} = b e_{n}.
$$
\item For $n \ge 2$, 
let $a$ be an element of $\abmw  {n}$ of the form
$$
\qquad\ a = z_{0} x_{n-1}^{r_{1}} z_{1} x_{n-1}^{r_{2}}\cdots z_{k-1} x_{n-
1}^{r_{k}} 
z_{k},
$$
where $z_{i} \in\{e_{n-1}, g_{n-1}\}$ for all $i$, and $r_{i}\ge1$.
There exists $b \in \abmw  {n-1}$
such that 
$$
\qquad\ e_{n} a e_{n} = b e_{n}.
$$
\end{enumerate}
\end{lemma}

\begin{proof}  

For $n = 1$, the first assertion follows from the relation (5 c) of Definition
\ref{definition affine BMW} and Corollary \ref{corollary -- trace of neg. powers of x1}.
  We take this as the base for an 
inductive proof
of both assertions.

For $n \ge 2$ and $s \ge 1$,  the first assertion follows from the second, since
$x_{n}^{s} = (g_{n-1} x_{n-1} g_{n-1})^{s}$.  The assertion for $s \le -1$ 
follows by applying the
symmetry~$\beta$.

The proof of statement (2) is by induction first  on $n$ and then on the number 
$k$  of factors
$x_{n-1}^{ r_{j}} z_{j}$  in the expression for $a$.   For $n \ge 2$ and $k = 
0$,
the statement is immediate from the BMW relations.  Fix $n \ge 2$, $k \ge 1$, 
and $a$ of the form above with $k$ factors; assume inductively that

\vskip4pt
(1) $e_{m} x_{m}^{s} e_{m} \in \abmw {{m-1}} e_{m}$ for $m < n$ and $s \ge 1$, 
and

(2) $e_{n} a' e_{n} \in \abmw {{n-1}} e_{n}$ if $a'$ is of the same form with 
fewer than $k$ factors.
\vskip4pt

We consider several cases:

\ods
(a) Two successive $z_{i}$ are equal to $e_{n-1}$.  Then 
$$
a = a' e_{n-1} x_{n-1}^{s} e_{n-1} a'' = a' b e_{n-1} a'',
$$
where $b \in \abmw  {n-2}$, by the induction hypothesis.  Thus
$$
e_{n}a e_{n} = b e_{n} a' e_{n-1} a'' e_{n} \in \abmw  {n-1} e_{n},
$$
since $a' e_{n-1} a''$ has fewer than $k$ factors.

\ods
(b)  For some $i$, $z_{i} = e_{n-1}$ and  $z_{i  + 1} = g_{n-1}$.  Then
$$
\begin{aligned}
e_{n}a e_{n} &= e_{n}a' e_{n-1} x_{n-1}^{r} g_{n-1} a''e_{n}  \\[2pt]
&=\la^{-2} e_{n}a' e_{n-1} x_{n-1}^{r-1} g_{n-1}\inv x_{n-1}\inv a''e_{n}
\quad\ \text{(by Proposition \ref{proposition- unwrapping 1})} \\[2pt]
&=  \la^{-2}  e_{n}a' e_{n-1} x_{n-1}^{r-1} [ g_{n-1} + ze_{n-1} -z ] x_{n-
1}\inv a'' e_{n} \\[2pt]
& \equiv  \la^{-2} e_{n}a' e_{n-1} x_{n-1}^{r-1} g_{n-1} x_{n-1}\inv a'' e_{n}  
\bmod  \abmw  {n-1} e_{n} ,
\end{aligned}
$$
using  case (a) and the induction hypothesis.  If $a'' = 1$, then the final 
expression is equal to
$$
 \la^{-2} e_{n}a' e_{n-1} x_{n-1}^{r-1} g_{n-1}  e_{n} x_{n-1}\inv.
$$
Repeating this step $r$ times in all, we get
$$
\begin{aligned}
e_{n}a e_{n}  &\equiv \la^{-2 r}  e_{n}a' e_{n-1}  g_{n-1}  e_{n} x_{n-1}^{-r}    
\bmod  \abmw  {n-1} e_{n}\\[2pt]
& = \la^{-2r-1} e_{n}a' e_{n-1}   e_{n}  x_{n-1}^{-r} .
\end{aligned} 
$$
By the induction hypothesis, this is in $\abmw  {n-1} e_{n}$.  If $a'' \ne 1$, 
then
$a'' = x_{n-1}^{s} z_{i+2} a''' $, so we get
$$
e_{n}a e_{n}  \equiv  \la^{-2} e_{n}a' e_{n-1} x_{n-1}^{r-1} g_{n-1} x_{n-1}^{s-
1} z_{i+2} a''' e_{n}
\bmod \abmw  {n-1} e_{n}.
$$
If $r  > s$, we repeat this step $s$ times in all, obtaining finally
$$
e_{n}a e_{n}  \equiv  \la^{-2s } e_{n}a' e_{n-1} x_{n-1}^{r-s}g_{n-1}  z_{i+2} 
a''' e_{n}
\bmod \abmw  {n-1} e_{n}.
$$
Now $g_{n-1} z_{i+2}$ is a linear combination of 
$e_{n-1}$,  $g_{n-1}$, and $1$,
so by the induction hypothesis, the latter expression is
 in $\abmw  {n-1} e_{n}$.   The cases
$r < s$ and $r = s$ are similar.

\ods
(c)  For some $i$, $z_{i} = g_{n-1}$ and $z_{i+1} = e_{n-1}$.  This is 
essentially the same as case~(b).

\ods
(d)  For all $i$, $z_{i} = g_{n-1}$.   
If $a = g_{n-1} x_{n-1}^{r} g_{n-1}$, then
$$
\begin{aligned}
e_{n}a e_{n} &= e_{n }g_{n-1} x_{n-1}^{r} g_{n-1} e_{n}   \\
& \equiv e_{n }g_{n-1} x_{n-1}^{r} g_{n-1}\inv e_{n}  \bmod \abmw  {n-1} 
e_{n} \quad\ \text{(as in case (b))}\\
&= e_{n} g_{n-1} g_{n} x_{n-1}^{r} g_{n} \inv g_{n-1}\inv e_{n}  \\
&= e_{n} e_{n-1} x_{n-1}^{r} e_{n-1} e_{n}  \in  \abmw  {n-1} e_{n}, \\
\end{aligned}
$$
by case (a).  Otherwise, $a = g_{n-1} x_{n-1}^{r} g_{n-1} x_{n-1}^{s} g_{n-1} a' 
$, and
$$
\begin{aligned}
e_{n}a e_{n} &=  e_{n}  g_{n-1} x_{n-1}^{r} g_{n-1} x_{n-1}^{s} g_{n-1} a'    
e_{n} \\
& \equiv  e_{n}  g_{n-1}\inv  x_{n-1}^{r} g_{n-1} x_{n-1}^{s} g_{n-1} a'    
e_{n} \bmod \abmw  {n-1} e_{n} \quad\ \text{(as in case (b))}.
\end{aligned}
$$
The affine braid relation of Proposition \ref{proposition- unwrapping 1} implies
$$
x_{n-1}^{r} g_{n-1} x_{n-1} g_{n-1} = g_{n-1} x_{n-1} g_{n-1} x_{n-1}^{r}  , 
$$
or
$$
g_{n-1}\inv x_{n-1}^{r} g_{n-1} x_{n-1}  
= x_{n-1} g_{n-1} x_{n-1}^{r}  g_{n-1} \inv.
$$
Applying this, we get
$$
e_{n}a e_{n}   \equiv  e_{n}x_{n-1}  g_{n-1}  x_{n-1}^{r}   g_{n-1}\inv 
 x_{n-1}^{s-1} g_{n-1} a'    e_{n} \bmod \abmw  {n-1} e_{n}.
$$
Now we can change the $ g_{n-1}\inv$ to $g_{n-1}$ while maintaining 
congruence mod
$\abmw  {n-1} e_{n} $,   as in case (b), so
$$
e_{n}a e_{n}  \equiv x_{n-1}  e_{n} g_{n-1}  x_{n-1}^{r}   g_{n-1} 
 x_{n-1}^{s-1} g_{n-1} a'    e_{n} \bmod \abmw  {n-1} e_{n}.
$$
Repeating this step a total of $s$ times, we get
$$
e_{n}a e_{n}  \equiv x_{n-1}^{s}  e_{n} g_{n-1}  x_{n-1}^{r}   g_{n-1}^{2} a'    
e_{n} \bmod \abmw  {n-1} e_{n}.
$$
Since $g_{n-1}^{2}$ is a linear combination of $e_{n-1}$, $g_{n-1}$, and $1$,  
the last expression is in $\abmw  {n-1} e_{n}$ by the induction assumption.
\end{proof}

For the remainder of this section, we maintain the following notation:
 $A_{1} $~de\-notes the linear span of $x_{1}^{r}$  for $r \in \Z$.  (Thus 
$A_{1} = \abmw  {1}$.)
For each $n \ge 2$,  $A_{n}$ denotes the linear span of $e_{n-1}$, $g_{n-
1}\pmone$, and
$x_{n}^{r}$  for $r \in \Z$.    

The following proposition is the analogue for the $\abmw  {n}$  of Lemma 3.1 
in~\cite{Birman-Wenzl}.  
This result is due to H\"aring-Oldenburg~\cite{H-O2}.

\begin{proposition}\label{proposition- bimodule 2}
 Every element of $\abmw  {n}$ is a linear combination of
elements of the  form $a \chi b$, where $a, b \in \abmw  {n-1}$ and 
$\chi  \!\in\! \{e_{n-1}, g_{n-1}\pmone\}
\cup \{x_{n}^{r} : r \!\in\! \Z\}$.
\end{proposition}

\begin{proof} 
We have to show that for all $n \ge 1$, $\abmw  {n} = \abmw  {n-1} A_{n}  \abmw  
{n-1}$.
Since $\abmw  {n} $ is generated as an algebra by $\abmw  {n-1}$ and $A_{n}$, it  
suffices to show that $A_{n} \abmw  {n-1} A_{n} \subseteq \abmw  {n-1} A_{n} 
\abmw  {n-1}$.

The assertion is evident for $n = 1$.  We assume it holds for a
 particular $n \ge 1$ and
prove that $A_{n+1}  \abmw  {n} A_{n+1} \subseteq \abmw  {n} A_{n+1}
\abmw  {n} $.
By the induction assumption, $\abmw  {n} =  \abmw  {n-1} A_{n}  
\abmw  {n-1}$, so
$A_{n+1}  \abmw  {n} A_{n+1} \!=\! A_{n+1} \abmw  {n-1} A_{n}  \abmw 
 {n-1}  A_{n+1}\break
=  \abmw  {n-1}  A_{n+1}A_{n}    A_{n+1} \abmw  {n-1}$.  Thus it suffices to 
show that
$$
A_{n+1}A_{n}    A_{n+1} \subseteq \abmw  {n} A_{n+1} \abmw  {n}.
$$
We consider several cases:

\ods
(a) $\chi_{n}' \chi_{n-1} \chi_{n}'' \in \abmw  {n} A_{n+1} \abmw  {n}$, where
$\chi_{n}', \chi_{n}'' \in \{e_{n}, g_{n}\pmone\}$ and $\chi_{n-1} \in 
 \{e_{n-1}, g_{n-1}\pmone\}$.  This follows easily from the BMW relations.
 
\ods
(b) $e_{n} x_{n}^{r} e_{n} \in  \abmw  {n-1} e_{n} \subseteq \abmw  {n} 
A_{n+1}$.   This follows from
Lemma \ref{lemma- bimodule 1}.

\ods
(c)  $g_{n}\pmone x_{n}^{r} e_{n} \in   \abmw  {n} e_{n} \subseteq \abmw  {n} 
A_{n+1}$, and
$ e_{n}  x_{n}^{r} g_{n}\pmone \in  e_{n}\abmw  {n}  \subseteq A_{n+1}\abmw  
{n}$.  The second statement follows from the first by applying the symmetry 
$\alpha$.
   Note that
$$
g_{n}\inv x_{n}^{r} e_{n}  \equiv g_{n} x_{n}^{r} e_{n}  
\bmod  \abmw  {n} e_{n},
$$
by the Kauffman skein relation and case (b).  Moreover, the assertion for $r \le 
-1$ follows from the assertion for $r \ge 1$ by applying the symmetry $\beta$.  
For $r \ge 1$, we have
$$
g_{n} x_{n}^{r} e_{n} = \la^{-2} x_{n}\inv g_{n}\inv  x_{n}^{r-1} e_{n}  \equiv
x_{n}\inv g_{n} x_{n}^{r-1} e_{n} \bmod \abmw  {n} e_{n},
$$
by the unwrapping relation of Proposition \ref{proposition- unwrapping 1}.
Repeating this step a total of $r$ times gives
$$
\begin{aligned}
g_{n} x_{n}^{r} e_{n}  &\equiv  \la^{-2r}x_{n}^{-r} g_{n} e_{n}  \bmod  \abmw  
{n} e_{n}\\
&= \la^{-2r -1}   x_{n}^{-r} e_{n}. \\
\end{aligned}
$$

(d)   $A_{n+1}  A_{n} x_{n+1}^{t } \subseteq  \abmw  {n} A_{n+1} \abmw  {n}$ and
 $ x_{n+1}^{t } A_{n} A_{n+1} \subseteq  \abmw  {n} A_{n+1} \abmw  {n}$.  The 
second assertion follows from the first by applying the symmetry $\alpha$.  
Since $x_{n+1}^{t } $ commutes with $A_{n}$, it suffices to show $A_{n+1}  
x_{n+1}^{t } \subseteq  \abmw  {n} A_{n+1} \abmw  {n}$.   Since $x_{n+1}^{s} 
x_{n+1}^{t} = x_{n+1}^{s+t} \in
 \abmw  {n} A_{n+1} \abmw  {n}$, we only have to check that
 $
 e_{n}  x_{n+1}^{t } \in  \abmw  {n} A_{n+1} \abmw  {n},
 $
 and
  $
 g_{n}\pmone  x_{n+1}^{t } \in  \abmw  {n} A_{n+1} \abmw  {n}.
 $
 
 We have
 $$
 e_{n}  x_{n+1}^{t } = \la^{-2t} e_{n} x_{n}^{-t} \in  \abmw  {n} A_{n+1} \abmw  
{n},
 $$
 by Proposition \ref{proposition-xr relations2}.
 It follows from this and the Kauffman skein relation that
 $$
 g_{n}  x_{n+1}^{t }  \equiv g_{n}\inv  x_{n+1}^{t }  \bmod \abmw  {n} A_{n+1} 
\abmw  {n}.
 $$
 Moreover, for $t \ge 1$,
 $$
 g_{n}\inv x_{n+1}^{t } = x_{n} g_{n} x_{n+1}^{t-1 } \equiv x_{n} g_{n}\inv 
x_{n+1}^{t-1 } 
 \bmod \abmw  {n} A_{n+1} \abmw  {n}.
 $$
 Repeating this step $t$ times gives
 $$
 g_{n}\inv x_{n+1}^{t } \equiv  x_{n}^{t} g_{n}\inv \bmod \abmw  {n} A_{n+1} 
\abmw  {n}.
 $$
 Thus for $t \ge 1$,
 $$
 g_{n}\pmone x_{n+1}^{t } \in  \abmw  {n} A_{n+1} \abmw  {n},
 $$
and the same statement for $t \le -1$ follows by applying
 the symmetry $\beta$.
\end{proof}

\begin{proposition}\label{proposition- bimodule 3}
 For $n\ge 1$,
 $e_{n} \abmw  {n} e_{n}  =  \abmw  {n-1} e_{n}.$
\end{proposition}

\begin{proof}  We only have to prove the containment
$e_{n} \abmw  {n} e_{n}  \subseteq  \abmw  {n-1} e_{n}$.
This is obvious for $n = 1$.  For $n \ge 2$, we have
$$
e_{n} \abmw  {n} e_{n}  = e_{n} \abmw  {n-1} A_{n} \abmw  {n-1} e_{n} =
\abmw  {n-1}  e_{n}A_{n} e_{n }\abmw  {n-1}.
$$
 But $e_{n} A_{n} e_{n} \subseteq \abmw {n-1} e_{n}$ by the BMW relations and 
Lemma \ref{lemma- bimodule 1}.
\end{proof}

\begin{corollary}\label{corollary- bimodule 3.5}
 For $n\ge 1$,
 $e_{n} \bmw {n} e_{n}  =  \bmw {n-1} e_{n}.$
\end{corollary}

\begin{lemma}\label{lemma- bimodule 4}
For $n \ge 1$,
$
g_{n}\pmone \abmw  {n} e_{n} \subseteq  \abmw  {n} e_{n}
$
and
$
x_{n+1}^{r} \abmw  {n} e_{n} \subseteq  \abmw  {n} e_{n}.
$
\end{lemma}

\begin{proof}
Let $\chi$ be one of $x_{n+1}^r$ or $g_{n}\pmone$.  We have
$$
\chi \abmw  {n} e_{n} =  \chi  \abmw  {n-1}A_{n} \abmw  {n-1} e_{n} =
 \abmw  {n-1}\chi  A_{n} e_{n}\abmw  {n-1}.
$$
But $g_n\pmone  A_{n} e_{n} \subseteq \abmw  {n} e_{n}$ by the BMW relations and 
case  (c) in the
proof of Proposition \ref{proposition- bimodule 2}.  Likewise, 
$x_{n+1}^r A_n e_n = A_n x_{n+1}^r e_n = \la^{-2r} x_n^{-r} e_n \subseteq \abmw 
n e_n$.
Thus
$
\chi \abmw  {n} e_{n}  \subseteq  \abmw  {n} e_{n}\abmw  {n-1} = \abmw  {n} 
e_{n}.
$
\end{proof}

\begin{proposition}\label{proposition- bimodule 5}
For $n \ge 1$, $\abmw  {n+1} e_{n} = \abmw  {n} e_{n}$.
\end{proposition}

\begin{proof}
We have
$
\abmw  {n+1} e_{n} = \abmw  {n} A_{n+1} \abmw  {n} e_{n} \subseteq \abmw  {n} 
e_{n},
$
by Proposition \ref{proposition- bimodule 3}  and Lemma \ref{lemma- bimodule 4}.
\end{proof}

\begin{corollary}\label{corollary- bimodule 5.5}
For $n \ge 1$, $\bmw {n+1} e_{n} = \bmw {n} e_{n}$.
\end{corollary}

Note that Corollaries \ref{corollary- bimodule 3.5} and  
\ref{corollary- bimodule 5.5} also follow directly from Proposition
\ref{proposition-Birman-Wenzl lemma} by similar reasoning.

\subsection{A homomorphism from $\abmw {n} $ to 
$\akt {n}$.}\label{section-homomorphism phi}
 Let  $R$ be a ring with distinguished
elements $\lambda$,
$z$, and $\delta$ as above, and $S \supseteq R$  a ring with additional
elements $q_1, q_2, \dots$.

Let $X_1$, $G_i$, and $E_i$ denote the affine tangle diagrams shown in the 
figure following Definition
\ref{definition affine tangle},
 regarded as elements of the affine Kauffman
tangle algebra $\akt {n, S}$.

\begin{proposition}\label{proposition-algebra homomorphism from
affine bmw to affine tangle}
The assignment  $x_1 \mapsto X_1$, 
$g_i \mapsto G_i$, and  
$e_i \mapsto E_i$ 
determines an $S$-algebra homomorphism $\varphi$ from 
$\abmw {n, S } $ to $\akt {n, S}$.
\end{proposition}

{\em Proof}.\hods  One needs to check that $X_1$, $G_i$, and $E_i$  satisfy
the relations of the affine BMW algebra 
(Definition~\ref{definition affine BMW}), using regular isotopy of affine
tangles and the relations of the affine Kauffman tangle algebra.  The
inverses of $G_i$ and $X_1$ in $\akt {n, S}$ are represented by the tangle
diagrams with the crossings reversed:
$$
X_1\inv = \inlinegraphic{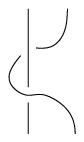}\quad\ 
G_i\inv = \inlinegraphic{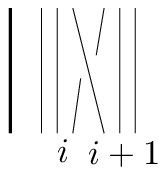}.
$$
The verification that these elements are in fact inverses involves
applications of Reidemeister II.

The Kauffman skein relation for $G_i$ depends on 
point (1) of Definition \ref{definition Kauffman
tangle algebra}.
The relation $E_i^2 = \delta E_i$  results from point (3) of 
Definition \ref{definition Kauffman
tangle algebra}.  The relation $E_1 X_1^{ r} E_1 = q_r E_1$ 
results from 
point (2) of Definition \ref{definition affine Kauffman
tangle algebra}.

  The affine braid relations result from applications of
Reidemeister moves II and III.  The remaining commutation relations and the
tangle relation (6a) only use planar isotopy, while (6b) also uses
Reidemeister move II.  The untwisting relations (7) follow from the
untwisting relation (2) for $\akt {n, S}$.

Finally, the unwrapping relation follows from the following
graphical ``computation'', which we give here for the sake of
illustration.  The first equality of pictures is just isotopy, the second 
a Reidemeister III move, the third an untwisting move, and the last two
are Reidemeister II moves.
$$
\begin{aligned}
E_1 X_1 G_1 X_1 &= \inlinegraphic{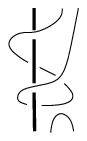} =
\inlinegraphic{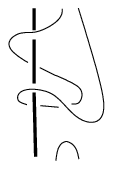} \\ 
&= \inlinegraphic{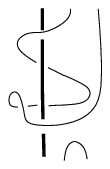} = \la \inv \inlinegraphic{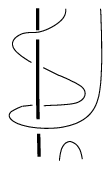} \\
&= \la \inv \inlinegraphic{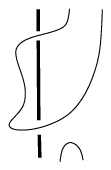} = \la \inv
\inlinegraphic{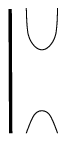} = \la \inv E_1.\squ
\end{aligned}
$$

\begin{proposition}\label{proposition-algebra homomorphism from bmw
to tangle}
 The assignment $g_i \mapsto G_i$ and $e_i \mapsto E_i$
determines an $R$-algebra homomorphism $\varphi$ from $\bmw {n, R}$ to $\kt {n, 
R}$.
\end{proposition}

\begin{proof}  One has to verify relations, as in the previous proof.
\end{proof}

\begin{remark}\rm
\mbox{} 
\begin{enumerate}
\item 
The homomorphism $\varphi: \abmw {n }  \rightarrow \akt {n}$
gives us a trace $\eps : \abmw {n }  \rightarrow S$ defined by
$\eps = \eps\circ\varphi$.
\item  The following diagram commutes (see Propositions 
\ref{proposition-algebra homomorphism from
affine bmw to affine tangle} and
\ref{proposition-algebra homomorphism from
bmw to tangle}   and  Remarks~\ref{remark-injections of tangle algebras in 
affine tangle algebras} and
\ref{remark- inclusions and symmetries affine BMW}):
$$\qquad\
\begin{diagram}
\abmw {n }      &\rTo^{\scriptstyle \varphi}    & \akt {n} \\
\uTo<{\scriptstyle i_n}                  &        & \uTo<{\scriptstyle i_n}  \\
\bmw {n}                 &\rTo^{\scriptstyle \varphi}    & \kt {n} \\
\end{diagram}
$$
 \item The following diagrams commute:
 $$\qquad\
 \begin{diagram}
\bmw {n}     &\rTo^{\scriptstyle\varphi}    & \kt {n} \\
\dTo<{\scriptstyle\iota}                  &        & \dTo<{\scriptstyle\iota}  
\\
\bmw {n+1}               &\rTo^{\scriptstyle\varphi}    & \kt {n+1} \\
\end{diagram}
\qquad
 \begin{diagram}
\abmw {n }      &\rTo^{\scriptstyle\varphi}    &\akt {n} \\
\dTo<{\scriptstyle\iota}                  &        & \dTo<{\scriptstyle\iota}  
\\
\abmw  {n+1}                 &\rTo^{\scriptstyle\varphi}    & \akt {n+1}\\
\end{diagram}
 $$
\item The following diagrams commute (see Section \ref{subsection-symmetries of 
kauffman tangle algebra},  
Remark  \ref{remark- injection-  shift- and reversal},  and Remark 
\ref{remark- inclusions and symmetries affine BMW}):
$$
\qquad\ \begin{diagram}
\abmw {n }      &\rTo^{\scriptstyle\varphi}    & \akt {n} \\
\dTo<{\scriptstyle\alpha}                  &        & \dTo<{\scriptstyle\alpha}  
\\
\abmw {n}                 &\rTo^{\scriptstyle\varphi}    & \akt {n} \\
\end{diagram}
\quad\
\begin{diagram}
\abmw {n }      &\rTo^{\scriptstyle\varphi}    & \akt {n} \\
\dTo<{\scriptstyle\beta}                  &        & \dTo<{\scriptstyle\beta}  
\\
\abmw {n}                 &\rTo^{\scriptstyle\varphi}    & \akt {n} \\
\end{diagram}
\quad\
\begin{diagram}
\bmw {n}     &\rTo^{\scriptstyle\varphi}    & \kt {n} \\
\dTo<{\scriptstyle\varrho_n}                  &        & 
\dTo<{\scriptstyle\varrho_n}  \\
\bmw {n}                 &\rTo^{\scriptstyle\varphi}    & \kt {n} \\
\end{diagram}\vadjust{\vskip4pt}
$$
\item  One has $\varphi(x_{r}) = X_{r}$  for all $r$. 
\end{enumerate}
\end{remark}

We will eventually show that $\varphi: \abmw {n, S } 
\rightarrow \akt {n, S}$ is an isomorphism for any $S$. 

\section{The affine Brauer algebra}\label{sec4}
\vskip0pt\vskip-10pt\vskip0pt
\subsection{The Brauer algebra.}
The Brauer algebra $\br_n$ is an algebra of planar $(n,n)$-tangle diagrams
in which {\em crossings are ignored}.  The precise definition follows.

Fix points $a_i $ in $I$, for $i \ge 0$, as in the
description of $(n,n)$-tangles.  For convenience write
$\p i = (a_i, 1)$ and $\overline{ \p i} = (a_i, 0)$.

\begin{definition}\rm\label{definition brauer diagram}
An $(n,n)$-{\em Brauer diagram}  (also called an $n$-{\em connector})
consists
of a collection of 
$n$ curves in the rectangle $R = I \times I$ such that
\begin{enumerate}
\item  The curves connect the points $\{\p 1,  \dots, \p n, \pbar 1,
\dots \pbar n\}$ in pairs.
\item  For each curve $C$ in the collection, the intersection of $C$ with
$\bdry(R)$ consists of the two endpoints of $C$.
\end{enumerate}
\end{definition} 

Consider the free $\Z[\deltabold^{\pm 1}]$-module $\br_n$ with basis the set
of 
$(n,n)$-Brauer diagrams.   The product of two Brauer diagrams is defined
to be a certain multiple of another Brauer diagram.  Namely, given two
Brauer diagrams $a, b$, first ``stack" $b$ over $a$ (as for tangle
diagrams).   Let $r$ denote the number of closed curves in the interior
of $R$ in the resulting planar ``tangle", and let $c$ be the Brauer
diagram obtained by removing all the closed curves.  Then
$$
a b = \deltabold^r c.
$$

\begin{definition}\rm
The {\em Brauer algebra} $\br_n$ over $\Z[\deltabold^{\pm
1}]$ is the free $\Z[\deltabold^{\pm 1}]$-module with basis the set of 
$(n,n)$-Brauer diagrams, with the bilinear product determined by the
multiplication of Brauer diagrams.
\end{definition}

The Brauer algebras were introduced by R. Brauer~\cite{Brauer} as a device
for studying the invariant theory of orthogonal and symplectic groups.
The generic structure of the Brauer algebras, and conditions for
semisimplicity of $\br_n \otimes_{\Zdeltabold} k$, where $k$ is a field,
were determined by H.  Wenzl~\cite{Wenzl-Brauer}.

Note that the Brauer diagrams with only vertical strands, that is,
diagrams in which upper points are paired only with lower points, are in
bijection with permutations of $\{1, \dots, n\}$, and that the
multiplication of two such diagrams coincides with the multiplication of
permutations.  Thus the  Brauer algebra contains the group algebra of
the permutation group $\S_n$ (over~$\Zdeltabold$).

\subsection{The affine Brauer algebra.}
We will define  the affine Brauer algebra as a sort of  ``wreath product''
of $\Z$ with the Brauer algebra (containing the wreath product of $\Z$
with~$\S_n$).

\begin{definition}\rm
A {\em colored $(n,n)$-Brauer diagram}, or  {\em colored
$n$-con\-nector} is a Brauer diagram in which each strand is labeled by an
integer.
\end{definition}

We will define the affine Brauer
algebra over the ring
$$
\Z[\deltabold^{\pm 1}]^\wedge= \Z[\deltabold^{\pm 1}, \qbold_1, 
\qbold_2,
\dots],
$$
where $\qbold_1, \qbold_2, \dots $ are indeterminates. 

Order the points $\{\p 1, \dots, \p n, \, \pbar 1, \dots, \pbar n\}$ by
$\p 1 < \cdots < \p n < \pbar n  < \cdots < \pbar 1$. 
The colors on the strands of a colored Brauer diagram should be regarded
as assigning integer values to {\em oriented strands} of the Brauer
diagram.  If a strand is colored  by the integer $r$, then the strand
endowed  with the positive orientation (i.e. the orientation from a lower
numbered vertex to a higher numbered vertex) takes the value $r$, but the
same strand endowed with the negative orientation takes the value~$-r$. 

Consider the free $\zdeltahat$-module $\abr_{n}$ with basis the set
of  colored
$(n,n)$-Brauer diagrams.   The product of two  colored Brauer diagrams
is defined to be a certain multiple of another  colored  Brauer
diagram, determined as follows.

Given two colored Brauer diagrams $a, b$, first ``stack" $b$ over $a$ (as
for tangle diagrams and ordinary Brauer diagrams).   In the resulting
``tangle'' there are three types of curves:\vskip6pt\vskip0pt
\begin{enumerate}

\item  {\em Vertical strands}.\hods  These are concatenations of one vertical 
strand
from $b$, an even number of horizontal strands from the bottom of $b$
and the top of $a$, and finally one vertical strand from $a$.  Traversing
such a composite strand in its standard orientation (from lower numbered
vertex to higher numbered vertex) sum the integer values of the {\em
oriented} strands encountered.  (We repeat for emphasis:  if a strand
colored by
$r$ is traversed in the negative direction, then it contributes $-r$ to
the sum.)  Color the strand by the resulting sum.

\vskip6pt
\item  {\em Horizontal strands}.\hods These include  horizontal strands 
remaining from the original diagrams, namely horizontal strands from the
top of~$b$, and from the bottom of~$a$.  These strands retain their
original coloring from the original diagrams.  The horizontal strands
also include concatenations of a vertical strand from $b$ or $a$ followed
by an odd number of horizontal strands from the bottom of $b$ and the
top of $a$, and finally a second vertical  strand from the same diagram
as the first vertical strand.  The color of the composite strand is
determined as for composite vertical strands.

\vskip6pt
\item {\em Closed strands}.\hods  These are concatenations of an even number of
horizontal strands from the bottom of $b$
and the top of $a$.  There is no preferred orientation on such a strand,
so pick an orientation arbitrarily, and obtain a color by summing the
integer values of the oriented strands encountered in traversing the
curve, as for vertical strands.   For each $i \in \{0, 1, \dots\}$, let
$m_i$ be the number of closed loops with color~$\pm i$.\vskip6pt\vskip0pt
\end{enumerate}

Let $c$ be the colored 
Brauer diagram obtained by removing all the closed curves.  Then
$$
a b = (\deltabold^{m_0} \qbold_1^{m_1} \qbold_2^{m_2} \cdots) c.
$$

\begin{definition}\rm\label{definition affine Brauer algebra}
 The {\em affine Brauer algebra}
$\abr_n$ over
$\zdeltahat$ is the free $\zdeltahat$-module with basis
the set of colored
$(n,n)$-Brauer diagrams, with the bilinear product determined by the
multiplication of colored Brauer diagrams.
\end{definition}

One can easily check that the multiplication is associative.  Note that
the subalgebra generated by colored Brauer diagrams with only vertical
strands is isomorphic to the wreath product of $\Z$ with~$\S_n$.

\subsection{Conditional expectation and trace for the affine Brauer
al\-gebras.}
Just as for affine Kauffman tangle algebras, one has a homomorphism
$\iota$ of 
$\abr_{n-1}$ into $\abr_n$ by attaching an additional strand on the right
of an affine Brauer diagram (colored by~$0$).  

Moreover, one has a conditional expectation $\eps_n : \abr_n \rightarrow
\abr_{n-1}$ defined as follows:  First define a map $\cl_n$ from
colored $n$-connectors to colored $(n-1)$-connectors by
 joining the rightmost
pair of vertices  $\pbar n, \p n$ of a colored $n$-connector  $d$ by a
new strand, with color~$0$:
$$
{\rm cl}_n: \quad \inlinegraphic{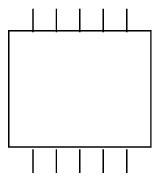} \quad \mapsto \quad 
\inlinegraphic{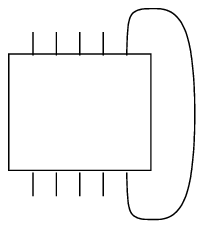}
$$

The new  strand is part of a concatenated (vertical, horizontal or closed)
strand.  If the concatenated strand is vertical or horizontal, orient it
positively (from lower numbered vertex to higher) and obtain its color by
adding the colors of its (oriented) components.  If the concatenated strand
is closed (which happens precisely if $d$ contains a strand connecting
$\p n$ and
$\pbar n$ with some color $r$) then remove the closed loop and multiply
the resulting colored
$(n-1)$-connector by $\qbold_{|r|}$ if $r \ne 0$ or $\deltabold$ if $r = 0$. 
For
example:
$$
{\rm cl}_n: \quad \inlinegraphicee{brauer_example} \quad \mapsto \quad 
\inlinegraphicee{brauer_example_closure} \quad \mapsto \quad 
\inlinegraphicee{brauer_example2}
$$
Define $\eps_n : \abr_n \rightarrow \abr_{n-1} $ by 
$$
\eps_n(d) =\deltabold\inv \cl_n(d).
$$

One can check that $\eps_n$ is a conditional expectation.  Since
$\eps_n \circ \iota(x) = x$,  the map $\iota$ is injective;  therefore, we
consider $\abr_{n-1}$ as a subalgebra of $\abr_n$.

Define $\eps = \eps_1 \circ \dots \circ \eps_n : \abr_n \rightarrow
\abr_0 = \zdeltahat$.  Alternatively, define the closure $\cl$ of a
colored $n$-connector by closing all the strands:
$$
{\rm cl}: \quad \inlinegraphic{tangle_box2} \quad \mapsto \quad 
\inlinegraphic{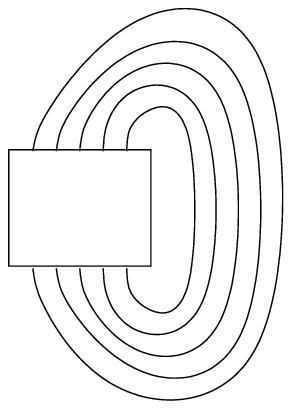}
$$
Compute the
color of each (closed) strand in the resulting diagram as before, and
replace each  closed strand by the appropriate factor
$\qbold_r$ or $\deltabold$.  Then $\eps(d) = \deltabold^{-n} \cl(d)$.

Using this picture for $\eps$, one can check that $\eps$ is a trace.

We can define a symmetric $\zdeltahat$-bilinear form on $\abr_n$ by
$(x, y) \mapsto \eps(xy)$.   

\def\ov#1{#1\llap{$\overline{\phantom{\rm #1}}$}}
\def\ovv#1#2{#1\llap{$\overline{\phantom{\rm #2}}$}}

Define the reflection $\ovv{d}{\hskip.7pt I\hskip-1pt}$ of a colored 
$n$-connector $d$  by reflecting the diagram vertically and changing the
color of each  strand to its opposite:
$$
 \inlinegraphic{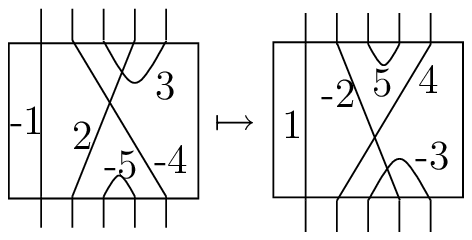}
$$
Note that for colored $n$-connectors $d, d'$, the closure
$\cl(d d')$ has at most $n$ (closed) strands, and that it has $n$
precisely when the underlying diagrams of $d$ and $d'$ are reflections of
one another.   Moreover, each of these closed loops has color $0$
precisely when $d' = \ovv{d}{\hskip.7pt I\hskip-1pt}$.  Consequently, we have:

\begin{lemma}\label{lemma-det A not zero}
\mbox{} 
\begin{enumerate}
\item For colored $n$-connectors $d, d'$, 
$$\qquad\ 
\eps(d d') =\deltabold^{-n} m,
$$
where $m$ is a  monomial in $\deltabold, \qbold_1, 
\qbold_2,
\dots$ of total degree  $\le n$. The degree of $\eps(d d')$  in $\deltabold$
is strictly negative unless $d' = \ovv{d}{\hskip.7pt I\hskip-1pt}$, while
$\eps(d \ovv{d}{\hskip.7pt I\hskip-1pt})  = 1$.
\item Let $S$ be a finite set of colored
$n$-connectors that is closed under the reflection $d \mapsto \ovv{d}{\hskip.7pt 
I\hskip-1pt}$. 
Consider the matrix
$A_S = (\eps(d d'))_{d, d' \in S}$.  The determinant of $A_S$ is
non-zero.
\end{enumerate}
\end{lemma}

\begin{proof} The first statement follows from the preceding discussion.
For the second statement,  each
row and column of
$A_S$ has exactly one entry equal to~1.  All other entries have strictly
negative degree in
$\deltabold$.  Therefore, $\det (A_S) = \pm 1 + \deltabold\inv p$, where $p$ is a
polynomial in $\deltabold\inv$ with coefficients in $\Z[\qbold_1, \qbold_2, 
\dots]$. 
\end{proof}
  
\subsection{A homomorphism from $\akt {n, \varLambdaHat}$ to 
$\abr_n$.}
We wish to define
a map from affine tangle diagrams to colored Brauer diagrams which
basically forgets the sense of crossings of ordinary strands, but
remembers the sense of crossings of ordinary strands with the flagpole.

\def\cconn{{\rm cconn}}
Define a map $c$ (the connector map) from affine $(n,n)$-tangle diagrams to 
$\abr_n$ as follows.  
Number and order the $2n$
vertices of affine
$(n,n)$-tangle diagrams by the same convention as for colored Brauer diagrams.
For each strand of an affine $(n,n)$-tangle  diagram $a$ that connects two
vertices, draw a curve connecting the corresponding vertices in 
$c(a)$.  Determine the color $r$ of the curve as follows:  with
the strand oriented  from the lower numbered vertex to the
higher numbered vertex, $r$ is the number of clockwise rotations of the
strand around the flagpole, viewed from above~$(^1)$\footnot{$(^1)$ The 
color $r$
can be determined combinatorially:  traversing the strand from lower
numbered vertex to higher numbered vertex, list the over-crossings $(+)$ and
under-crossings $(-)$ of the strand with the flagpole.  Cancel any
two successive  $+$'s or $-$'s in the list, so  the list now consists of
alternating $+$'s and $-$'s.  Then $r$ is $\pm (1/2)$ the length of the
list, $+$ if the list begins with a $+$, and $-$ if the list begins with
a~$-$.}.

Let $d$ be the resulting colored $n$-connector.  
Give each closed strand of $a$ an arbitrary orientation, and determine the color of the strand as above;   for each $r \ge 1$,  let $m_r$ be the number of closed loops in $a$
with color $\pm r$.
Finally, define
$$
c(a) = (\deltabold^{m_0} \qbold_1^{m_1} \qbold_2^{m_2} \cdots) d.
$$

The map $c$ respects regular isotopy (in fact, ambient isotopy), so induces a 
map
$c: \uhat n n \rightarrow \abr_n$.  One can check that $c$ is a monoid map, so 
we extend it linearly to 
a $\zdeltahat$-algebra map 
$$
c: \zdeltahat \uhat n n \rightarrow \abr_n.
$$
  This 
algebra map respects the 
relations of $\akt {n, \zdeltahat}$  (see Remark \ref{remark - Brauer specialization of affine algebra} at the the end of Section
\ref{subsection- affine KT algebra}), so induces a map
$$
c: \akt {n, \zdeltahat} \rightarrow \abr_n.
$$
  Finally, we have the
composition
$$
c: \akt {n, \varLambdaHat}  \rightarrow \akt {n, \zdeltahat} \rightarrow 
\abr_n.
$$
This map $c$ is given simply by the formula $c(\sum \alpha_i T_i) = \sum 
e(\alpha_i) c(T_i)$, where
$\alpha_i \in \varLambdaHat$, the $T_i$ are affine tangles representing elements 
of $\akt {n, \varLambdaHat}$, and $e$ is the
homomorphism of $\varLambdaHat$ to $\zdeltahat$ determined by $\labold \mapsto 
1$, $\zbold \mapsto 0$, $\deltabold \mapsto \deltabold$,
and $\qbold_i \mapsto \qbold_i$.

\begin{proposition}
The following diagrams commute:
$$
\begin{diagram}
\akt {n, \varLambdaHat}&\rTo^{\scriptstyle c}    & \abr_n \\
\dTo<{\scriptstyle\eps_n}                  &        & \dTo<{\scriptstyle\eps_n}  
\\
\akt {n-1, \varLambdaHat}                 &\rTo^{\scriptstyle c} & \abr_{n-1}\\
\end{diagram}
\quad\ 
\begin{diagram}
\akt {n, \varLambdaHat}&\rTo^{\scriptstyle c}    & \abr_n \\
\dTo<{\scriptstyle\eps}                  &        & \dTo<{\scriptstyle\eps}  \\
\varLambdaHat                 &\rTo^{\scriptstyle e}    & \zdeltahat\\
\end{diagram}
$$
\end{proposition}

\begin{proof}  Left to the reader.
\end{proof}

The following proposition is adapted from~\cite{Morton-Traczyk}
and~\cite{Morton-Wassermann}.

\begin{proposition}\label{proposition- distinct affine connectors linear 
independence}
For each colored $n$-connector $d$, let $T_d$ be an affine
$(n,n)$-tangle diagram with $c(T_d) = d$.  Then the set $\{T_d : d$
 a colored $n$-connector$\}$ is linearly independent over
$\varLambdaHat$.
\end{proposition}

\begin{proof}  Let $S$ be a finite set of colored $n$-connectors
which is closed under reflection $d \mapsto \ovv{d}{\hskip.7pt I\hskip-1pt}$. 
Consider the matrices $B_S = (\eps(T_d T_{d'}))_{d, d' \in
S}$ and $A_S = (\eps(d d'))_{d, d' \in S}$. One has 
$$
e(\det(B_S)) = \det(e \circ \eps(T_d T_{d'})) = \det (\eps \circ c(T_d
T_{d'})) = \det(\eps( d d')) = \det(A_S).
$$
By Lemma \ref{lemma-det A not zero}, $\det(A_S) \ne 0$, so $\det(B_S) \ne
0$.  Since
$\varLambdaHat$ is an integral domain,  $B_S$ is invertible over the
field of fractions of $\varLambdaHat$. It follows that 
$\{T_d : d \in S\}$ is linearly
independent over $\varLambdaHat$.  

Since $S$ is arbitrary, it follows that $\{T_d : 
d \text{ a colored $n$-connector}\}$ is linearly independent over
$\varLambdaHat$.
\end{proof}

\begin{corollary}\label{corollary-independence of tangles with different 
connectors}
For each ordinary $n$-connector $d$, let  $T_d$ be an ordinary $(n,n)$-tangle 
diagram with $c(T_d) = d$.  Then the set $\{T_d : d
\text{ an $n$-connector}\}$ is linearly independent over~$\varLambda$.
\end{corollary}

\begin{proof}  Use the injection $i_n: \kt {n, \varLambda} \rightarrow \akt 
{n, \varLambdaHat}$ of
Remark \ref{remark-injections of tangle algebras in affine tangle algebras}.
\end{proof}

\setcounter{ssection}{4}
\ssection{Isomorphism of}{ordinary BMW and Kauffman tangle 
algebras}\label{sec5}

This section is devoted to exhibiting an $R$-basis of $\kt {n, R}$ and
to proving that $\varphi: \bmw {n, R} \rightarrow \kt {n, R}$ is  an 
isomorphism.  
The results and arguments of this section are taken 
from~\cite{Morton-Wassermann}.

All the tangle diagrams in this section will be ordinary $(n,n)$-tangle diagrams 
for some~$n$.

\vfill\eject
\setcounter{section}{5}  
\subsection{Totally descending tangles and freeness 
  of $ \kt n$}\label{subsection-totally descending
tangles}
 
\begin{lemma}
$\kt {n, R}$ is spanned by tangle diagrams without closed 
strands.
 \end{lemma}
 
 \begin{proof}  The proof is by induction on the number of crossings.  

If a tangle diagram $T$  has no crossings, then, by Definition \ref{definition 
Kauffman tangle
algebra}(3),   $T = \delta^k T'$, where $k$ is the number of closed loops of 
$T$ and $T'$ is the
tangle diagram obtained by removing all the closed loops of~$T$.
 
 Let $T$ be a tangle diagram with $ l  \ge 1$ crossings.  Assume that any 
tangle diagram with fewer than
$ l $ crossings is in the span of tangle diagrams without closed strands.
 
 By the Kauffman tangle relation, if  $S$ is a  tangle diagram which differs 
from $T$ only by reversing
one or more crossings, then $T$ and $S$ are congruent modulo the span of tangle 
diagrams with fewer
crossings.  Hence $T$ and $S$ are congruent modulo the span of tangle diagrams 
without closed strands, by
the induction hypothesis.
 
 Suppose $T$ has a closed strand $s$.  By changing crossings, one can suppose 
that
 $s$ has only over-crossings with other strands of $T$ and that $s$ is 
unknotted.
 Then $T$ is ambient isotopic to a tangle diagram in which the closed strand 
corresponding to $s$ has no
crossings with other strands  and no self-crossings.   If $T'$ is the tangle 
diagram with $s$ removed,
then $T = \delta\lambda^k T'$ for some $k$, by
 Definition~\ref{definition Kauffman tangle algebra}(2,3).
 \end{proof}
 
\begin{definition}\rm
An {\em orientation} of an  affine or ordinary 
  $(n,n)$-tangle diagram  is
  a linear ordering of the strands,
 a choice of an orientation of each strand, and a choice of an initial point on 
each closed loop.
   \end{definition}
    
Order the boundary points  $\{\p 1, \dots, \p n, \pbar 1, \dots,
\pbar n\}$  of $(n,n)$-tangle diagrams by
$$
\p 1 < \p 2 < \cdots < \p n < \pbar n  < \cdots < \pbar 2 < \pbar 1,
$$
as in the discussion of colored $n$-connectors. 

\begin{definition}\rm
A {\em standard orientation} of an ordinary or affine 
$(n,n)$-tangle diagram is one in
which 
\begin{enumerate}
\item each non-closed strand is oriented from its lower numbered endpoint to its 
higher numbered endpoint,
\item the non-closed strands are ordered according to the order of their initial 
endpoints,
\item  the closed loops follow the non-closed strands in the ordering of the 
strands.
\end{enumerate}
\end{definition}

If a tangle diagram has no closed loops, then it has a unique standard 
orientation.

An orientation determines a way of traversing the tangle diagram;\break namely,  
the 
strands are traversed successively, in
the given order and orientation (the closed loops being traversed starting at 
the assigned initial point).

  \begin{definition}\rm
  An  oriented ordinary $(n,n)$-tangle diagram   is {\em 
totally descending} with respect to its orientation if, as the tangle diagram is 
traversed, each crossing is encountered first as an over-crossing.
  \end{definition}
  
\begin{proposition}\label{proposition-totally descending tangles
span}
$\kt {n, R}$ is spanned by ordinary $(n,n)$-tangle diagrams without closed 
loops that are  totally descending with respect to the standard orientation.
  \end{proposition}
  
  \begin{proof}  We already know that  $\kt {n, R}$ is spanned by $(n,n)$-tangle 
diagrams without closed loops.  
  
  If a tangle diagram has no crossings, it is already totally descending.
  Let $T$ be a tangle diagram  with $ l  \ge 1$ crossings. Assume that any 
tangle diagram  with fewer
than $ l $ crossings is in the span of totally descending tangle diagrams.   
The totally descending
tangle diagram $S$ which differs from $T$ only by reversing some number of 
crossings
  is congruent to $T$ modulo the span of tangle diagrams with fewer crossings, 
hence modulo the span of
totally descending tangle diagrams.
  \end{proof}

\begin{corollary}\label{corollary-totally descending tangles span}
$\kt {n, R}$ is spanned by the set
of  $(n,n)$-tangle diagrams $T$ without closed loops that are  totally 
descending with respect to the
standard orientation, and such that no strand of $T$ has self-crossings.
\end{corollary}

\begin{proof}  If $T$ is a totally descending  $(n,n)$-tangle diagram, then $T$ 
is layered; that is, $T$
can be drawn with different strands in different levels above the plane of $R = 
\R \times I$.   The
individual strands are unknotted, so they can be changed by level-preserving 
ambient isotopy to arcs
without self-crossings. Thus $T = \lambda^k T'$, where $T'$ is a totally 
descending tangle whose strands
have no self-crossings.
\end{proof}

\begin{proposition}\label{proposition-uniqueness of tangles with given 
connector}
Suppose $S$ and $T$ are two $(n,n)$-tangle diagrams without closed 
loops such that
\begin{enumerate}
\item $S$ and $T$  have the same connector, 
\item $S$ and $T$ are both totally descending (with respect to the same 
orientation),
\item the strands of $S$ and $T$ have no self-crossings.
\end{enumerate}
Then $S$ and $T$   are regularly isotopic, so they represent the same element of 
$\kt 
{n, R}$.
\end{proposition}

\begin{proof}  Since $S$ and $T$  have the same connector and are both totally 
descending with respect to the same ordering of the strands,  they can be 
layered, with the strands connecting corresponding endpoints in the two 
diagrams lying at the same level above the plane of $R$.    
Each strand of $S$ can then be deformed by a level-preserving isotopy to 
coincide with the corresponding strand of $T$; this deformation corresponds to 
regular isotopy of the diagrams.  
\end{proof}

For each $n$-connector $d$,  let $T_d$ be an $(n,n)$-tangle diagram without 
closed loops such that 
 \begin{enumerate}
 \item  $c(T_d) = d$,
 \item $T_d$ is totally descending with respect to the standard orientation,
 \item the strands of $T_d$  have no self-crossings.
 \end{enumerate}
By  Proposition~\ref{proposition-uniqueness of tangles with given connector}, 
 $T_d$ is unique up to regular isotopy.    The tangle diagrams $T_d$ (or rather 
the regular isotopy classes which they represent) can be regarded as elements of
$\kt {n, R}$ for any $R$.  For any~$R$,
$$
B_0 = \{ T_d : d \text{ is an $n$-connector} \}
$$
spans  $\kt {n, R}$, by  Corollary \ref{corollary-totally descending tangles 
span}.
By Corollary~\ref{corollary-independence of tangles with different connectors},  
$B_0$ is linearly independent in $\kt {n, \varLambda}$.   Thus we have:

\begin{theorem}\label{theorem- freeness of ktn}
 $\kt {n, \varLambda}$ is free over $\varLambda$ with basis  $B_0$.
\end{theorem}

\begin{corollary}\label{corollary- freeness of ktn and isomorphism of 
specializations}
 For each $R$,   $\kt {n, R}$ is free over $R$ with basis $B_0$.
   Moreover, $\kt {n, R} \cong \kt {n, \varLambda} \otimes_{\varLambda}  R$.
\end{corollary}

\begin{proof}
$\kt {n, \varLambda} \otimes_{\varLambda}  R$ is free over $R$ with basis 
$B_0 \otimes 1$.  On the other hand, $B_0$ spans
$\kt {n, R}$ by  Corollary~\ref{corollary-totally descending tangles span}.  
There is an $R$-algebra homomorphism
from $\kt {n, R}$ to $\kt {n, \varLambda} \otimes_{\varLambda}  R$ which sends a 
tangle $T$ to $T \otimes 1$.  Since this
map takes a spanning set to a basis, it is an isomorphism.
\end{proof}

\begin{corollary}\label{corollary-imbedding of ordinary BMW
algebras}
Let $R$ be a ring with distinguished elements
$\lambda$, $z$ and~$\delta$, and let
$S \supseteq R$ be a ring containing $R$.  Then $\kt {n, R}$ imbeds in  $\kt{n, 
S}$.
\end{corollary}

\begin{corollary}
 Let  $R$ be a ring with distinguished elements
$\lambda$, $z$, and $\delta$, and $S \supseteq R$  a ring with additional
elements $q_1, q_2, \dots$.  The $R$-algebra homomorphism 
$i_n: \kt {n, R} \rightarrow \akt {n, S}$  of Remark 
{\rm\ref{remark-injections of tangle algebras in affine tangle
algebras}} is injective.
\end{corollary}

\begin{proof}  This follows from the previous corollary and point (1) of Remark 
\ref{remark-injections of tangle
algebras in affine tangle algebras}.
\end{proof}

\subsection{Positive permutation braids}

\begin{definition}\rm
An $n$-{\em braid diagram} is an $(n,n)$-tangle diagram 
all of whose strands are monotone.  That is, each strand decreases monotonically 
from a top vertex to a bottom vertex. 
\end{definition}

\begin{definition}\rm
The (geometric) {\em braid group} $B_n$ is the group of 
$n$-braid diagrams modulo ambient isotopy. 
\end{definition}

The group $B_n$ has the well known presentation (due to Artin) with generators
 $\sigma_1, \ldots, \sigma_{n-1}$ and relations 
\begin{enumerate}
\item $\sigma_i \sigma_{i+1} \sigma_i = \sigma_{i+1} \sigma_i
\sigma_{i+1}$,
\item $\sigma_i \sigma_j = \sigma_j \sigma_i $  if $|i-j| \ge 2$. 
\end{enumerate}
The generator $\sigma_i$ is the $n$-braid diagram with a single positive 
crossing between the $i$th and
$(i+1)$st strand, 
$$
\sigma_i = \inlinegraphicee{sigma_i}.
$$
Since the generators $g_i$ of the 
ordinary BMW algebra satisfy the braid
relations,  $\psi : \sigma_i \mapsto g_i$ determines a group homomorphism from 
$B_n$ into the group of invertible elements of $\bmw n$.  Denote by perm the 
homomorphism from $B_n$ to the symmetric
group $\S_n$:  $\perm(\beta)(i) = j$ if  the braid diagram $\beta$ connects the 
top vertex $\p i$ with the bottom
vertex $\pbar j$.  In particular $\perm(\sigma_i)$ is the adjacent transposition 
$s_i = (i , i+1)$.  Note that
$\perm = c \circ \varphi \circ \psi$,  where we identify permutations with their 
diagrams in the Brauer algebra,
$$
\perm : \sigma_i \mapsto g_i \mapsto G_i \mapsto c(G_i) = s_i.
$$

\begin{proposition}\label{proposition-positive permutation braid}
The following are equivalent for an element $\beta$ of the braid group:
\begin{enumerate}
\item Two strands of $\beta$ cross at most once, and all crossings are positive 
(that is, $\beta$ is in the monoid generated by the  $\sigma_i$).
\item  $\beta$ is the product of $r$ generators $\sigma_i$, where $r$ is  the 
length of $\perm(\beta)$.
\end{enumerate}
Moreover, if $\beta$ satisfies these conditions, and $s_{i_r} \cdots s_{i_1}$ is 
any reduced expression for $\perm(\beta)$, then $\beta = \sigma_{i_r} \cdots 
\sigma_{i_1}$.
\end{proposition}

\begin{proof}   Suppose that $\beta =  \sigma_{i_r} \cdots \sigma_{i_1}$, but 
that two strands of $\beta$ cross twice.   Then
$\perm(\beta) = s_{i_r} \cdots s_{i_1}$ has a subword $s_a \pi_0 s_b$ with the 
property that 
$\pi_0(b) = a$ and $\pi_0(b+1) = a+1$.  But then $\pi_0 s_b = s_a \pi_0$, so 
$s_a \pi_0 s_b = \pi_0$, and the length of $\perm(\beta)$ is less than $r$.   
This proves (2)$\Rightarrow$(1).

Now suppose that $\beta$ satisfies (1).  Let $\perm(\beta) = s_{i_r} \cdots 
s_{i_1}$ be a reduced expression for $\perm(\beta)$ and set $\beta' =  
\sigma_{i_r} \cdots \sigma_{i_1}$.  Then $\beta$ and
$\beta'$ are two braid diagrams both satisfying condition 
(1) with $\perm(\beta) = \perm(\beta')$.
But a braid diagram satisfying condition (1) is totally descending 
(with respect to the orientation
in which the strands are oriented from top to bottom and ordered in the reversed 
order of their
top vertices).  Therefore, by Proposition~\ref{proposition-uniqueness of tangles 
with given connector},  such a braid diagram $\beta$ is determined up to ambient 
isotopy by its connector, that is, by $\perm(\beta)$.
\end{proof}

\begin{definition}\rm\label{definition-positive permutation braid}
A braid diagram satisfying the conditions of the previous proposition is called 
a {\em positive permutation braid}.
\end{definition}

For each $\pi \in \S_n$, there is a unique positive permutation braid $\beta_\pi 
\in B_n$ with $\perm(\beta_\pi) = \pi$.     We  write $g_\pi$ for the image of 
$\beta_\pi$ in $\bmw n$ and
$G_\pi$ for the image of $\beta_\pi$ in $\kt n$.   We also call these elements 
(which are determined by~$\pi$) positive permutation braids.  
Note that $g_{\pi\inv} = \alpha(g_\pi)$.

For $\pi \in \S_n$  and $1 \le i < n$, we have $ \ell (\pi s_i) =  \ell (\pi) + 1 
\Leftrightarrow \pi(i) < \pi(i+1)$.  In this case, $g_{\pi s_i} = g_\pi 
g_i$.  Otherwise, $ \ell (\pi s_i) =  \ell (\pi) - 1$, $\pi(i) > \pi(i+1)$, and
$g_{\pi s_i} = g_\pi g_i\inv$.  Likewise,  $ \ell (s_i\pi)  =  \ell (\pi) + 1 
\Leftrightarrow \pi\inv(i) < \pi\inv(i+1)$. In this case,
 $g_{s_i \pi} = g_i g_\pi $.  Otherwise, $ \ell (\pi s_i) =  \ell (\pi) - 1$, 
$\pi\inv(i) > \pi\inv(i+1)$, and
 $g_{s_i \pi} = g_i\inv g_\pi$. 
 
\begin{definition}\rm
An $(a, b)$-{\em shuffle} is an element $ \pi \in 
\S_{a+b}$ such that
 $\pi(i) < \pi(j)$ if  $1 \le i < j \le a$ or $a+1 \le i < j \le a+b$.
 \end{definition}
 
\begin{lemma}\label{lemma- factorization of positive permutation
braid}
 If $\pi \in \S_{a+b}$ then $\pi = \pi_1 \pi_2$ where $\pi_1$ is an 
$(a, b)$-shuffle,
 $\pi_2 \in \S_a \times \S_b \subseteq \S_{a+b}$, and $\ell (\pi) = \ell (\pi_1) + 
\ell (\pi_2)$.  It follows that
 $g_\pi = g_{\pi_1} g_{\pi_2}$.
 \end{lemma}
 
 \begin{proof}  The proof is by induction on the length of $\pi$.  The result is 
evident if $\pi$ is the identity permutation.  If $\pi$ is not already an $(a, 
b)$-shuffle, then there exists an $i$ with $i \ne a$ such that
 $\pi(i) > \pi(i+1)$.  Consequently, $\pi = \pi' s_i$ with $\ell (\pi) = 
\ell (\pi') + 1$.  The result follows by applying the induction hypothesis
to~$\pi'$.
 \end{proof}
 
\begin{lemma}\label{lemma- transport gi by parallel strands}
Let $1 \le i \le n-1$ and suppose
$\pi \in \S_n$ satisfies  $\pi\inv(i+1) = \pi\inv(i) + 1$.
Then $g_i g_\pi  = g_\pi g_{\pi\inv(i)}$ and
$e_i g_\pi  = g_\pi e_{\pi\inv(i)}$.
\end{lemma}

{\em Proof}.\hods  Put $j = \pi\inv(i)$.
We have $s_i \pi = \pi s_{j}$, and $\ell (s_i \pi ) = \ell (\pi) + 1$.
Hence 
$$
g_\pi g_{j}=  g_{ \pi s_j}= g_{s_i \pi} = g_i g_\pi.
$$ 

The second equality is proved by induction on the length of $\pi$.  If 
$\ell (\pi) = 0$, the assertion is trivial.  Suppose that $\ell (\pi) > 1$.  
Choose $k$ such that
$\pi = s_k \pi_1$, and $\ell (\pi) = \ell (\pi_1) + 1$.
To prove the induction step, we consider three cases:

\begin{enumerate}
\item $k = i+1$.  In this case $\pi\inv(i +2) < \pi\inv(i+1) = \pi\inv(i) +1$.  
It follows that $\pi\inv(i+2) < \pi\inv(i)$, and $\pi = s_{i+1} s_i \pi'$, with 
$\ell (\pi) = \ell (\pi') + 2$.  Therefore
$g_\pi = g_{i+1} g_{i} g_{\pi'}$ and
$e_i g_\pi = e_i g_{i+1} g_i g_\pi' = g_{i+1} g_{i} e_{i+1} g_{\pi'}$.  Since
${\pi'}\inv(i+2) = j+1 ={ \pi'}\inv(i+1) + 1$,  we have $ e_{i+1} g_{\pi'} = 
g_{\pi'} e_j$, by the induction hypothesis.
\item $k = i-1$.  This case is similar.
\item $|k-i| \ge 2$.  Then $e_i g_\pi = e_i g_k g_{\pi_1} = g_k e_i g_{\pi_1}$. 
Now 
${\pi_1}\inv(i+1) = j+1 = {\pi_1}\inv(i) + 1$, so by the induction hypothesis,
$e_i g_{\pi_1} = g_{\pi_1} e_j$.\squ
\end{enumerate}

\subsection{Surjectivity of  $\varphi: W_n \rightarrow \kt n$.}
We will prove that  $\kt {n, R}$ is generated as a unital algebra by $\{E_i, 
G_i\pmone :  1 \le i \le n-1\}$, which is equivalent to the surjectivity of     
$\varphi: \bmw {n, R} \rightarrow \kt {n, R}$.

The tensor product $T_1 \otimes T_2$ of a $(k, k)$-tangle diagram and an 
$(l ,l )$-tangle diagram is the $(k + l , k + l )$-tangle diagram 
obtained by placing $T_1$ and $T_2$ side by side.   

The tensor product of tangle diagrams clearly extends to 
 a bilinear product $\kt {k, R} \times \kt {l , R} \rightarrow \kt {k + 
l , R}$.  If $T_1$ and 
 $T_2$ are both in the unital subalgebra generated by the $E_i$'s and $G_i$'s, 
then
 so is $T_1 \otimes T_2$. 

\begin{lemma}
Any element of $\kt {n, R}$ which is represented by an  $n$-braid 
diagram is in the group generated by $\{G_i\pmone: 1 \le i \le n-1\}$.
\end{lemma}
 
 \begin{proof}  Induction on the number of crossings.
 \end{proof}
 
\begin{theorem}\label{theorem- surjectivity of varphi on Wn}
$\varphi: \bmw 
{n, R} \rightarrow \kt {n, R}$ is
surjective.
 \end{theorem}
 
 \begin{proof}
 For $n = 0$ and  $n = 1$, $\kt {n, R} \cong R$ by  
Corollary~\ref{corollary-  freeness of ktn and isomorphism of specializations}, 
and the 
statement is trivially valid.  

 Fix $n \ge 2$.  We have to show that
 $\kt {n, R}$ is generated as a unital algebra by 
 $\{E_i, G_i\pmone :  1 \le i \le n-1\}$.     
  By Theorem \ref{theorem- freeness of ktn}, it suffices to show that
 each totally descending $(n,n)$-tangle diagram  $T$ whose  strands have no 
self-crossings is in the unital  subalgebra generated by the $E_i$'s and 
$G_i$'s. 

 If the connector $c(T)$ is a permutation diagram  (i.e. top vertices are 
connected only to bottom vertices), then $T$ is regularly isotopic to an 
$n$-braid diagram, and thus
 $T$ is in the monoid generated by  $\{G_i\pmone: 1 \le i \le n-1\}$.

Otherwise,  for some $k \ge 1$,  the connector of $T$ has $k$ horizontal strands
connecting pairs of vertices at the top and $k$ horizontal strands connecting 
vertices at the bottom.  In this case there exist $n$-braid diagrams $B_1$ and 
$B_2$ and
an $(n-2k)$-braid diagram $B$ such that
$$
T = B_1[ (E_1 E_3 \cdots E_{2k-1}) \otimes B] B_2.
$$
Since each of $B_1, B_2$, and $B$ are in the group generated by the 
$\{G_i\pmone\}$,
it follows that $T$ is in the monoid generated by $\{E_i, G_i\pmone\}$.
\end{proof}

\subsection{The elements $f_k$ and a filtration of $\bmw n$.}
Fix a ring $R$ and write $\bmw n$ for $\bmw {n, R}$ and 
 $\kt n$ for $\kt {n, R}$.

Consider the element $F_k \in \kt {2k}$  represented by the tangle diagram with 
no crossings, in which the points
$\p i $ and $\p {2k} + \p{1} - \p{i} $  at the top of the diagram are connected, 
and 
likewise
the points $\pbar i $ and $\overline{\mathbold {2k} + \mathbold{1} -
\mathbold{i}} $  at the bottom of 
the diagram are connected.  For example
$$
F_3 =  \inlinegraphic{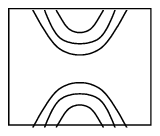}
$$
The element $F_k$ is evidently fixed by the symmetries $\alpha$, $\beta$, and    
$\varrho_{2k}$ of $\kt {2k}$.

The following proposition is from~\cite{Morton-Wassermann} and is easily 
verified by 
picture proofs.

\begin{proposition}\label{proposition-properties of Fk}
\mbox{}
\begin{enumerate}
\item For all $i < k$,  $G_i\pmone F_k = G_{2k - i}\pmone F_k$ and $F_k 
G_i\pmone= F_k G_{2k -i}\pmone$.
\item  For all $i < k$, $E_i F_k = E_{2k - i} F_k$ and $F_k E_i = F_k E_{2k -
i}$.
\item $F_k = (G_1 G_2 \cdots G_{2k-1}) F_{k-1} E_{2k -1} (G_{2k-1} \cdots G_2 
G_1)$.
\end{enumerate}
\end{proposition}

Following~\cite{Morton-Wassermann}, 
we recursively define  elements $f_k$ of $\bmw {2k}$ such that $\varphi(f_k)\break 
= 
F_k$, as follows:
 
\begin{definition}\rm
Define $f_1 = e_1$ and 
$$
f_k = (g_1 g_2 \cdots g_{2k-2}) f_{k-1} e_{2k-1} 
(g_{2k-2} \cdots g_2 g_1) 
$$
for $k \ge 2$.
\end{definition}

We want to find an expression for  $f_k$ that does not involve the $g_i$'s  
(since the corresponding tangle diagrams have no crossings), and that makes 
manifest the symmetries
$\alpha(f_k) = \beta(f_k) = \varrho_{2k}(f_k) = f_k$ (which the definition does 
not).

\begin{lemma}  
\begin{eqnarray*}
(g_1 g_2 \cdots g_{2k}) 
 (e_1 e_3 \cdots e_ {2k -1})  &=& 
 (e_2 e_4 \cdots e_{2k}) (e_1 e_3 \cdots e_ {2k -1}) , \\
 (e_1 e_3 \cdots e_ {2k -1}) (g_{2k} \cdots  g_2 g_1)  &=&
 (e_1 e_3 \cdots e_ {2k -1})  (e_2 e_4 \cdots e_{2k}) .
\end{eqnarray*}
\end{lemma}

\begin{proof}  The second equation follows from the first by applying the 
anti-automorphism $\alpha$.   To prove the first equation, rewrite the left hand 
side 
as
$$
(g_1 g_2 e_1) (g_3 g_4 e_3) \cdots (g_{2k-1} g_{2k} e_{2k-1}),
$$
which equals
$$
(e_2 e_1) (e_4 e_3)  \cdots (e_{2k} e_{2k-1}),
$$
by use of (5b) of Definition \ref{definition affine BMW}.
Finally, the even terms can be shuffled to the left.
\end{proof}

\begin{proposition}\label{proposition-diamond expression for fk}
For $k \ge 1$, 
$$
f_{k} = (e_{k})(e_{k-1}e_{k+1})\cdots (e_{1}e_{3}\cdots e_{2k-1})
 \cdots (e_{k-1}e_{k+1}) (e_{k}).
$$
\end{proposition}

For example,
$$
f_4 = e_4 (e_3 e_5) (e_2 e_4 e_6) (e_1 e_3 e_5 e_7) (e_2 e_4 e_6) (e_3 e_5) e_4.
$$

\begin{proof}   The proof goes by induction on $k$,  the  assertion being 
evident for $k=1$.  
Consider $k >2$.  We have
$$\displaylines{\quad\ 
\eqalign{
f_k &=  (g_1 g_2 \cdots g_{2k-2}) f_{k-1} e_{2k-1} (g_{2k-2} \cdots g_2g_1)\cr
&= (g_1 g_2 \cdots g_{2k-2})[ e_{k-1} (e_{k-2} e_{k}) \cdots    
(e_1 e_3 \cdots e_ {2k-3}) 
\cdots (e_{k-2} e_{k})  e_{k-1} ) ] \cr}\hfill\cr
\hfill \cdot\, e _{2k-1}  (g_{2k-2} \cdots g_2 g_1)\quad\ \cr} 
$$
by the induction hypothesis.    Move $e_{2k-1}$ to the left to get
$$\displaylines{
= (g_1 g_2 \cdots g_{2k-2})[ e_{k-1} (e_{k-2} e_{k}) \cdots   
(e_1 e_3 \cdots e_ {2k-3}e _{2k-1}) \cdots 
(e_{k-2} e_{k})  e_{k-1} ) ] \hfill\cr 
\hfill \cdot\,  (g_{2k-2} \cdots g_2 g_1). \cr} 
$$
Now move  all of the  $e_i$'s, except those in the product
 $(e_1 e_3 \cdots e_{2k-1})$, to the left or right  past  a string of $g_i$'s, 
making use of Lemma \ref{lemma-commutation relations 1}; this yields
$$\displaylines{\quad\ 
= e_k (e_{k-1} e_{k+1}) \cdots 
(e_3 e_5 \cdots e_{2k - 3}) (g_1 g_2 \cdots g_{2k-2}) 
(e_1 e_3 \cdots e_ {2k -3}e _{2k-1})\hfill\cr 
\hfill \cdot\, (g_{2k-2} \cdots g_2 g_1)(e_3 e_5 \cdots e_{2k - 3})  \cdots(e_{k-
1} e_{k+1})  
e_k. \quad\ \cr} 
$$
Both strings of $g_i$'s in the middle of the expression can be replaced by
$e_2 e_4 \cdots e_{2k -2}$, by the previous lemma, and this gives the desired 
expression.
\end{proof}

\begin{remark}\rm
The elements in each group (indicated by parentheses) commute.  
It can be helpful to view the  entire expression as a diamond-shaped grid, for 
example
$$
f_{3} = 
\begin{aligned}  
&{}&{}&e_{1}&{}&{}  \\
&{}&e_{2}\  &{}&e_{2}\ &{} \\
&e_{3}&{}&e_{3}&{} &e_{3}\\
&{}&e_{4}\  &{}&e_{4}\ &{}\\
&{}&{}&e_{5}&{}&{} \\
\end{aligned}
$$
This is read from left to right by columns, with the elements in each column 
commuting.  
This expression for $f_k$ makes  evident the invariance of $f_k$ under the  maps 
$\varrho_{2k}$,  
$\alpha$, and $\beta$ of $\bmw {2k}$.

Analogous  elements in the Temperley--Lieb algebras were introduced in subfactor
theory  by M. Pimsner and S. Popa~\cite{Pimsner-Popa-iterating} to study 
iterations of the Jones basic construction;  see  also~\cite{Jones-Sunder}.  The 
diamond grid representation of these elements is due to A. Ocneanu (personal 
communication).

Note that one can also read the diamond grid by diagonals, so we have
$$
f_k = (e_k e_{k+1} \cdots e_{2k-1}) (e_{k-1} e_k \cdots e_{2k-2}) \cdots (e_1 
e_2 \cdots e_k).
$$
Moreover,   $f_k$ satisfies
$$
f_k = (e_k e_{k+1} \cdots e_{2k-2}) f_{k-1} e_{2k-1}  (e_{2k-2} \cdots e_{k+1} 
e_k).
$$
\end{remark}

The following proposition from~\cite{Morton-Wassermann} establishes the analogue
of  Proposition \ref{proposition-properties of Fk} for the elements $f_k$.  In 
points (3) and (4), $S$ denotes the shift homomorphism (see 
Remark \ref{remark- injection-  shift- and reversal}).

\begin{proposition}\label{proposition-transfer property of fk 1}
\mbox{}
\begin{enumerate}
\item For $i < k$,  $g_i\pmone f_k = g_{2k-i}\pmone f_k$  and $f_k g_i\pmone= 
f_k g_{2k-i}\pmone$.
\item For $i  < k$,  $e_i f_k = e_{2k-i} f_k$  and $f_k e_i = f_k e_{2k-i}$.
\item For $2 \le i < k+1$, $g_i\pmone S(f_k) = g_{2k+2 -i}\pmone S(f_k)$  and 
$S(f_k) g_i\pmone = S(f_k)g_{2k+2-i}\pmone$.
\item For $2 \!\le\! i  \!<\! k+1$, $e_i S(f_k) \!=\! e_{2k+2 -i} S(f_k)$  and 
$S(f_k) e_i \!=\! 
S(f_k)e_{2k+2-i}$.
\end{enumerate}
\end{proposition}

{\em Proof}.\hods The second part of each assertion follows from the first part by 
applying the anti-auto\-mor\-phism~$\alpha$.

 Note that $i$ and $2k-i$ have the same parity, so $g_{i}$ and $g_{2k-i}$ 
commute.  Therefore,
 multiplying the equality
 \begin{equation}\label{equation: fk transfer property1}
 g_i f_k = g_{2k-i} f_k
\end{equation}
 by $g_{i}\inv g_{2k-i}\inv$ gives
\begin{equation}\label{equation: fk transfer property2}
 g_i\inv f_k = g_{2k-i}\inv f_k.
 \end{equation}
 Thus, for statement (1), it suffices to prove equation (\ref{equation: fk 
transfer property1}).
 
 We prove statement (1) by induction on $k$ (following the proof 
in~\cite{Morton-Wassermann}).      
For $k = 2$, the only instance to check is $i = 1$. Using 
Lemma \ref{lemma-commutation relations 1}(3),  we get
\[
g_1 f_2 = g_1 (g_1  g_2 e_1 e_3) g_2 g_1 = g_1 (g_3 g_2 e_3 e_1) g_2 g_1
= g_3 (g_1 g_2 e_1 e_3 g_2 g_1) = g_3 f_k.
\]
Consider $k > 2$.  If $ 2 \le i  < k$, then, using Lemma  \ref{lemma-commutation 
relations 1}, we have
\begin{align*}
g_i f_k &= g_i (g_1 g_2  \cdots  g_{2k-2}) f_{k-1} e_{2k-1} (g_{2k-2} \cdots g_2 
g_1) \\
&= (g_1 g_2  \cdots  g_{2k-2}) g_{i-1} f_{k-1} e_{2k-1} (g_{2k-2} \cdots g_2 
g_1) \\
&= (g_1 g_2  \cdots  g_{2k-2}) g_{2k -i -1}  f_{k-1} e_{2k-1} (g_{2k-2} \cdots 
g_2 g_1),
\end{align*}
which by the induction hypothesis equals
\[
 g_{2k -i } (g_1 g_2  \cdots  g_{2k-2})  f_{k-1} e_{2k-1} (g_{2k-2} \cdots g_2 
g_1) =  g_{2k -i } f_k.
\]
The last case to check is $k > 2$ and $i = 1$. We have
\begin{align*}
g_1 f_k &= g_1 (g_1 g_2  \cdots  g_{2k-2} )  f_{k-1} e_{2k-1} (g_{2k-2} \cdots 
g_2 g_1) \\
&= g_1 (g_1 g_2  \cdots  g_{2k-2} )(
(g_1 g_2  \cdots  g_{2k-4} )  f_{k-2} e_{2k-3} (g_{2k-4} \cdots g_2 g_1)
) \\
&\phantom{========================} \cdot\,e_{2k-1} (g_{2k-2} \cdots g_2 g_1) \\
&= g_1  (g_2 g_3  \cdots  g_{2k-3} )(g_1 g_2  \cdots  g_{2k-2} )
  f_{k-2} e_{2k-3} (g_{2k-4} \cdots g_2 g_1)
 \\
 &\phantom{========================}\cdot\,e_{2k-1} (g_{2k-2} \cdots g_2 g_1),
 \end{align*}
by repeated  use of Lemma  \ref{lemma-commutation relations 1}(4).  Moving
$e_{2k-3}$ and  $e_{2k-1}$ to the left and applying  
Lemma~\ref{lemma-commutation relations 1}(3), 
we get
$$\displaylines{\quad\
g_1  (g_2 g_3  \cdots  g_{2k-3} )(g_1 g_2  \cdots g_{2k-4}  (g_{2k-3} g_{2k-2} 
e_{2k-3}  e_{2k-1}) )\hfill\cr
\hfill \cdot\,f_{k-2}  (g_{2k-4} \cdots g_2 g_1)
  (g_{2k-2} \cdots g_2 g_1) \quad\ \cr
\quad\ = g_1  (g_2 g_3  \cdots  g_{2k-3} )(g_1 g_2  \cdots g_{2k-4} (g_{2k-1} 
g_{2k-2} e_{2k-3}  e_{2k-1}) )\hfill\cr
\hfill
 \cdot\, f_{k-2}  (g_{2k-4} \cdots g_2 g_1)  (g_{2k-2} \cdots g_2 g_1).\quad\
\cr}
$$
Now $g_{2k-1}$ and $g_{2k-2}$  can be moved to the left  and 
  $e_{2k-1}$ and $e_{2k-3}$ to the right to yield
$$\displaylines{\quad\
g_{2k-1} g_1  (g_2 g_3  \cdots  g_{2k-3}g_{2k-2} )(g_1 g_2  \cdots g_{2k-4})
  f_{k-2} (  e_{2k-3}  e_{2k-1}) \hfill\cr
\hfill
\cdot\, (g_{2k-4} \cdots g_2 g_1) (g_{2k-2} \cdots g_2 g_1)\quad \ \cr
\hfill =\,  g_{2k-1} f_k.\quad\ \cr}
$$
This completes the proof of (1).

Statement (2) is also proved by induction on $k$.  For $k = 2$ and $i = 1$, we
have
$$
e_1 f_2 = e_1 e_2 e_1 e_3 e_2 = e_1 e_3 e_2,
$$
and similarly $e_3 f_2 = e_1 e_3 e_2.$

For $2 \le i < k$, exactly the same induction step can be used as in the proof 
of state\-ment~(1).  The only remaining case to check is $k >2$ and $i = 1$.

We claim that $e_1 f_k$ is the word $w_k$ in the $e_i$'s obtained by deleting 
the two columns to the left of the middle column in the diamond expression 
for~$f_k$.  For example,
$$
e_1 f_{4} =  e_1
\left [
\begin{aligned}  
&{}\ &{}\ &{}\ &e_{1}\ &{}\ &{} \ &{}\\
&{}\ &{}\ &e_{2}\  &{}\ &e_{2}\ &{} \ &{}\\
&{}\ &e_{3}\ &{}\ &e_{3}\ &{} \ &e_{3}\ &{}\\
&{e_4}\ &{}\ &e_{4}\  &{}\ &e_{4}\ \ &{}\ &{e_4}\\
&{}\ &{e_5}\ &{}\ &e_{5}\ &{}\ &{e_5}\ &{} \\
&{}\ &{}\ &e_{6}\  &{}\ &e_{6}\ &{} \ &{}\\
&{}\ &{}\ &{}\ &e_{7}\ &{}\ &{}\\
\end{aligned}  \right ] \quad
=  \quad e_4 \left [
\begin{aligned}  
&e_{1}\ &{}\ &{} \ &{}\\
&{}\ &e_{2}\ &{} \ &{}\\
&e_{3}\ &{} \ &e_{3}\ &{}\\
&{}\ &e_{4}\ \ &{}\ &{e_4}\\
&e_{5}\ &{}\ &{e_5}\ &{} \\
&{}\ &e_{6}\ &{} \ &{}\\
&e_{7}\ &{}\ &{}\\
\end{aligned} 
\right ].
$$
This claim can be proved by induction on $k$.   Moreover, both $f_k$ and $w_k$ 
are 
invariant under the automorphism $\varrho_{2k}$ of $\bmw {2k}$.  Applying this 
automorphism
to the equation $e_1 f_k = w_k$ gives
$$
e_{2k-1} f_k = w_k = e_1 f_k.
$$ 
This completes the proof of (2).

Statements (3) and (4) follow from (1) and (2) by applying the shift 
homomorphism $S$, for example
$$
e_{i}S(f_{k}) = S(e_{i-1} f_{k}) = S(e_{2k + 1 -i} f_{k}) = e_{2k+2 -i}
S(f_{k}).\squ
$$

\begin{remark}\rm
One can give a similar explicit expression for $e_i f_k$  for 
all $i < k$.
One obtains $e_{k-1} f_k$  by deleting the leftmost $e_k$ from the diamond 
expression for  $f_k$.  
Moreover,  if $i < k-1$, then  $e_i f_k $ is the word in the $e_j$'s obtained by 
deleting the  two columns
 just to the left of the first column beginning with
$e_i$ in the diamond expression for $f_k$.  For example,
$$
e_2 f_4 =  \left [
\begin{aligned}  
 &{}\ &e_{1}\ &{}\ &{} \ &{}\\
&e_{2}\  &{}\ &e_{2}\ &{} \ &{}\\
&{}\ &e_{3}\ &{} \ &e_{3}\ &{}\\
 &e_{4}\  &{}\ &e_{4}\ \ &{}\ &{e_4}\\
 &{}\ &e_{5}\ &{}\ &{e_5}\ &{} \\
&e_{6}\  &{}\ &e_{6}\ &{} \ &{}\\
&{}\ &e_{7}\ &{}\ &{}\\
\end{aligned}   \right ].
$$
The  proof of statement (2) in the previous proposition can be based entirely on 
this observation.
\end{remark}

\begin{proposition}\label{proposition- fk left ideal}
\mbox{}
\begin{enumerate}
\item For all $k$,  $\bmw {2k} f_k = \bmw {k} f_k$ and $f_k \bmw {2k} = f_k 
\bmw {k}$.
\item For all $k$,  $\bmw {2k + 1} S(f_k) = \bmw {k + 1} S(f_k)$ and $S(f_k)  
\bmw {2k + 1} =   S(f_k) \bmw {k + 1}$.
\end{enumerate}
\end{proposition}

\begin{proof}  The second part of each statement follows from the first  part by 
applying the symmetry $\alpha$.

To prove (1), it suffices to show that $ \bmw {k} f_k$ is a left ideal in $\bmw 
{2k}$.  For this, it is enough to show that
$e_i  \bmw {k} f_k \subseteq \bmw {k} f_k$ and $g_i\pmone \bmw {k} f_k \subseteq 
\bmw {k} f_k $ for $k \le i \le 2k-1$.
For $k < i < 2k-1$, this follows from Proposition \ref{proposition-transfer 
property of fk 1}, since $e_i$ and $g_i\pmone$ commute with $\bmw {k}$.  
Finally, by Proposition \ref{proposition-diamond expression for fk}, we can 
write
$f_k = e_k r_k$ for some $r_k \in \bmw {2k}$.  Hence for $\chi \in \{e_k, 
g_k\pmone\}$, 
$$
\chi  \bmw {k} f_k  =  \chi  \bmw {k} e_k r_k \subseteq \bmw {k + 1} e_k r_k  =  
\bmw {k} e_k r_k = \bmw {k} f_k,
$$
using Corollary \ref{corollary- bimodule 5.5}.    The proof of part (2) is 
similar.
\end{proof}

\begin{definition}\rm\label{definition- wnr}
Let $n \ge 2k$ and put $r = n-2k$.  Define 
$\w n r$ to be the ideal in $\bmw {n}$ generated by $e_1 e_3 \cdots e_{2k-1}$. 
\end{definition}

 Evidently, one has
$$\eqalign{
\w n 0 \subseteq \w n 2\subseteq \cdots \subseteq \w n n = \bmw n\quad\ &
\mbox{if $n$ is even,}\cr
\w n 1 \subseteq \w n 3 \subseteq \cdots \subseteq \w n n = \bmw n\quad\
&\mbox{if $n$ is odd.}\cr}
$$

Morton and Wassermann  show that 
 $\varphi(\w n r)$ is the ideal  $\k n r$ in $\kt {n}$ spanned by tangle 
diagrams  of rank no more than $r$, 
 that is, tangle diagrams of the form
 $S T$, where $T$ is an $(n,r)$-tangle diagram and $S$ is an $(r, n)$-tangle 
diagram.   We will not need this observation, but it provides the motivation for 
the definition of $\w n r$.

\begin{lemma}\label{lemma- w n r ideal}
 Let $n \ge 2k$ and put $r = n-2k$.
\begin{enumerate}
\item  $\w n r$ is the ideal in $\bmw {n}$ generated by $f_k$.
\item  If $n > 2k$ then $\w n r$ is the ideal in $\bmw {n}$ generated by 
$S(f_k)$. 
\end{enumerate}
\end{lemma}
 
\begin{proof}  It follows from the definition of $f_k$, and induction on $k$, 
that $f_k = a (e_1 e_3 \cdots e_{2k-1}) b$ where
$a$ and $b$ are in the monoid generated by $\{g_j\}$, and in particular are 
invertible.  Likewise, if $n > 2k$, then
$S(f_k) =  {\rm Ad}(a)(f_k)$ where $a$ is invertible.
\end{proof}

\begin{corollary}\label{corollary- expression for w2k 0}
\mbox{}
\begin{enumerate}
\item  $\w {2k} 0 =  \bmw {k} f_k \bmw {k}$.
\item  $\w {2k + 1} 1 =  \bmw {k + 1} S(f_k) \bmw {k + 1}$.
\end{enumerate}
\end{corollary}

\begin{proof}  Immediate from Proposition \ref{proposition- fk left ideal} and 
Lemma \ref{lemma- w n r ideal}.
\end{proof}

\subsection{Injectivity of $\varphi: \bmw n \rightarrow \kt n$.}\label{subsection-
injectivity of varphi}
Fix a ring $R$ and write $\bmw n$ for $\bmw {n, R}$ and 
$\kt n$ for $\kt {n, R}$.
This section contains Morton and Wassermann's  proof of the  injectivity of 
$\varphi: \bmw n \rightarrow \kt n$.  
The strategy is to show that $\varphi$ is injective on $\w n r$ for all $n$ and 
$r$, by double induction.

\begin{proposition}\label{proposition- inductive step from w n-1 to w n
0}
 Suppose $\varphi |_{  \bmw {n-1}}$ is injective for some $n$.  Then
$\varphi |_{\w n 0}$ is injective  if $n$ is even, and 
$\varphi |_{\w n 1}$ is injective  if $n$ is odd.
\end{proposition}

\begin{proof}   Consider the case that $n = 2k$.  By Corollary \ref{corollary- 
expression for w2k 0},  $\w n 0 =\break  \bmw {k} f_k \bmw {k}$.  By Theorem 
\ref{theorem- freeness of ktn}, $\kt k$ has a basis $\{T_c:   c \ 
\text{a $k$-connector}\}$ consisting of totally descending tangle diagrams with 
distinct 
connectors.  
Since by hypothesis,
$\varphi : \bmw {k}  \rightarrow \kt k$ is an isomorphism,  it follows that
$\{\varphi\inv(T_c):   c \ \text{a $k$-connector}\}$ is a basis of $\bmw {k}$.
Therefore
$\w n 0$ is spanned by the set  
$$
\{\varphi\inv(T_c) f_k \varphi\inv(T_d) : c, d \ \text{$k$-connectors}\}.
$$
The image of this spanning set under $\varphi$ is $\{ T_c F_k T_d :  c, d \ 
\text{$k$-connectors}\}$.  This set of $(n,n)$-tangle diagrams has distinct  
connectors, and therefore is linearly independent by Corollary 
\ref{corollary-independence of tangles with different connectors}.  It follows 
that
$\varphi$ is injective on $\w n 0$.  The proof for the case $n = 2k + 1$ is 
similar.
\end{proof}

Next we introduce a linear complement of $\w n {r-2}$ in $\w n r$.
 
\begin{definition}\rm\label{definition- vnr}
Let $n > 2k$ and put $r = n - 
2k$. 
Then $ \V n r$ is the span of $\{ g_\pi [ w \otimes g_\tau]  \alpha(g_\sigma)\}$, 
where $w \in \w {2k} 0$, 
$\pi$ and $\sigma$ are $(2k, r)$-shuffles, $\tau$ is an $r$-permutation,  and 
$g_\pi$, $g_\sigma$, and $g_\tau$ are the corresponding positive permutation 
braids.
\end{definition}

\begin{lemma}
Assume that $\varphi |_{\bmw {n-1}}$ is injective for a particular 
$n$. 
\begin{enumerate}
\item $ \varphi  |_{\V n r}$ is injective for $r > 0$.
\item $\V n r \cap \w n {r-2} = (0)$.
\end{enumerate}
\end{lemma}

\begin{proof}  By assumption, $\varphi |_{\w {2k} 0}$ is injective, and 
moreover, $\w {2k} 0$ has a basis
consisting of elements $\varphi\inv(T_c) f_k \varphi\inv(T_d)$, where $T_c$ and 
$T_d$ vary independently over a set of totally descending $(k, k)$-tangle 
diagrams with distinct connectors.  Thus
$\V n r$ has a spanning set  
$$
\{ g_\pi [\varphi\inv(T_c) f_k \varphi\inv(T_d)  
\otimes g_\tau]  \alpha(g_\sigma)\}.
$$
The image of this set,  
$$
\{ G_\pi [T_c F_k T_d  \otimes G_\tau]  \alpha(G_\sigma)\},
$$
is a family of tangle diagrams with distinct connectors
$
\pi[ c F_k d \otimes \tau] \sigma\inv,
$
and is therefore linearly independent by 
Corollary \ref{corollary-independence of tangles with different connectors}.  
This proves statement~(1).

For (2), $\varphi(\w n {r-2})$ is spanned by totally descending tangle 
diagrams with 
no more than $r -2$ through strands (i.e. strands connecting top to bottom), 
while 
$\varphi(\V n r)$ is spanned by totally descending tangle diagrams with exactly 
$r$ through strands.
By Corollary \ref{corollary-independence of tangles with different connectors}, 
$\varphi(\V n r) \cap \varphi (\w n {r-2}) = (0)$, so also $\V n r \cap \w n {r-
2} = (0)$.
\end{proof}

\begin{corollary}\label{corollary- inductive step w n r-2 to v nr plus w n
r-2}
 If $\varphi$ is injective on $\bmw {n-1}$ and also on $\w n {r-2}$, then
$\varphi$ is injective on $\w n {r-2} \oplus \V n r$.
\end{corollary}

The remainder of the proof consists of showing that  $\w n {r-2} \oplus \V n r = 
\w n r$.  Since 
$e_1 e_3 \cdots e_{2k-1} \in \V n r$, it suffices to show that $\w n {r-2} 
\oplus \V n r$ is an ideal
in $\bmw n$.  Because of the invariance of $\w n {r-2} \oplus \V n r$ under the 
anti-automorphism $\alpha$, it is enough to show that $\w n {r-2} \oplus \V n r$ 
is a left ideal.  It suffices to show that $\chi \V n r \subseteq \w n {r-2} 
\oplus \V n r$ for $\chi$ a generator of $\bmw n$, because
 $\w n {r-2}$ is already an ideal.    Moreover, since
$g_i\inv = g_i - z e_i + z$, it suffices to show this for  $\chi \in \{e_i, g_i 
: 1 \le i \le n-1\}$.

Note that $\V n n$ is the span of $\{g_\pi : \pi \in \S_n\}$ and $\w n {n-2} = 
\bmw n e_1 \bmw n$, so   a particular case of our assertion is that $\V n n 
\oplus \bmw n e_1 \bmw n = \bmw n$.  

\begin{lemma}\label{lemma- absorption in Vn r}
Let $n > 2k$ and $r = n - 2k$.
\begin{enumerate}
\item  $\V n n \oplus \bmw n e_1 \bmw n = \bmw n$.
\item $\w {2k} 0 \otimes \bmw r \subseteq  \w n {r-2} +  \w {2k} 0 \otimes \V r 
r \subseteq   \w n {r-2} + \V n r$.
\item $ g_\pi ( \w {2k} 0 \otimes \bmw r)  g_\sigma \subseteq   \w n {r-2} + \V 
n r$ for any  positive permutation braids $g_\pi, g_\sigma $.
\end{enumerate}
\end{lemma}

{\em Proof}.\hods   For (1), we have to show that 
 $\chi g_\pi \in \V n n + \bmw n e_1 \bmw n$ for $\chi \in \{ e_i, g_i : 1 \le i 
\le n-1\}$.   We have
 $e_i g_\pi  \in \bmw n e_i \bmw n = \bmw n e_1 \bmw n$ for all~$i$.    If 
$\ell (s_i \pi) = \ell (\pi) + 1$, then
 $g_i g_\pi = g_{s_i \pi} \in \V n n$;  on the other hand, if $\ell (s_i \pi_i) = 
\ell (\pi) - 1$, then
 $g_\pi = g_i g_{s_i \pi}$, and 
$$
g_i g_\pi = (g_i^2) g_{s_i \pi}= g_{s_i \pi}  - z g_\pi + z \la\inv e_i 
g_{s_i \pi}  \in \V n n + \bmw n e_1 \bmw n.
$$

Applying (1) to $\bmw r$, we get $\bmw r = \bmw r e_1 \bmw r + \V r r$, so
$$
\w {2k} 0 \otimes \bmw r = \w {2k} 0 \otimes \bmw r e_1 \bmw r +  \w {2k} 0 
\otimes \V r r \subseteq\w n {r-2} +  \V n r.
$$
This shows (2).

Using Lemma \ref{lemma- factorization of positive permutation braid}, the 
positive permutation braid
$g_\pi$ can be written as $g_{\pi_1} w_1$, where $\pi_1$ is a $(2k, r)$-shuffle 
and
$w_1 \in \bmw {2k} \otimes \bmw {r}$.  Applying the same result to 
$\alpha(g_\sigma)$, we get
$g_\sigma = w_2 \alpha(g_{\sigma_1})$, where $\sigma_1$ is a $(2k, r)$-shuffle 
and
$w_2 \in  \bmw {2k} \otimes \bmw {r}$.  Thus
$$
g_\pi( \w {2k} 0 \otimes \bmw r)  g_\sigma  \subseteq g_{\pi_1} (\w {2k} 0 
\otimes \bmw r) \alpha(g_{\sigma_1}).
$$
By statement (2), this lies in 
$$
 \w n {r-2} +  g_{\pi_1} (\w {2k} 0 \otimes \V r r) \alpha(g_{\sigma_1}) 
\subseteq  \w n {r-2}  + \V n r.\squ
$$

\begin{lemma}\label{lemma- invariance of V n r}
$e_i \V n r \subseteq \w n {r-2} + \V n r$ and
$g_i \V n r \subseteq \w n {r-2} + \V n r$ for $1 \le i \le n-1$.
\end{lemma}

\pr
Consider $x = g_\pi [w \otimes g_\tau] \alpha(g_\sigma)$, where
$\pi$ and $\sigma$ are $(2k, r)$-shuffles, $w \in \w {2k} 0$, 
and $\tau \in \S_r$.  We have to show
that $e_i x$ and $g_i x$ lie in $\w n {r-2} + \V n r$ for all $1 \le i \le 
n-1$.

Suppose that $\pi\inv(i) > \pi\inv(i+1)$.  Then $g_\pi = g_i g_{\pi_1}$, where
${\pi_1}\inv(i) < \pi_1\inv(i+1)$, and $\pi_1$ is also a $(2k, r)$-shuffle.
Thus $e_i g_\pi = e_i g_i g_{\pi_1} = \la\inv e_i g_{\pi_1}$.  Likewise,
$g_i g_\pi = (g_i)^2 g_{\pi_1} = g_{\pi_1} -z g_\pi + z \la\inv e_i g_{\pi_1}$.
We are therefore reduced to considering the case that $\pi\inv(i) < 
\pi\inv(i+1)$.

If  $\pi\inv(i) < \pi\inv(i+1)$, then 
$$
g_i  x  = g_{s_i \pi} [w \otimes g_\tau] \alpha(g_\sigma) \subseteq \w n {r-2} 
+ \V n r,
$$
by Lemma \ref{lemma- absorption in Vn r}.  

It remains to consider
$e_i x$ when $\pi\inv(i) < \pi\inv(i+1)$.

If $\pi\inv(i+1) \le 2k$ or $\pi\inv(i) \ge 2k+1$, then 
$\pi\inv(i+1) = \pi\inv(i) + 1$.
By Lemmas 
\ref{lemma- transport gi by parallel strands}    and
\ref{lemma- absorption in Vn r},
$$
e_i x = g_\pi e_{\pi\inv(i)}[  w \otimes g_\tau] \alpha(g_\sigma) \in
g_\pi [ \w {2k} 0 \otimes  \bmw r ]  \alpha(g_\sigma)  \subseteq  \w n {r-2} + 
\V n r.
$$

The remaining case to consider is $\pi\inv(i) \le 2k$ and $\pi\inv(i+1) \ge 
2k+1$.
Define a permutation $\varrho$ by
$$
\varrho(j) = \cases{
j &if  $j < \pi\inv(i)$,\cr
j+1  &if $\pi\inv(i) \le j < 2k$,\cr
\pi\inv(i) & if $j = 2k$,\cr
\pi\inv(i+1) &if $j = 2k + 1$,\cr
j -1  &if $2k +1 < j \le \pi\inv(i+1)$,\cr
j & if $ j > \pi\inv(i+1)$.\cr}
$$
Since $\varrho \in \S_{2k} \times \S_r$, $\ell (\pi \varrho) = \ell (\pi) + 
\ell (\varrho)$ and
$g_{\pi \varrho} = g_\pi g_\varrho$.  The permutation $\pi \varrho$ has the 
following properties:
$\pi \varrho(2k) = i$;  $\pi \varrho(2k+1) = i+1$;  if $1 \le a < b \le 2k-1$ or 
$2k + 2 \le  a < b \le n$, then $\pi\varrho(a) < \pi\varrho(b)$.  We have
$$
e_i  g_\pi = e_i g_\pi  g_\varrho g_\varrho\inv = e_i g_{\pi \varrho} g_\varrho 
\inv
=  g_{\pi \varrho} e_{2k} g_\varrho \inv.
$$
Therefore
$$
e_i x \in  g_{\pi \varrho} e_{2k} (\w {2k} 0 \otimes \bmw r)  \alpha(g_\sigma).
$$
We concentrate on $e_{2k} (\w {2k} 0 \otimes \bmw r) $.

\begin{lemma}\label{lemma- factorization of w 2k 0}
 Suppose $\varphi$ is injective on $\w {2k} 0$.  Then
$\w {2k} 0$ is spanned by elements of the form
$$
g_m g_{m+1} \cdots g_{2k-2} e_{2k-1} w,
$$
where $1 \le m \le 2k-2$ and $w \in \w {2k} 0$.
\end{lemma}

\begin{proof} By assumption, $\varphi$ is an isomorphism from $\w {2k} 0$  to 
the
ideal $\k {2k} 0$ of
$\kt {2k}$ spanned by tangle diagrams  of rank~0. Now, $\k {2k} 0$ is spanned by 
totally descending tangle diagrams, and we can choose the order of the strands 
so that the strand incident with the vertex $\overline {\mathbold  {2k}}$ lies 
above all other strands.  Such a totally descending tangle diagram is isotopic 
to a tangle diagram of the following form:

\centerline{$T = \inlinegraphicee{w-2k-0-tangle}$}

\noindent Such a tangle diagram has the factorization $T = G_m \cdots G_{2k-2} 
E_{2k-1} 
T'$. 
 \end{proof}

\begin{lemma}\label{lemma- last MW lemma}
$g_{\pi \varrho} e_{2k} g_m g_{m+1} \cdots g_{2k-2} e_{2k-1} =
g_{\pi'} e_{2k-1}$ for some positive permutation braid $g_{\pi'}$.
\end{lemma}

\begin{proof}
First,  $e_{2k}$ can be moved to the right of $g_m g_{m+1} \cdots g_{2k-2}$ and
 $ e_{2k} e_{2k-1}$ can be written as $g_{2k-1} g_{2k} e_{2k-1}$ or as
$g_{2k-1}\inv g_{2k}\inv e_{2k-1}$.
Therefore
\begin{eqnarray*}
g_{\pi \varrho} e_{2k} g_m g_{m+1} \cdots g_{2k-2} e_{2k-1} 
&=& g_{\pi \varrho}  (g_m g_{m+1} \cdots g_{2k-2}) (g_{2k-1} g_{2k})e_{2k-1} \\
&=&g_{\pi \varrho}  (g_m g_{m+1} \cdots g_{2k-2}) (g_{2k-1}\inv g_{2k} \inv) 
e_{2k-1}.
\end{eqnarray*}
Since $\pi \varrho(a) < \pi \varrho(b)$ if $1 \le a < b \le 2k-1$, 
$g_{\pi \varrho}  g_m g_{m+1} \cdots g_{2k-2}$ is the  
positive permutation braid corresponding to the permutation
$\pi \varrho s_m \cdots s_{2k-2}$.
Moreover, 
$$
\begin{aligned}
\pi \varrho s_m \dots s_{2k-2} (2k-1) &= \pi\varrho (m), \\ 
\pi \varrho s_m \dots s_{2k-2} (2k) &= \pi \varrho(2k) = i,\\
\pi \varrho s_m \dots s_{2k-2} (2k+1) &= \pi \varrho(2k+1) = i+1.
\end{aligned}
$$
It follows that if $\pi \varrho(m) < i$, then 
$g_{\pi \varrho}  (g_m g_{m+1} \cdots g_{2k-2}) (g_{2k-1} g_{2k})$ is a positive
permutation braid, while if $\pi \varrho(m) > i+1$, then 
$g_{\pi \varrho}  (g_m g_{m+1} \cdots g_{2k-2})\cdot\break (g_{2k-1}\inv g_{2k} 
\inv)$ is 
a positive permutation braid.
In either case, the desired conclusion follows.
\end{proof}

{\em End of the proof of Lemma \ref{lemma- invariance of V n r}}.\hods  
We have to show that   
$$
g_{\pi \varrho} e_{2k} (\w {2k} 0 \otimes \bmw r)  
\alpha(g_\sigma)\subseteq \w n {r-2} + \V n r.
$$
By Lemma \ref{lemma- factorization of w 2k 0}, it is enough to show
$$
g_{\pi \varrho} e_{2k}(g_m \cdots g_{2k-2}) e_{2k-1} (\w {2k} 0 \otimes \bmw r)  
\alpha(g_\sigma)
\subseteq \w n {r-2} + \V n r.  
$$
But by Lemmas \ref{lemma- last MW lemma}
 and \ref{lemma- absorption in Vn r}, 
$$
\displaylines{\quad\
g_{\pi \varrho} e_{2k}(g_m \cdots g_{2k-2}) e_{2k-1} (\w {2k} 0 \otimes \bmw 
r)  \alpha(g_\sigma) \hfill\cr
\hfill =
g_{\pi'} e_{2k-1}  (\w {2k} 0 \otimes \bmw r)  \alpha(g_\sigma)
\subseteq \w n {r-2} + \V n r.\squ\quad\ \cr}
$$

\begin{proposition}\label{v n r complement of w n r-2}
 If $\varphi$ is injective on $\bmw {n-1}$ then
$$
\V n r \oplus \w n {r-2} = \w n r.
$$
\end{proposition}

\begin{proof} \hskip-3pt We have shown that $\V n r \!\oplus \w n {r-2}\!$ is an 
ideal 
containing
$e_1 \!\cdots e_{2k-1}$.  Therefore $\V n r \oplus \w n {r-2} = \w n r$.
\end{proof}

\begin{theorem}\label{theorem- isomorphism for ordinary BMW}
For all $n$, $\varphi : \bmw n \rightarrow \kt n$ is an isomorphism.
\end{theorem}

\begin{proof}  Surjectivity was shown in Theorem \ref{theorem- surjectivity of 
varphi on Wn}.  
Injectivity is evident for $n = 0, 1$.  We show $\varphi$ is injective on all 
ideals
$\w n r$ by double induction on $n$ and $r$.  By Proposition \ref{proposition- 
inductive step from w n-1 to w n 0}, if $\varphi$ is injective on
$\bmw {n-1}$, then $\varphi$ is injective on $\w n 0$ (if $n$ is even) or on
$\w n 1$ (if $n$ is odd).  By Corollary 
\ref{corollary- inductive step w n r-2 to v nr plus w n r-2} and 
Proposition \ref{v n r complement of w n r-2},  if
$\varphi$ is injective on $\bmw {n-1}$ and on $\w n {r-2}$, then
$\varphi$ is injective on $\w n r$.
\end{proof}

\begin{corollary}
For all $R$ and $n$,  $\bmw {n, R}$ is a free $R$-module with 
basis $\{\varphi\inv(T_c) : c \text{ is an $n$-connector}\}$.
\end{corollary}

\begin{corollary}\label{corollary- isomorphic specializations ordinary
BMW}
For all $R$ and $n$, $\bmw {n, R}
\cong\bmw {n,\varLambda}\otimes_{\varLambda} R$.
\end{corollary}

\begin{proof}  $\bmw  {n, R} \cong \kt {n, R}$ and 
 $\bmw  {n, \varLambda} \cong \kt {n, \varLambda}$.  But
 $\kt {n, R} \cong \kt {n, \varLambda} \otimes_{\varLambda} R$, by
 Corollary \ref{corollary-  freeness of ktn and isomorphism of specializations}.
\end{proof}

\begin{corollary}
If $R$ is a ring with distinguished elements $\lambda$, $z$ 
and $\delta$, and $S
\supseteq R$ is a ring containing $R$, then $\bmw {n, R} $ imbeds in $\bmw {n, 
S} $.
\end{corollary}

\begin{corollary}\label{corollary- imbedding of bmw in affine bmw}
 If $S$ is a ring with distinguished elements  $\lambda$, $z$, $\delta$, and 
$q_r$ ($r\ge 1$), then $i_n : \bmw {n, S} \rightarrow \abmw {n, S} $ 
is injective.
\end{corollary}

\begin{proof}  Consider the commuting diagram
\begin{diagram}
\abmw {n, S }      &\rTo^{\scriptstyle\varphi}    & \akt {n, S} \\
\uTo<{\scriptstyle i_n}                  &        & \uTo<{\scriptstyle i_n'}  \\
\bmw {n, S}                 &\rTo^{\scriptstyle\varphi}    & \kt {n, S} \\
\end{diagram}
Since $\varphi : \bmw {n, S} \rightarrow \kt {n, S}$ is an isomorphism, 
$i_n' = \varphi \circ i_n \circ
\varphi\inv$.  But $i_n'$ is injective by Remark \ref{remark-injections of 
tangle algebras in affine tangle algebras},
and therefore $i_n$ is  injective as well.
\end{proof}

\begin{corollary}
For all $R$ and $n$,
\begin{enumerate}
\item  There is a conditional expectation $\eps_n : \bmw {n, R} \rightarrow
\bmw {n-1, R}$ such that $e_n x e_n = \delta \eps_n(x)$ for $x \in \bmw {n, R}$.
\item The map $\iota : \bmw {n-1, R} \rightarrow \bmw {n, R}$ is injective
(and $\eps_n \circ \iota(x) = x$ for $x \in \bmw {n-1, R}$).
\item  There is a trace $\eps : \bmw {n, R} \rightarrow R$ such that
$\eps\circ \eps_{n+1} (x) = \eps(x)$  for $x \in \bmw {n+1, R}$ and
$\eps \circ \iota (x) = \eps(x)$ for $x \in  \bmw {n-1, R}$.
\end{enumerate}
\end{corollary}

\begin{proof}  This follows from the isomorphism $\bmw {n, R} 
\cong \kt {n, R}$ and the corresponding properties of $\kt {n, R}$.
\end{proof}

\ssection{Isomorphism of affine BMW}{and Kauffman tangle algebras}\label{sec6}
\setcounter{equation}{0}
\setcounter{theorem}{0}
\setcounter{section}{6} 
This section contains our main result, the isomorphism of $\abmw  {n, S}$ and
$\akt {n, S}$ for any ring $S$ with distinguished elements $\lambda$, $z$, 
$\delta$, and $q_r$ ($r\ge
1$).  We also give a basis of $\akt {n,S}$ over~$S$.

\subsection{Surjectivity of $\varphi: \abmw  {n}  \rightarrow \akt
{n}$.}
Fix a ring $S$ with distinguished elements $\lambda$, $z$, $\delta$, and $q_r$ 
($r \ge 1$).

\begin{definition}\rm
A {\em simple winding} is a piece of an affine tangle
diagram with one ordinary strand, without self-crossings, 
 regularly isotopic to the intersection of one of the
affine tangle diagrams
$X_1$ or $X_1\inv$ with a small neighborhood of the flagpole, as in the 
following figure:
\end{definition}

\centerline{$\inlinegraphic{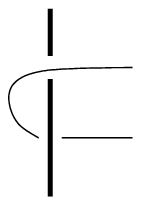}$ }

\begin{definition}\rm
An affine tangle diagram is in {\em standard
position}  if:
\begin{enumerate}
\item  It has no crossings to the left of the flagpole.
\item  There is a neighborhood of the flagpole whose intersection with
the tangle diagram is a union of simple windings.
\item  The simple windings have no crossings and are not nested.  That is,
between the two crossings of a simple winding with the flagpole, there is
no other crossing of a strand with the flagpole.
\end{enumerate}
\end{definition}

\centerline{$\inlinegraphicee{standard_position}$}

\begin{lemma}
Any affine tangle diagram is regularly isotopic to 
an affine tangle diagram in standard position.
\end{lemma}

{\em Proof}.\hods  An affine $(n,n)$-tangle can be viewed as the union of 
\begin{enumerate}
\item an ordinary tangle diagram  $T_{1}$with
$2n + 2 l $ boundary points to the right of the flagpole,
\item an ordinary tangle  diagram $T_{2}$ with $2l $ boundary points to the 
left of the flagpole, and
\item  $2 l $  horizontal strands connecting $T_{1}$ and $T_{2}$ and crossing 
above or below the
flagpole.
\end{enumerate}

\centerline{$\inlinegraphic{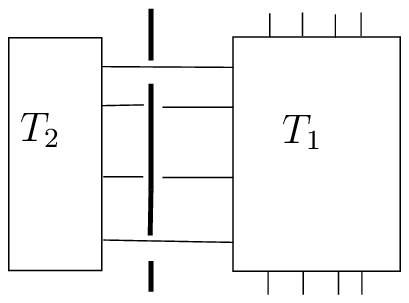}$}

From each horizontal strand crossing over the flagpole,  pull a finger over the 
tangle diagram~$T_{2}$.

\centerline{$\inlinegraphic{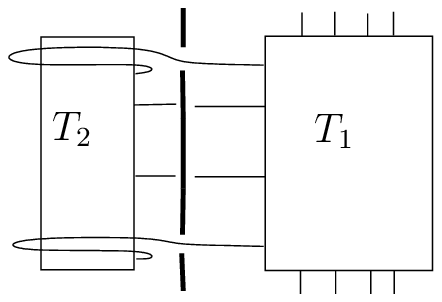}$}

\noindent Now slide $T_{1}$ and $T_{2}$ to the right, with the diagram $T_{2}$ 
sliding 
under the flagpole.  The result is a diagram in standard position.
$$
\displaylines{\hfill\inlinegraphic{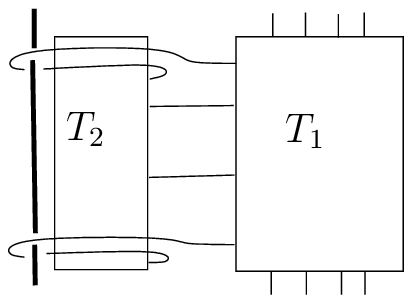}\hfill\cr
\noalign{\vskip-14pt}
\hfill\vrule height4pt width4pt depth0pt\cr}
$$

We give each affine tangle diagram in standard position a standard orientation
(see Section \ref{subsection-totally descending tangles}).   

In the oriented diagram, the simple windings are of the four possible types 
pictured below:
$$
\begin{aligned}[c]\setlength{\unitlength}{0.00062500in}
\begingroup\makeatletter\ifx\SetFigFont\undefined%
\gdef\SetFigFont#1#2#3#4#5{%
  \reset@font\fontsize{#1}{#2pt}%
  \fontfamily{#3}\fontseries{#4}\fontshape{#5}%
  \selectfont}%
\fi\endgroup%
{\renewcommand{\dashlinestretch}{30}
\begin{picture}(815,1281)(0,-10)
\thicklines
\path(278,783)(278,33)
\path(278,1233)(278,933)
\thinlines
\path(353,408)(803,408)
\path(803,858)(800,858)(794,858)
	(782,858)(765,858)(741,857)
	(712,857)(678,857)(641,856)
	(602,855)(563,855)(524,854)
	(486,853)(451,852)(417,850)
	(387,849)(358,848)(333,846)
	(309,844)(287,843)(267,841)
	(249,838)(232,836)(216,833)
	(193,828)(171,823)(151,817)
	(133,811)(115,804)(99,796)
	(84,787)(70,778)(58,768)
	(48,758)(39,748)(32,738)
	(26,727)(22,717)(18,706)
	(16,695)(13,681)(12,666)
	(13,649)(15,632)(18,615)
	(22,597)(28,580)(34,563)
	(41,548)(49,533)(57,520)
	(66,508)(74,497)(84,487)
	(95,477)(108,467)(123,456)
	(140,445)(158,434)(175,424)(203,408)
\path(122.998,419.163)(203.000,408.000)(152.766,471.258)
\end{picture}
} \\ {\rm(a)} \\ \end{aligned} 
\qquad \qquad
\begin{aligned}[c] \setlength{\unitlength}{0.00062500in}
\begingroup\makeatletter\ifx\SetFigFont\undefined%
\gdef\SetFigFont#1#2#3#4#5{%
  \reset@font\fontsize{#1}{#2pt}%
  \fontfamily{#3}\fontseries{#4}\fontshape{#5}%
  \selectfont}%
\fi\endgroup%
{\renewcommand{\dashlinestretch}{30}
\begin{picture}(815,1281)(0,-10)
\thicklines
\path(278,783)(278,33)
\path(278,1233)(278,933)
\thinlines
\path(803,858)(800,858)(794,858)
	(782,858)(765,858)(741,857)
	(712,857)(678,857)(641,856)
	(602,855)(563,855)(524,854)
	(486,853)(451,852)(417,850)
	(387,849)(358,848)(333,846)
	(309,844)(287,843)(267,841)
	(249,838)(232,836)(216,833)
	(193,828)(171,823)(151,817)
	(133,811)(115,804)(99,796)
	(84,787)(70,778)(58,768)
	(48,758)(39,748)(32,738)
	(26,727)(22,717)(18,706)
	(16,695)(13,681)(12,666)
	(13,649)(15,632)(18,615)
	(22,597)(28,580)(34,563)
	(41,548)(49,533)(57,520)
	(66,508)(74,497)(84,487)
	(95,477)(108,467)(123,456)
	(140,445)(158,434)(175,424)
	(189,416)(198,411)(202,408)(203,408)
\path(473.000,438.000)(353.000,408.000)(473.000,378.000)
\path(353,408)(803,408)
\end{picture}
} \\ {\rm(b)} \\ 
\end{aligned}
\qquad \qquad
\begin{aligned}[c]\setlength{\unitlength}{0.00062500in}
\begingroup\makeatletter\ifx\SetFigFont\undefined%
\gdef\SetFigFont#1#2#3#4#5{%
  \reset@font\fontsize{#1}{#2pt}%
  \fontfamily{#3}\fontseries{#4}\fontshape{#5}%
  \selectfont}%
\fi\endgroup%
{\renewcommand{\dashlinestretch}{30}
\begin{picture}(692,1281)(0,-10)
\thicklines
\path(230,258)(230,33)
\path(230,1233)(230,408)
\thinlines
\path(155,783)(154,783)(150,780)
	(141,775)(128,766)(111,756)
	(94,744)(78,732)(64,720)
	(53,708)(44,697)(36,684)
	(30,670)(25,658)(21,644)
	(18,628)(15,612)(14,595)
	(12,577)(12,558)(12,539)
	(14,521)(15,504)(18,488)
	(21,472)(25,458)(30,445)
	(35,435)(40,425)(46,416)
	(53,407)(61,399)(70,391)
	(81,384)(92,377)(105,371)
	(118,366)(133,361)(148,357)
	(164,354)(181,350)(199,348)
	(218,345)(233,344)(250,342)
	(269,341)(289,340)(311,339)
	(336,338)(364,337)(394,336)
	(427,336)(462,335)(499,335)
	(536,334)(571,334)(604,334)
	(632,333)(653,333)(668,333)
	(676,333)(679,333)(680,333)
\path(680,783)(305,783)
\path(395.000,813.000)(305.000,783.000)(395.000,753.000)
\end{picture}
} \\ {\rm(c)} \\ 
\end{aligned} 
\qquad \qquad
\begin{aligned}[c]\setlength{\unitlength}{0.00062500in}
\begingroup\makeatletter\ifx\SetFigFont\undefined%
\gdef\SetFigFont#1#2#3#4#5{%
  \reset@font\fontsize{#1}{#2pt}%
  \fontfamily{#3}\fontseries{#4}\fontshape{#5}%
  \selectfont}%
\fi\endgroup%
{\renewcommand{\dashlinestretch}{30}
\begin{picture}(692,1281)(0,-10)
\thicklines
\path(230,258)(230,33)
\path(230,1233)(230,408)
\thinlines
\path(680,783)(305,783)
\path(69.436,693.675)(155.000,783.000)(37.468,744.449)
\path(155,783)(128,766)(111,756)
	(94,744)(78,732)(64,720)
	(53,708)(44,697)(36,684)
	(30,670)(25,658)(21,644)
	(18,628)(15,612)(14,595)
	(12,577)(12,558)(12,539)
	(14,521)(15,504)(18,488)
	(21,472)(25,458)(30,445)
	(35,435)(40,425)(46,416)
	(53,407)(61,399)(70,391)
	(81,384)(92,377)(105,371)
	(118,366)(133,361)(148,357)
	(164,354)(181,350)(199,348)
	(218,345)(233,344)(250,342)
	(269,341)(289,340)(311,339)
	(336,338)(364,337)(394,336)
	(427,336)(462,335)(499,335)
	(536,334)(571,334)(604,334)
	(632,333)(653,333)(668,333)
	(676,333)(679,333)(680,333)
\end{picture}
} \\ {\rm(d)} \\ 
\end{aligned}.
$$

\begin{lemma}\label{lemma-affine surjectivity}
Let $T$ be an affine tangle diagram  
in standard position
and with a standard orientation.  Suppose T has
$l $ simple windings and
$K$ crossings of ordinary strands.    If $T$ has a simple winding on a 
non-closed strand, then $T$ can be
written as
$T =  A T_0 + B$  or $T = T_0 A + B$, where $A$ is in the monoid generated by 
$\{X_{1} \pmone,
G_{i}\pmone\}$,  $T_0$  is an affine tangle diagram in standard position with 
$l-1$ simple windings, and
$B$ is in the span of affine tangle diagrams  in standard position with fewer 
than $K$ crossings of
ordinary strands and at most $l $ simple windings.
\end{lemma}

\begin{proof}  Let $s$ be a non-closed strand of $T$ with a simple winding.
Suppose that the initial endpoint of $s$ is on the top boundary of  $R$, 
say the $j$th vertex on the top  boundary.  (Otherwise, the final endpoint
of $s$ is on the bottom boundary of $R$, and this case can be handled
similarly.)   Consider the first simple winding on~$s$.  

If the simple winding is of type (a) or (d), then consider the affine tangle 
diagram $T'$ obtained from
$T$ by reversing  crossings of $s$ as necessary so that all crossings of $s$ 
with other ordinary strands
are over-crossings, and all self-crossings of $s$ are encountered first as 
over-crossings. By the
Kauffman skein relation, $T$ and $T'$ are congruent modulo the span of affine 
tangle diagrams in standard
position with fewer than
$K$ crossings (and at most $l $ simple windings).   Consider the diagram
$T_0 = T' G_{j}\inv \cdots G_{1}\inv X_{1} \inv   $. It has two simple 
windings, one of type (c) and 
the other of type (a) or (d), connected by an unknotted arc which has only 
over-crossings with other
arcs;  thus $T_0$ is  regularly isotopic to an affine tangle diagram in standard 
position with $l  -1$ simple windings:
$$
\inlinegraphic{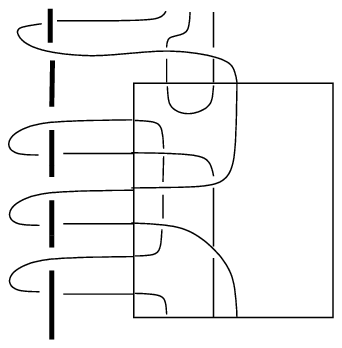}
$$

We have $T = T' + B = T_0 X_{1} G_{1} \cdots G_{j}  + B$, where $B$ is in the 
span of affine tangle
diagrams in standard position with fewer than 
$K$ crossings and at most $\ell $ simple windings.

If the first simple winding on $s$ is of type (b) or (c), then take $T'$ to be 
the diagram obtained from
$T$ by changing crossings of $s$ so that  all crossings of~$s$ with other 
ordinary strands are
under-crossings, and all self-crossings of $s$ are encountered first as 
under-crossings, and take
$T_0 = T' G_{j} \cdots G_{1} X_{1} $.   Now $T_0$ has two simple windings, one 
of type (a) and the other of
type (b) or (c), connected by an unknotted arc which has only  under-crossings 
with other arcs;  again, 
$T_0$ is  regularly isotopic to an affine tangle diagram in standard position 
with $l-1$ simple
windings.
\end{proof}

Say that the oriented affine tangle diagram $T$ is {\em totally descending} if, 
as  $T$ is traversed in
accordance with the orientation, each crossing of ordinary strands is 
encountered the first time as an
over-crossing.  By the inductive argument of Proposition 
\ref{proposition-totally descending tangles
span}, 
$\akt {n}$ is spanned by totally descending affine tangle diagrams in standard 
position.

\begin{theorem}\label{theorem-surjectvity of varphi on affine Wn}
$\varphi : \abmw  {n, S} \rightarrow \akt {n, S}$ is surjective.
 \end{theorem}

\begin{proof}  We have to show that each totally descending 
oriented affine\break tangle diagram $T$ in standard position 
 is in the subalgebra of $\akt {n}$ generated by $X_{1}\pmone$ and $\{E_{i}, 
G_{i}\pmone\}$.

The proof is by induction 
 first on the number $l $ of simple windings of $T$ on non-closed strands, and
then on the number $K$  of crossings of ordinary strands of $T$.

Suppose that $T$ has no simple windings on non-closed strands.
 Because $T$ is totally descending and all simple windings are located on closed 
loops of~$T$,
 $T$ is regularly isotopic to  a diagram in which the closed loops are confined 
to a neighborhood of the
flagpole and the non-closed strands do not enter this neighborhood.  It follows 
that $T \in \akt {0, S}
\otimes KT_{n, S}$, which is isomorphic to $KT_{n, S}$ by Turaev's theorem 
\cite{Turaev-Kauffman-skein}.
By Theorem~\ref{theorem- surjectivity of varphi on Wn},  $T$ is in the 
subalgebra generated by $\{E_{i},
G_{i}\pmone\}$.
 
Suppose $T$ has
$l  \ge 1$ simple windings on non-closed strands and that any totally 
descending affine tangle diagram
in standard position with fewer than $l $ simple windings  on non-closed 
strands is in the subalgebra
of $\akt {n, S}$ generated by $X_{1}\pmone$ and $\{E_{i}, G_{i}\pmone\}$. If the 
number $K$ of crossings
of ordinary strands is zero, then by Lemma \ref{lemma-affine surjectivity},
$T = A T_0$ or $T = T_0 A$, where $A$ is in the monoid generated by $\{X_{1} 
\pmone, G_{i}\pmone\}$,  and
$T_0$  is an affine tangle diagram in standard position with $l-1$ simple 
windings.  By the induction
assumption on $l $, the affine tangle  $T_0$ and hence $T$ is in the 
subalgebra of $\akt {n, S}$
generated by $X_{1}\pmone$ and $\{E_{i}, G_{i}\pmone\}$.   

Now suppose $T$ has $l \ge 1$ simple windings on non-closed strands and $K 
\ge 1$ crossings of ordinary
strands.    Then, by Lemma \ref{lemma-affine surjectivity}, $T$ can be written 
as
$T = A T_0 + B$ or $T = T_0 A + B$ where $B$  is has fewer than $K$ crossings, 
$T_0$ has fewer than
$l $ simple windings, and $A$ is in the monoid generated by  $\{X_{1} \pmone, 
G_{i}\pmone\}$.  It
follows from the induction assumptions on $l $ and $K$ that $T$  is in the 
subalgebra of $\akt {n}$
generated by $X_{1}\pmone$ and
$\{E_{i}, G_{i}\pmone\}$.
\end{proof}

\subsection{Freeness of  $\abmw  n$ and 
injectivity of $\varphi: \abmw  n \rightarrow \akt n$.}
Fix a ring $S$  with distinguished elements $\lambda$, $z$, $\delta$, and $q_r$ 
($r \ge 1$).  Write $\bmw n$ for
$\bmw {n, S}$, $\kt n$ for 
$\kt {n, S}$, $\abmw n$ for $\abmw {n, S}$, and $\akt n$ for
$\akt {n, S}$.

We can identify $\bmw n$ with $\kt n$ via $\varphi$, and we can regard $\bmw n$ 
as imbedded in $\abmw n$ by
Corollary~\ref{corollary- imbedding of bmw in affine bmw} and $\kt n$ imbedded 
in $\akt n$ by Remark
\ref{remark-injections of tangle algebras in affine tangle algebras}.

\begin{lemma}\label{lemma- elementary transport lemma1}
Let $T' \in \kt n$ be a totally descending $(n,n)$-tangle diagram
whose strands have no self-crossings.  Suppose that $T'$ has a strand 
connecting a top vertex $\p a$ with a bottom vertex $\pbar b$, and that this 
strand lies above all other strands (i.e. has only over-crossings with other 
strands).  Let\/ $T''$ be the otherwise identical tangle diagram in which the 
strand connecting $\p a$ and~$\pbar b$ lies under all other strands. Then
$$
x_b \varphi\inv(T') = \varphi\inv(T'') x_a.
$$
\end{lemma}

\begin{proof}  If $T'$ is a braid diagram, then the assertion follows from
Propositions \ref{proposition-xr relations1}  and \ref{proposition-xr 
relations2}
 and induction on the number of crossings.  
 
Otherwise, $T'$ is regularly isotopic to a tangle diagram of the form
 $$
 T' = B_1[ (E_1 E_3 \cdots E_{2k-1}) \otimes B] B_2,
 $$
where $B_1$, $B_2$, and $B$ are totally descending braid diagrams,
 as described in the proof of  Theorem \ref{theorem- surjectivity of varphi on 
Wn}.
The
 $x_*$ can be successively passed through $\varphi\inv(B_1)$, 
 $\varphi\inv(B)$, and $\varphi\inv(B_2)$, by the result for braid diagrams.
\end{proof}

\begin{remark}\rm
Applying $\varphi$, we get 
$$
X_b T' = T'' X_a.
$$
This statement can also be verified  by the following picture proof (which 
provides the motivation for Lemma \ref{lemma- elementary transport
lemma1}):
$$
\begin{aligned}
\inlinegraphic{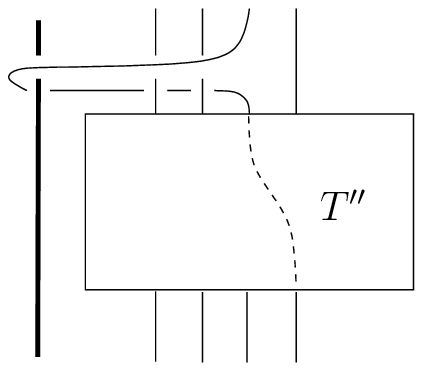} \quad  &\leftrightarrow 
\quad\inlinegraphic{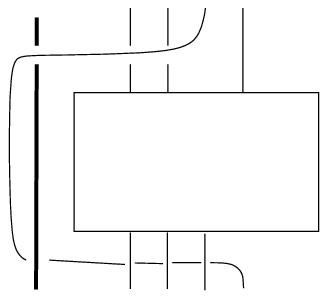}\quad \\
&\leftrightarrow  \quad \inlinegraphic{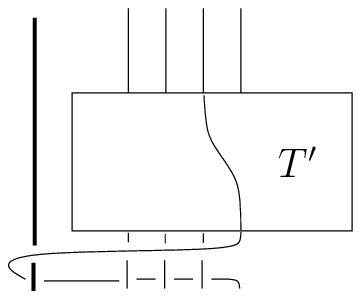}
\end{aligned}
$$
\end{remark}

\begin{lemma}\label{lemma- elementary transport lemma2}
Let $T' \in \kt n$ be a totally descending $(n,n)$-tangle diagram
whose strands have no self-crossings.  Suppose that $T'$ has a strand 
connecting a bottom vertex $\pbar a$ with a bottom vertex $\pbar b$, and that 
this strand lies above all other strands (i.e. has only over-crossings with 
other strands).  Let\/ $T''$ be the otherwise identical tangle diagram in which 
the strand connecting $\pbar a$ and $\pbar b$ lies under all other strands. Then
$$
x_a \varphi\inv(T') =\la^{-2} x_b\inv  \varphi\inv(T'') .
$$
\end{lemma}

{\em Proof}.\hods  Assume without loss of generality that $a < b$.   As in the 
proof 
of
Lemma \ref{lemma- factorization of w 2k 0}, $T'$ has a factorization $T' = G_a 
G_{a+1} \cdots G_{b-2} E_{b-1} T_0$.  Thus
$$
\begin{aligned}
x_a \varphi\inv(T') &= x_a g_a \cdots g_{b-2} e_{b-1} \varphi\inv(T_0)\\ 
&= g_a\inv \cdots g_{b-2}\inv x_{b-1} e_{b-1} \varphi\inv(T_0)\\
&= \la^{-2}g_a\inv \cdots g_{b-2}\inv x_{b}\inv e_{b-1}\varphi\inv(T_0)\\  
&=  \la^{-2}x_{b}\inv g_a\inv \cdots g_{b-2}\inv  e_{b-1}\varphi\inv(T_0)  \\
&=  \la^{-2}x_{b}\inv \varphi\inv(T'').\squ
\end{aligned}
$$

\begin{remark}\rm
Applying $\varphi$, we get 
$$
X_a T' = \la^{-2} X_b\inv  T''.
$$
This statement can also be verified   by the following picture proof (which 
provides the motivation for Lemma \ref{lemma- elementary transport
lemma2}):
$$
\begin{aligned}
\inlinegraphic{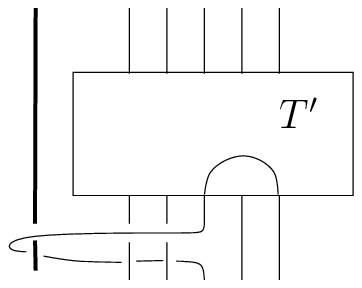} \quad  &\leftrightarrow \quad \lambda\inv 
\inlinegraphic{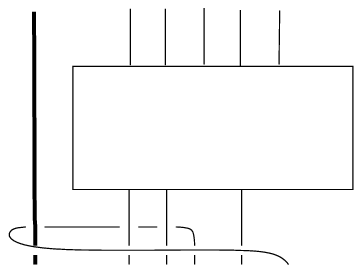}\quad\\ &\leftrightarrow \quad \lambda^{-2} 
\inlinegraphic{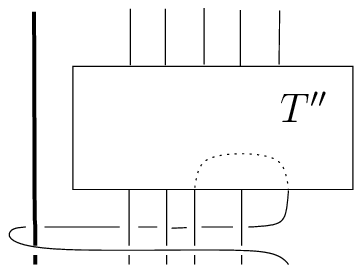}
\end{aligned}
$$
\end{remark}

In the following we write $x^\mu = x_1^{\mu_1} \cdots x_n^{\mu_n}$ for a Laurent 
monomial in the commuting elements  $x_i$.

\begin{lemma}[Transport Lemma]\label{lemma-transport}
Let $T \in \kt n$ be an ordinary  $(n,n)$-tangle diagram.
\begin{enumerate}
\item Suppose $T$ has a strand connecting a top vertex $\p a$ with a bottom 
vertex $\pbar b$.
Then 
$$\qquad\ 
x_b \varphi\inv(T) \equiv \varphi\inv(T) x_a,  \quad\ 
x_b\inv \varphi\inv(T) \equiv \varphi\inv(T) x_a\inv  
$$
modulo the span of elements of the form
$
x^\mu \varphi\inv(T_0)x^\nu, 
$
where $T_0$ is an ordinary $(n,n)$-tangle 
diagram with fewer crossings than~$T$.
\item Suppose $T$ has a strand connecting two  bottom vertices $\pbar a$ and 
$\pbar b$.
Then 
$$\qquad\ 
x_a \varphi\inv(T)  \equiv \la^{-2} x_b\inv \varphi\inv(T) 
$$
 modulo the 
span of elements  of the form 
$
x^\mu \varphi\inv(T_0)x^\nu,
$
where $T_0$ is an ordinary $(n,n)$-tangle diagram with fewer crossings than $T$.
\item Suppose $T$ has a strand connecting two  top vertices $\p a$ and $\p b$.
Then 
$$\qquad\ 
\varphi\inv(T) x_a \equiv \la^{-2} \varphi\inv(T) x_b\inv
$$
modulo the  span of elements  of the form
$
x^\mu \varphi\inv(T_0)x^\nu,
$
where $T_0$ is an ordinary $(n,n)$-tangle diagram with fewer crossings than $T$.
\end{enumerate}
\end{lemma}

\begin{proof}  By the Kauffman skein relation,   if $T'$ is an ordinary 
$(n,n)$-tangle diagram which differs from $T$ by changing one or more crossings, 
then 
$T$ and $T'$ are congruent modulo the span of 
$(n,n)$-tangle diagrams with fewer crossings. 

 Let $T'$ be a diagram obtained from $T$ by changing crossings,  in which the 
strand connecting $\p a$ and $\pbar b$ has only over-crossings with other 
strands,
and let $T''$ be the otherwise identical diagram in which 
 the strand connecting $\p a$ and $\pbar b$ has only under-crossings with other 
strands.
Then $T \equiv T' \equiv T''$ modulo the span of tangle diagrams with fewer 
crossings.
By Lemma \ref{lemma- elementary transport lemma1}, $x_b \varphi\inv(T') 
=\varphi\inv( T'' )x_a$, and  the first assertion in statement (1) follows.  The 
second assertion results from multiplying the congruence on the left 
by~$x_b\inv$ and on the right by $x_a\inv$.

Statement (2) follows similarly using Lemma \ref{lemma- elementary transport 
lemma1}, and statement (3) is obtained from (2) by applying the anti-automorphism 
$\alpha$.
\end{proof}

\begin{proposition}\label{proposition-XTX ideal}
$\abmw  {n}$ is spanned by 
 elements of the form $x^\mu  w x^\nu$, where
$x^\mu$, $x^\nu$ are Laurent monomials in the elements $x_i$ and  
$w \in \bmw n$.
\end{proposition}

\begin{proof}  Let $M$ denote the span of elements of the form $x^\mu w x^\nu$ .
Because $M$ contains the identity element, it suffices to show that $M$ is a 
left ideal in~$\abmw  {n}$.

Since $\abmw  {n}$ is generated as a unital  algebra by
$x_{1}\pmone$ and $\{g_{i}\pmone, e_{i}\}$,  it suffices to show that $M$ is 
invariant under left multiplication by these  generators.  Invariance under left 
multiplication by $x_1\pmone$ is obvious (since the
$x_i$'s commute).    Because of Propositions~\ref{proposition-xr relations1}
and \ref{proposition-xr relations2},   invariance under left multiplication by
$g_i\pmone$ will follow from invariance under left multiplication by~$e_i$.

For $k \ge 0$, let $M_k$ be the span of elements of the form $x^\mu 
\varphi\inv(T) x^\nu$, where
$x^\mu$, $x^\nu$ are monomials in the elements $x_i$ and $T$ is an ordinary 
$(n,n)$-tangle diagram with
no more than $k$ crossings.    Let $M_{-1} = (0)$.
We show by induction first on $i$ and then on $k$ that $ e_{i }M_k \subseteq M$.

Suppose $T$ is an ordinary  $(n,n)$-tangle diagram with no more than $k$ 
crossings and 
$\pbar b$ and $\p a$ are connected by a strand of $T$.  Then, by Lemma 
\ref{lemma-transport},
$x_b \varphi\inv(T) \equiv \varphi\inv(T) x_a$ and 
$x_b \inv \varphi\inv(T) \equiv \varphi\inv(T) x_a \inv $ modulo $M_{k-1}$.  
In particular, if $k = 0$, then $x_b \varphi\inv(T)  = \varphi\inv(T) x_a$  and   
$x_b\inv \varphi\inv(T)  = \varphi\inv(T) x_a\inv$.
Likewise, if  $\pbar b$ and $\pbar a$ are connected by a strand of $T$, then
$x_b \varphi\inv(T) \equiv \la^{-2}x_a \inv \varphi\inv(T) $ modulo $M_{k-1}$, 
and in particular, if $k = 0$, then
$x_b \varphi\inv(T) = \la^{-2}x_a \inv \varphi\inv(T) $.

It follows from these observations that, for any $i$, 
$$
x^{\mu} \varphi\inv(T) \equiv  \varrho x^{\mu'} \varphi\inv(T) x^{\nu'}\bmod 
M_{k -1},
$$
where $\varrho \in S$ and
\begin{enumerate}
\item if $\pbar i$ and $\overline  { \mathbold {i}+\mathbold{1}}$ are not
connected by a strand of $T$,  then  $\mu'_{i} = \mu'_{i+1} = 0$,
\item  if 
$\pbar i$ and $\overline  { \mathbold {i}+\mathbold{1}}$ are connected by a strand 
of $T$, 
then
 $\mu'_{i+1}  = 0$.
\end{enumerate}
   
For all $i$, if $\pbar i$ and $\overline  { \mathbold {i}+\mathbold{1}}$ are not 
connected 
by a strand of $T$, then
\begin{equation}\label{equation: XTX ideal 1}
e_{i} x^{\mu} \varphi\inv(T)  \equiv  \varrho e_{i} x^{\mu'} \varphi\inv(T) 
x^{\nu'} =  \varrho x^{\mu'} e_{i}
\varphi\inv(T) x^{\nu'}   \bmod   e_{i} M_{k -1}.
\end{equation}
Suppose $\pbar i$ and $\overline  { \mathbold {i}+\mathbold{1}}$ are  connected by 
a strand 
$s$ of $T$.  We have
\begin{eqnarray} 
\qquad\ e_{i} x^{\mu} \varphi\inv(T)  \equiv   \varrho  e_{i} x^{\mu'} x_{i}^{r} 
\varphi\inv(T) \bmod e_{i} M_{k-1}
=   \varrho  x^{\mu'}  e_{i}  x_{i}^{r} \varphi\inv(T), 
 \end{eqnarray}
where $\mu'_{i}= \mu'_{i+1} = 0$. Modulo the span of diagrams with fewer 
crossings, $T$ is congruent to a
totally descending tangle diagram $T'$.  Moreover,  $T' = E_{i} T''$  (where 
$T''$ may have more than $k$
crossings).  In fact,  $T'$ must have a strand  $s'$ connecting two points on 
the top boundary of $\R
\times I$.   By regular isotopy, $s$~can be contracted and $s'$ pulled down so 
that 
the minimum on $s'$ lies
just above the maximum on $s$ and both lie below any other crossing or local 
extremum of the tangle
diagram.   Thus
\begin{eqnarray}\label{equation: XTX ideal 3}
e_{i} x^{\mu} \varphi\inv(T)  &\equiv&   \varrho  x^{\mu'} e_{i}  x_{i}^{r} 
e_{i}\varphi\inv(T'') \bmod e_{i} M_{k-1}\\
 &=&   \varrho  x^{\mu'}  b e_{i} \varphi\inv(T''),\nonumber
\end{eqnarray}
where $b \in \abmw  {n-1}$, by Proposition \ref{proposition- bimodule 3}.
It follows from equations (\ref{equation: XTX ideal 1}) and 
(\ref{equation: XTX ideal 3}) that for all $i$ and~$k$,
\begin{equation}\label{equation: XTX ideal 4}
e_{i}M_{k} \subseteq M_{k} + e_{i}M_{k-1} +  A  \abmw  {i-1} M,
\end{equation}
where $A$ is the algebra  of Laurent polynomials in  $\{x_j : 1 \le j \le n\}$.
For $i = 1$, we have 
$$
e_{1}M_{k}  \subseteq M_{k} + e_{1} M_{k-1}.
$$
It follows by induction on $k$ that $e_{1} M_{k} \subseteq M_{k}$ for all $k$; 
that is, 
$e_{1} M \subseteq M$.
Now fix $i > 1$ and assume  that $e_{j}  M  \subseteq M$ for all
$j < i$; it follows that $\abmw  {i-1}  M \subseteq M$, and thus from
(\ref{equation: XTX ideal 4}), for all~$k$,
$$
e_{i}M_{k} \subseteq M + e_{i}M_{k-1}.
$$  
By induction on $k$,  $e_{i}M_{k} \subseteq M$ for all $k$.
\end{proof}

Let $B_{0} = \{T_{d} : \text{$d$ is an $n$-connector}\}$ be the common  basis of 
$\kt {n, R}$ for all $R$, as
described at the end of Section~\ref{subsection-totally descending tangles}.
Let $B$ be the set of elements
$X^{\mu} T X^{\nu}$, where $T \in B_{0}$ 
and $\nu_{i} = 0$ unless $T$ has a strand with initial
point~$\p i$, and $\mu_{j} = 0$ unless  $T$ has a strand with initial
point $\pbar j$.     The set $B$ may be regarded as a subset of
$\akt {n, S}$ for any $S$. 

For any $S$, we lift $B$ to a subset of $\abmw {n, S}$ as follows:
$\varphi : \bmw {n, S} \rightarrow \kt {n, S}$ is an isomorphism.
We take $A = A(n, S)$ to be the set of elements of $\abmw {n, S}$ of the form
$x^{\mu} \varphi\inv(T) x^{\nu}$, where $T \in B_{0}$ 
and $\nu_{i} = 0$ unless $T$ has a strand with initial
point $\p i$, and $\mu_{j} = 0$ unless  $T$ has a strand with initial
point~$\pbar j$.  

\begin{proposition}\label{proposition- A n R spans}
For every $S$ and $n$, 
\begin{enumerate}
\item$A(n, S)$ spans $\abmw {n, S}$.
\item $B$ spans
$\akt {n, S}$.
\end{enumerate}
\end{proposition}

\begin{proof}   Because  of Proposition \ref{proposition-XTX ideal},  to prove 
statement (1) it suffices
to show that every element of the form $x^\mu \varphi\inv(T) x^\nu$, where $T$ 
is an ordinary
$(n,n)$-tangle diagram, is in the span of $A = A(n, S)$.    The basis $B_0$ has 
the following triangular
property:  any ordinary $(n,n)$-tangle diagram is congruent to a multiple of an 
element of $B_0$ modulo
the span of diagrams with fewer crossings. Combining this with the Transport 
Lemma
\ref{lemma-transport}, we can prove that
$x^\mu \varphi\inv(T) x^\nu$ is in the span of $A$ by induction on the number of 
crossings of~$T$.  

Statement (2) follows from statement (1) by applying the surjective map 
$\varphi : \abmw {n, S} \rightarrow \akt {n, S}$.
\end{proof}

\begin{theorem}\label{theorem- main result}
For all $n$,  
\begin{enumerate}
\item
$\abmw {n, \varLambdaHat }$ is  a free $\varLambdaHat$-module
with basis~$A$.
\item $\akt {n, \varLambdaHat }$ is  a free $\varLambdaHat$-module
with basis~$B$.
\item
$\varphi: \abmw  {n, \varLambdaHat}  \rightarrow \akt {n, \varLambdaHat}$ is 
an isomorphism.
\end{enumerate}
\end{theorem} 

\begin{proof}  $ A = A(n, \varLambdaHat)$ spans $\abmw  {n, \varLambdaHat}$,
and $B = \varphi(A)$ spans $\akt {n, \varLambdaHat}$,  by Proposition 
\ref{proposition- A n R spans}.   Since  $c$ is
injective on $B$,  Proposition \ref{proposition- distinct affine connectors 
linear independence} implies that
$B$ is linearly independent over $\varLambdaHat$.  All the conclusions follow.
\end{proof}

\begin{corollary}\label{corollary- freeness affine specializations}
\mbox{}
\begin{enumerate}
\item
$\akt {n, S}$  is a free $S$ module with basis
$B$ and $\akt {n, S} \cong \akt {n, \varLambdaHat} \otimes_{\varLambdaHat} S$.
\item  $\varphi : \abmw {n, S} \rightarrow \akt {n, S}$ is an isomorphism.
\item  $\abmw {n, S} \cong \abmw {n, \varLambdaHat} \otimes_{\varLambdaHat} S$.
\end{enumerate}
\end{corollary}

{\em Proof}.\hods   By Proposition \ref{proposition- A n R spans}, $B$ spans 
$\akt {n, S}$.  On the other hand, $B \otimes 1$ is a basis of
$ \akt {n, \varLambdaHat} \otimes_{\varLambdaHat} S$.  Moreover, we have an 
$S$-algebra homomorphism
from $\akt {n, S}$ to $ \akt {n, \varLambdaHat} \otimes_{\varLambdaHat} S$ which 
takes a tangle $T$ to $T \otimes 1$.
Since this map takes a spanning set to a basis, it is an isomorphism.  This 
proves~(1).  

Once we know that $B$ is a basis of $\akt {n, S}$,  it follows at once that  
$\varphi : \abmw {n, S}
\rightarrow \akt {n, S}$ is an isomorphism.  Namely, $\varphi$ carries the 
spanning set $A(n, S)$ to the basis
$B$, and so  is an isomorphism.

Finally, (3) follows by combining the isomorphisms 
$$
\abmw {n, S} \cong \akt {n, S},   \quad\ 
\akt  {n, S} \cong  \akt {n, \varLambdaHat} \otimes_{\varLambdaHat} S,
\quad\ 
\abmw {n, \varLambdaHat} \otimes_{\varLambdaHat} S \cong \akt {n, \varLambdaHat} 
\otimes_{\varLambdaHat} S.\squ
$$

\begin{corollary}\label{corollary- injectivity of 
iota from affine bmw n-1 to affine bmw n}
The map $\iota: \abmw {n-1, S}\rightarrow \abmw   {n, S}$ is injective.
\end{corollary}

\begin{corollary}\label{corollary- conditional expectation  and trace for 
affine BMW}
\mbox{}
\begin{enumerate}
\item  For each $n\ge 1$, there is a  conditional expectation $\eps_n : \abmw      
{n, S} \rightarrow \abmw  {n-1,S}$ satisfying
$$
e_n x e_n = \delta \eps_n(x) e_n
$$
for $x \in \abmw  {n, S}$.
\item  There is a trace $\eps: \abmw     {n, S}  \rightarrow S$ defined by
$$
\eps = \eps_1 \circ \cdots \circ \eps_n.
$$
\item
The trace $\eps$ has the Markov property:  for $b \in \abmw {n-1, S}$,
\begin{enumerate}
\item[\rm(a)] $\eps(b g_{n-1}^{\pm 1} ) = (\la^{\pm 1}/\delta) \eps(x) $,
\item[\rm(b)] $\eps(b e_{n-1}) = (1/\delta) \eps(x)$,

\item[\rm(c)\hskip1.2pt]  $\eps(b (x'_n)^r) = q_{r} \eps(b)$, and
$\eps(b (x'_n)^{-r}) = f_{r}(q_1, \dots, q_r) \eps(b)$,
\end{enumerate}
for $r \ge 1$, 
where $x'_n = (g_{n-1} \cdots g_1) x_1 (g_1\inv \cdots g_{n-1}\inv)$.

\end{enumerate}
\end{corollary}

\begin{proof} All statements follow from the isomorphism $\abmw  {n, S} \cong 
\akt n$ and the properties
of $\akt {n, S}$  discussed in Section \ref{subsection- inclusions and trace}.  
\end{proof}

\begin{remark}\rm
The trace $\eps$ is the unique  trace on $\abmw {n, S}$ 
satisfying $\eps(1) \!=\! 1$ and the Markov properties
enumerated in statement (3) of the previous corollary.  To prove this uniqueness 
statement, one needs an analogue
of   Proposition \ref{proposition- bimodule 2} with $x_n$ replaced by $x_n'$.
\end{remark}

\begin{corollary}
For $x \in \abmw  {n+1}$, 
$$
x e_n =\delta \eps_{n+1}(x e_n) e_n.
$$
\end{corollary}

\begin{proof}  By Proposition \ref{proposition- bimodule 5},  there exists a $y 
\in \abmw  n$ such that
$x e_n = y e_n$.   Applying the conditional expectation $\eps_{n+1}$ to both 
sides gives
$\eps_{n+1}(x e_n) = y \eps_{n+1}(e_n) = \delta\inv y$.
\end{proof}

\begin{corollary}
Every element of $\akt {n, S}$ is a linear combination of 
elements of the  form $a \chi
b$, where $a, b \in \akt {n-1,S}$ and 
$$
\chi  \in \{E_{n-1}, G_{n-1}\pmone\} \cup \{X_{n}^{r} : r \in \Z\}.
$$
\end{corollary}

\begin{proof}  This follows from the isomorphism $\abmw  {n, S} \cong \akt  {n, 
S} $ and  Proposition
\ref{proposition- bimodule 2}.
\end{proof}

\afil{
\ladres{\adr1\email1}\padres{\adr2\email2}
}{20 November 2004}{2 November 2005}{587}

\end{document}